\documentclass{article}
\usepackage{hyperref}
\usepackage{amsthm,amsmath}
\usepackage{amssymb}
\usepackage{bm}
\usepackage{bbm}
\usepackage{url}
\usepackage{comment}
\usepackage{multirow}
\usepackage{graphicx}
\usepackage{color}
\usepackage{tikz}
\usepackage{caption}
\usepackage{float}
\usepackage{mathtools}
\usepackage{pifont}

\usepackage{amsmath, enumerate, algorithm,algorithmic,caption,subcaption,graphicx,color}
\usepackage{fullpage}

\newcommand{\half}{ \mbox{\small$\frac{1}{2}$}}
\newcommand{\four}{ \mbox{\small$\frac{1}{4}$}}

\newcommand{\be}{\begin{eqnarray}}
\newcommand{\ee}[1]{\label{#1}\end{eqnarray}}
\newcommand{\nn}{\nonumber \\}
\newcommand{\ese}{\end{eqnarray*}}
\newcommand{\bse}{\begin{eqnarray*}}
\newcommand{\rf}[1]{~(\ref{#1})}
\def\qed{\hfill$\square$}

\def\cX{{\cal X}}
\def\cS{{\cal S}}

\def\reg{{\hbox{\scriptsize\rm con}}}
\def\lasso{{\hbox{\scriptsize\rm pen}}}

\def\Tr{\mathop{\hbox{\rm Tr}}}

\def\E{{\mathbf{E}}}

\def\N{{\mathcal{N}}}
\def\CN{{\C\mathcal{N}}}
\def\S{{\mathcal{S}}}

\def\C{\mathbb{C}}
\def\R{\mathbb{R}}
\def\Z{\mathbb{Z}}
\def\Prob{{\hbox{\rm Prob}}}
\def\e{{\rm e}}
\newcommand{\wh}[1]{{\widehat{#1}}}

\def\Argmin{\mathop{\hbox{\rm Argmin}}}

\newcommand{\ci}[2]{{\color{blue}#2}}

\newcommand{\prob}{\text{Prob}}

\newtheorem{theorem}{Theorem}[section]
\newtheorem{lemma}{Lemma}[section]
\newtheorem{corollary}{Corollary}[section]

\allowdisplaybreaks

\title{Structure-Blind Signal Recovery}

\author{
\begin{minipage}[t]{0.25\textwidth}
\center
Dmitry Ostrovsky\thanks{LJK, Universit\'e Grenoble Alpes, 700 Avenue Centrale, 38401 Domaine Universitaire de Saint-Martin-d'H\`eres, France. Email: {\tt \{dmitrii.ostrovskii, anatoli.juditsky\}@imag.fr}}\\
\end{minipage}
\begin{minipage}[t]{0.24 \textwidth}
\center
Zaid Harchaoui\thanks{University of Washington, Seattle, WA 98195, USA. Email: {\tt zaid@uw.edu}}
\end{minipage}
\begin{minipage}[t]{0.24\textwidth}
\center
Anatoli Juditsky\footnotemark[1]
\end{minipage}
\begin{minipage}[t]{0.25\textwidth}
\center
Arkadi Nemirovski\thanks{Georgia Institute of Technology, Atlanta, GA 30332, USA. Email: {\tt nemirovs@isye.gatech.edu}}\\
\end{minipage}
}

\begin{document}

\maketitle

\begin{abstract}
We consider the problem of recovering a signal observed in Gaussian noise.
If the set of signals is convex and compact, and can be specified beforehand, one can use classical linear estimators that achieve a risk within a constant factor of the minimax risk. However, when the set is unspecified, designing an estimator
that is blind to the hidden structure of the signal remains a challenging problem.
We propose a new family of estimators to recover signals observed in Gaussian noise.
Instead of specifying the set where the signal lives, we assume the existence of a well-performing linear estimator.
Proposed estimators enjoy exact oracle inequalities
and can be efficiently computed through convex optimization. We present several numerical illustrations that show
the potential of the approach.
\end{abstract}

\section{Introduction}
We consider the problem of recovering a {\em complex-valued signal} $(x_t)_{t\in \Z}$ from the noisy observations
\begin{equation}\label{prob_model}
y_\tau = x_\tau + \sigma \zeta_\tau, \quad -n \le \tau \le n.
\end{equation}
Here $n\in \Z_+$, and $\zeta_\tau \sim \CN(0,1)$ are i.i.d. standard complex-valued Gaussian random variables, meaning that $\zeta_0= \xi_0^1 + \imath \xi_0^2$ with i.i.d. $\xi_0^1, \xi_0^2 \sim \N(0, 1)$. Our goal is to recover $x_t$, $0 \leq t \leq n$, given the sequence of observations $y_{t-n},..., y_t$ up to instant $t$, a task usually referred to as (pointwise)~\emph{filtering} in machine learning, statistics, and signal processing~\cite{haykin1991adaptive}.

The traditional approach to this problem considers \emph{linear estimators}, or linear filters, which write as
\[
\wh{x}_t=\sum_{\tau=0}^{{n}}\phi_\tau y_{t-\tau}, \quad {0 \leq t \leq n}.
\]
Linear estimators have been thoroughly studied in various forms,
they are both theoretically attractive~\cite{ibragimov1988, donoho1992, donoho1994, tsybakov_mono, wasserman2006all, kailath2000linear, mallat2008wavelet} and easy to use in practice. If the set $\mathcal{X}$ of signals is well-specified, one can usually compute a (nearly) minimax on $\mathcal{X}$ linear estimator in a closed form.
In particular, if $\mathcal{X}$ is a class of smooth signals, such as a H\"older or a Sobolev ball, then the corresponding estimator is given by the kernel estimator with the properly set bandwidth parameter~\cite{tsybakov_mono} and is minimax among all possible estimators. Moreover, as shown by~\cite{ibragimov1984,donoho1994}, if only $\cX$ is convex, compact, and centrally symmetric, the risk of the best linear estimator of $x_t$ is within a small constant factor of the minimax risk over $\cX$. Besides, if the set $\cX$ can be specified in a computationally tractable way, which clearly is still a weaker assumption than classical smoothness assumptions, the best linear estimator can be efficiently computed by solving a convex optimization problem on $\mathcal{X}$. In other words, given a computationally tractable set $\mathcal{X}$ on the input, one can \textit{compute} a nearly-minimax linear estimator and the corresponding (nearly-minimax) risk  over $\mathcal{X}$. The strength of this approach, however, comes at a price:  the set $\cX$ still must be known. Therefore, when one faces a {recovery problem} \textit{without any prior knowledge of $\mathcal{X}$}, this approach cannot be implemented.

We adopt here a novel approach {to filtering, which we refer to as \emph{structure-blind recovery}}. While we do not require $\cX$ to be specified beforehand, we assume that
there exists a \textit{linear oracle} -- a well-performing linear estimator of $x_t$. Previous works~\cite{jn1-2009, jn2-2010, harchaoui2015adaptive}, following a similar philosophy, proved
that one can efficiently adapt to the linear oracle filter of length $m=O(n)$ if the corresponding filter $\phi$ is \textit{time-invariant}, i.e. it  recovers the target signal uniformly well in the $O(n)$-sized neighbourhood of $t$, and if its $\ell_2$-norm is small -- bounded by $\rho/\sqrt{m}$ for a moderate $\rho\geq 1$. The adaptive estimator is computed by minimizing the $\ell_\infty$-norm of the filter discrepancy, in the Fourier domain,
under the constraint on the $\ell_1$-norm of the filter in the Fourier domain. Comparing to the oracle linear filter, the price for adaptation is proved to be $O(\rho^3\sqrt{\ln n})$, with the lower bound of $O(\rho\sqrt{\ln n})$~\cite{jn1-2009,harchaoui2015adaptive}.

We make the following contributions:
\begin{itemize}
\item we propose  a new family of recovery methods, obtained by solving a least-squares problem constrained or penalized by the $\ell_1$-norm of the filter in the Fourier domain;
\item we prove exact oracle inequalities for the $\ell_2$-risk of these methods;
\item we show that the price for adaptation improves upon previous works~\cite{jn1-2009,harchaoui2015adaptive} to $O(\rho^2 \sqrt{\ln n})$ for the point-wise risk and to $O(\rho\sqrt{\ln n})$ for the $\ell_2$-risk.
\item we present numerical experiments that show the potential of the approach on synthetic and real-world images and signals.
\end{itemize}

Before presenting the theoretical results, let us introduce the notation we use throughout the paper.
\paragraph{Filters}
Let $\C(\Z)$ be the linear space of all two-sided complex-valued sequences $x = \left\{x_t \in \mathbb{C}\right\}_{t \in \mathbb{Z}}$. For $k, k' \in \Z$ we consider finite-dimensional subspaces $$\C(\Z_{k}^{k'}) = \left\{x \in \C(\Z): \quad x_t = 0, \quad t \notin [k,k'] \right\}.$$ It is convenient to  identify $m$-dimensional complex vectors, $m=k'-k+1$, with elements of $\C(\Z_{k}^{k'})$ by means of the notation:
$$ x_{k}^{k'}:=[x_k;\, ...;\, x_{k'}]\in \C^{k'-k+1}.$$
We associate to linear mappings $\C(\Z_{k}^{k'}) \to \C(\Z_{j}^{j'})$ $(j'-j+1)\times (k'-k+1)$ matrices with complex entries.
The \textit{convolution} $u * v$ of two sequences $u,v \in \C(\Z)$ is a sequence with elements
\[
[u*v]_t = \sum_{\tau \in \mathbb{Z}} u_\tau v_{t-\tau}, \quad t \in \Z.
\]
Given observations \eqref{prob_model} and $\varphi\in \mathbb{C}(\Z_{0}^{m})$ consider the {\em (left) linear estimation} of  $x$ associated with {\em filter} $\varphi$:
\[
\wh{x}_t = [\varphi*y]_t
\]
($\wh{x}_t$ is merely a kernel estimate of $x_t$ by a kernel $\varphi$ supported on $[0,...,m]$). 
\paragraph{Discrete Fourier transform}
We define the unitary \emph{Discrete Fourier transform} (DFT) operator $F_{n}:\, \C^{n+1} \to \C^{n+1}$ by
\begin{equation*}
{z\mapsto F_{n} z},\quad [F_{n} z]_k =
{(n+1)^{-{1/2}}} \sum\limits_{t = 0}^n z_{t}\, e^{\frac{2\pi\imath k t}{n+1}}, \quad 0 \le k \le n.
\end{equation*}
The \emph{inverse Discrete Fourier transform} (iDFT) operator $F^{-1}_n$ is given by $F^{-1}_n:=F^H_n$ (here $A^H$ stands for Hermitian adjoint of $A$).
By the Fourier inversion theorem, $F^{-1}_{n}\, (F^{\vphantom{-1}}_{n}\, z) = z$.
\par
We denote $\|\cdot\|_p$ usual $\ell_p$-norms on $\C(\Z)$:
$\|x\|_{p} = (\sum_{t \in \Z} \left| x_t \right|^p)^{1/p}$, $p \in [1, \infty]$. Usually, the argument will be finite-dimensional -- an element of $\C(\Z_k^{k'})$; we reserve the special notation $$\|x\|_{n,p} := \|x_0^n\|_p.$$
Furthermore, DFT allows to equip $\C(\Z_{0}^n)$  with the norms associated with $\ell_p$-norms in the spectral domain:
\[
\|x\|^*_{n,p} := \|x_0^n\|^*_{p} :=  \|F_{n}x_0^n\|_p, \quad p \in [1, \infty];
\]
note that unitarity of the DFT implies the Parseval identity: $\|x\|^{\vphantom*}_{n,2} = \|x\|^*_{n,2}$.

Finally, $c$, $C$, and $C'$ stand for generic absolute constants.
\section{Oracle inequality for constrained recovery}\label{sec:oreg}
Given observations \eqref{prob_model} and $\overline{\varrho}>0$, we first consider the {\em constrained recovery} $\wh{x}_\reg$ given by
\[%begin{equation}
[\wh{x}_\reg]_t=[\wh{\varphi}*y]_t, \quad t=0,...,n,
\]%end{equation}
where $\wh\varphi$ is an optimal solution of the constrained optimization problem
\begin{equation}\label{eq:l2-1}
\min_{\varphi\in \C(\Z^{n}_0)}\big\{\|y-\varphi*y\|_{n,2}:\;
\|\varphi\|_{n,1}^*\leq \overline{\varrho}/\sqrt{n+1}\big\}.
\end{equation}

The constrained recovery estimator minimizes a least-squares fit criterion under a constraint on $\|\varphi\|_{n,1}^*=\|F_{n} \varphi_0^n \|_{1}$,
that is an $\ell_1$ constraint on the discrete Fourier transform of the filter. While the least-squares objective naturally follows from the Gaussian noise assumption,
the constraint can be motivated as follows.

\paragraph{Small-error linear filters}
Linear filter $\varphi^o$ with a small $\ell_1$ norm in the spectral domain and small recovery error exists, essentially, whenever there exists a linear filter with small recovery error~\cite{jn1-2009,harchaoui2015adaptive}. Indeed, let us say that $x\in \C(\Z_0^n)$ is {\em simple}~\cite{harchaoui2015adaptive} {\em with parameters $m\in \Z_+$ and $\rho\geq 1$}
if there exists $\phi^o\in \C(\Z_0^m)$ such that for all $-m\leq \tau\leq 2m$,
\begin{equation}\label{eq:ss1}
\left[\E\left\{|x_\tau - [\phi^{o} * y]_\tau|^2\right\}\right]^{1/2} \le {\frac{\sigma\rho}{\sqrt{m+1}}}.
\end{equation}
  In other words, $x$ is $(m,\rho)$-simple if there exists a hypothetical filter $\phi^o$ of the length at most $m+1$ which recovers $x_\tau$ with squared risk uniformly bounded by $\frac{\sigma^2\rho^2}{m+1}$  in the interval $-m\leq \tau\leq  2m$.
  Note that \eqref{eq:ss1} clearly implies that $\|\phi^o\|_2\leq\rho/\sqrt{m+1}$, and that
$|[x - \phi^{o} * x]_\tau| \leq \sigma \rho / \sqrt{m+1}$ $\forall\tau,\,-m\leq \tau\leq 2m$. Now, let $n=2m$, and let \[
\varphi^o=\phi^o*\phi^o\in \C^{n+1}.
\]
As proved in Appendix~\ref{sec:prdisc}, we have
\begin{equation}\label{eq:l1f}
\|\varphi^o\|_{n,1}^*\leq 2\rho^2/\sqrt{n+1},
\end{equation}
and, for a moderate absolute constant $c$,
\begin{equation}\label{eq:var2f}
\left\|x - \varphi^{o} * y\right\|_{n,2} \leq c\sigma\rho^2\sqrt{1+\ln[1/\alpha]}
\end{equation}
with probability $1-\alpha$.
To summarize, {\em if $x$ is $(m,\rho)$-simple, i.e., when there exists a filter $\phi^o$ of length $\leq m+1$ which recovers $x$ with small risk on the interval $[-m, 2m]$, then the filter $\varphi^o=\phi^o*\phi^o$ of the length at most $n+1$, with $n=2m$, has small norm $\|\varphi^o\|_{n,1}^*$ and recovers the signal $x$ with (essentially the same) small risk on the interval $[0,n]$.}
\paragraph{Hidden structure}
The constrained recovery estimator is completely blind to a possible hidden structure of the signal, yet can seamlessly adapt to it when such a structure exists, in a way that we can rigorously establish. Using the right-shift operator on $\C(\Z)$, $[\Delta x ]_t = x_{t-1}$, we formalize the hidden structure  as an unknown shift-invariant linear subspace of $\C(\Z)$, $\Delta \mathcal{S} = \mathcal{S}$, of a small dimension $s$. We do not assume that $x$ belongs to that subspace. Instead, we make a more general assumption that $x$ is \emph{close} to this subspace,
that is, it may be decomposed into a sum of a component that lies in the subspace and a component whose norm we can control.
\paragraph{Assumption A}
{\em We suppose that $x$ admits the decomposition
\[
x=x^{\mathcal{S}}+\varepsilon, \quad x^{\mathcal{S} }\in \mathcal{S},
\]
where $\mathcal{S}$ is an (unknown) shift-invariant, $\Delta \mathcal{S} = \mathcal{S}$, subspace of $\C(\Z)$ of dimension $s$, $1 \le s\le n+1$, and $\varepsilon\,$ is ``small'', namely,
\[
\left\|\Delta^{\tau} \varepsilon\right\|_{n,2} \leq \sigma\varkappa, \quad 0 \le \tau \le n.
\]}
\par\noindent Shift-invariant subspaces of $\C(\Z)$ are exactly the sets of solutions of homogeneous linear difference equations with polynomial operators. This is summarized by the following lemma (we believe it is a known fact; for completeness we provide a proof in Appendix~\ref{sec:prdisc}.
\begin{lemma}\label{lemma:shift_invariant}
Solution set of a homogeneous difference equation with a polynomial operator $p(\Delta)$,
\begin{equation}\label{eq:diff_eq}
[p(\Delta) x]_t = \left[\sum_{\tau=0}^s p_\tau x_{t-\tau}\right] = 0, \quad t \in \Z,
\end{equation}
with $\deg(p(\cdot))=s$, $p(0)=1$, is a shift-invariant subspace of $\C(\Z)$ of dimension $s$. Conversely, any shift-invariant subspace $\S \subset \C(\Z)$, $\Delta \S \subseteq \S$, $\dim(\S)=s < \infty$, is the set of solutions of some homogeneous difference equation \eqref{eq:diff_eq} with $\deg(p(\cdot))=s$, $p(0)=1$. Moreover, such $p(\cdot)$ is unique.
\end{lemma}

\noindent On the other hand, for any polynomial $p(\cdot)$, solutions of \eqref{eq:diff_eq} are exponential polynomials~\cite{jn1-2009} with frequencies determined by the roots of $p(\cdot)$. For instance, discrete-time polynomials $x_t=\sum_{k=0}^{s-1} c_k t^k$, $t\in \Z$ of degree $s-1$ (that is, exponential polynomials with all zero frequencies) form a linear space of dimension $s$ of solutions of the equation \eqref{eq:diff_eq} with a polynomial $p(\Delta) = (1-\Delta)^s$ with a unique root of multiplicity $s$, having coefficients $p_k=(-1)^k {s \choose k}$. Naturally, signals which are close, in the $\ell_2$ distance, to discrete-time polynomials are Sobolev-smooth functions sampled over the regular grid \cite{jn2-2010}.
Sum of harmonic oscillations $x_t = \sum_{k=1}^s c_k \e^{\imath \omega_k t}$, $\omega_k \in [0,2\pi)$ being all different, is another example; here, $p(\Delta) = \prod_{k=1}^s (1-e^{\imath \omega_k}\Delta)$.\\

\par\noindent We can now state an oracle inequality for the constrained recovery estimator; see Appendix~\ref{sec:proracle} for the proof.
\begin{theorem}\label{th:l2con}
Let $\overline{\varrho}\geq 1$, and let $\varphi^o\in\C(\Z^{n}_0)$ be such that
\[
\|\varphi^o\|^*_{n,1}\leq {\overline{\varrho}/\sqrt{n+1}}.
\]
Suppose that Assumption A holds for some $s\in \Z_+$ and $\varkappa<\infty$. Then for any {$\alpha$, $0<\alpha\leq 1,$} it holds with probability at least $1-\alpha$:
\begin{equation}\label{eq:th1}
\|x - \wh{x}_{\reg} \|_{n,2} \le \|x-\varphi^o*y\|_{n,2}
+ C\sigma
{\sqrt{s + \overline\varrho\,\big(\varkappa\sqrt{\ln\left[{1/ \alpha}\right]}+\ln\left[{n/ \alpha}\right]\big)}}.
\end{equation}
\end{theorem}
When considering \emph{simple} signals, Theorem~\ref{th:l2con} gives the following.
\begin{corollary}\label{cor:l2con}
Assume that signal $x$ is $(m,\rho)$-simple, $\rho\geq 1$ and $m\in \Z_+$. Let $n=2m$, $\overline{\varrho}\geq 2\rho^2$, and let Assumption A hold for some $s\in \Z_+$ and $\varkappa<\infty$.
Then for any $\alpha$, $0<\alpha\leq 1$, it holds with probability at least $1-\alpha$:
\[
\|x - \wh{x}_{\reg} \|_{n,2} \le C\sigma\rho^2\sqrt{\ln[1/\alpha]}+ C'\sigma \sqrt{s + \overline\varrho \, \big(\varkappa\sqrt{\ln\left[{1/ \alpha}\right]}+\ln\left[{n/ \alpha}\right]\big)}.
\]
\end{corollary}
\paragraph{Adaptation and price}
The price for adaptation in Theorem~\ref{th:l2con} and Corollary~\ref{cor:l2con} is determined by three parameters: the bound on the filter norm $\overline{\varrho}$, the deterministic error $\varkappa$, and the subspace dimension $s$.
Assuming that the signal to recover is simple, and that $\overline{\varrho} = 2\rho^2$, let us compare the magnitude of the oracle error to the term of the risk which reflects ``price of adaptation''. Typically (in fact, in all known to us cases of recovery of signals from a shift-invariant subspace), the parameter $\rho$ is at least $\sqrt{s}$. Therefore, the bound~\rf{eq:var2f} implies the ``typical bound'' $O(\sigma\sqrt{\gamma}\rho^2)=\sigma s\sqrt{\gamma}$ for the term $\| x - \varphi^o * y\|_{n,2}$ (we denote $\gamma = \ln(1/\alpha)$).
As a result, for instance, in the ``parametric situation'', when the signal belongs or is very close to the subspace, that is when $\varkappa = O(\ln(n))$, the price of adaptation $O\left(\sigma[s + \rho^2(\gamma + \sqrt{\gamma} \ln n)]^{1/2}\right)$ is much smaller than the bound on the oracle error. {In the ``nonparametric situation'', when $\varkappa = O(\rho^2)$, the price of adaptation has the same order of magnitude as the oracle error.}

Finally, note that under the premise of Corollary~\ref{cor:l2con} we can also bound the pointwise error. We state the result for $\overline{\varrho}= 2\rho^2$ for simplicity; the proof is in Appendix~\ref{sec:proracle}.
\begin{theorem}\label{th:point}
Assume that signal $x$ is $(m,\rho)$-simple, $\rho\geq 1$ and $m\in \Z_+$. Let $n=2m$, $\overline{\varrho}= 2\rho^2$, and let Assumption A hold for some $s\in \Z_+$ and $\varkappa<\infty$.
Then for any $\alpha, 0<\alpha\leq 1$, the constrained recovery $\wh{x}_\reg$ satisfies
\[
|x_n - [\wh{x}_{\reg}]_n| \le C{\frac{\sigma\rho}{\sqrt{m+1}}}
\left[\rho^2\sqrt{\ln[n/\alpha]}+ \rho\sqrt{ \varkappa\sqrt{\ln\left[{1/ \alpha}\right]}}+
\sqrt{s}\right].
\]
\end{theorem}

\section{Oracle inequality for penalized recovery}\label{sec:open}
%The constrained recovery estimator inherently assumes that the parameter $\varrho = 2\rho^2$ in th is known. If the noise variance is known (or can be estimated from data), then we can build a more practical estimator that still enjoys an oracle inequality.
To use the constrained recovery estimator with a provable guarantee, see e.g. Theorem~\ref{th:l2con}, one must know the norm of a small-error linear filter $\varrho$, or at least have an upper bound on it. However, if this parameter is unknown, but instead the noise variance is known (or can be estimated from data), we can build a more practical estimator that still enjoys an oracle inequality.

The {\em penalized recovery} estimator $[\wh{x}_\lasso]_t=[\widehat{\varphi}*y]_t$ is an optimal solution to a regularized least-squares minimization problem, where the regularization penalizes the $\ell_1$-norm of the filter in the Fourier domain:
\begin{equation}
\label{optim2pen}
\widehat{\varphi}\in \Argmin\limits_{\varphi \in \C(\Z_0^n) }
\left\{\|y-\varphi*y\|^2_{n,2} + \lambda \sqrt{n+1}\, \|\varphi\|^*_{n,1}\right\}.
\end{equation}
Similarly to Theorem \ref{th:l2con}, we establish an oracle inequality for the penalized recovery estimator.
\begin{theorem}\label{th:l2pen}
Let Assumption A hold for some $s\in \Z_+$ and $\varkappa<\infty$,  and let $\varphi^o \in \C(\Z_0^n)$ satisfy $\|\varphi^o \|_{n,1}^*\leq \varrho/\sqrt{n+1}$ for some $\varrho \ge 1$.
\paragraph{$1^o$.} Suppose that the regularization parameter of penalized recovery $\wh{x}_\lasso$ satisfies $\lambda\geq \underline{\lambda},$
\[
\underline{\lambda} := 60\sigma^2\ln[63n/\alpha].
\]
Then, for $0 < \alpha \le 1$, it holds with probability at least $1-\alpha$:
\[
\|x - \widehat{x}_\lasso\|_{n,2} \le \|x-\varphi^{o}*y\|_{n,2}
+ C\sqrt{\varrho\lambda} + C'\sigma \sqrt{s + (\widehat\varrho + 1) \varkappa\sqrt{\ln[1/\alpha]}},
\]
where $\widehat\varrho := \sqrt{n+1}\, \|\widehat\varphi\|^*_{n,1}$.
\paragraph{$2^o$.} Moreover, if $\varkappa\leq \bar\varkappa,$
\[\bar\varkappa := \frac{10\ln[42 n/\alpha]}{\sqrt{\ln\left[16/\alpha\right]}},
\] and $\lambda \ge 2\underline{\lambda}$, one has
\[
\|x - \widehat{x}_\lasso\|_{n,2} \le \|x-\varphi^{o}*y\|_{n,2}
+ C\sqrt{\varrho\lambda} +C'\sigma \sqrt{s}.
\]
\end{theorem}
The proof closely follows that of Theorem \ref{th:l2con} ans is also provided in Appendix~\ref{sec:proracle}.

\section{Discussion}\label{sec:disc}
\def\X{{\cal X}}
There is some redundancy between ``simplicity'' of a signal, as defined by \eqref{eq:ss1}, and Assumption A. Usually a simple signal or image $x$ is also close to a low-dimensional subspace of $\C(\Z)$ (see, e.g., \cite[section 4]{jn2-2010}), so that Assumption A holds ``automatically''. Likewise,
$x$ is ``almost'' simple when it is close to a low-dimensional time-invariant subspace. Indeed, if $x\in \C(\Z)$ {\em belongs} to $\S$, i.e. Assumption A holds with $\varkappa= 0$, one can easily verify that for $n\geq s$ there exists a filter $\phi^o\in \C(\Z_{-n}^n)$ such that
\begin{equation}\label{eq:interps}
\|\phi^o\|_2\leq \sqrt{s/(n+1)}, \;\;\mbox{and}\;\;
x_\tau=[\phi^o*x]_\tau, \;\; \tau \in \Z \; ,
\end{equation}
see Appendix~\ref{sec:prdisc} for the proof. This implies that $x$ can be recovered efficiently from observations \eqref{prob_model}:
\[
\left[\E\big\{|x_\tau-[\phi^o*y]_\tau|^2\big\}\right]^{1/2}\leq \sigma\sqrt{\frac{s}{n+1}}.
\]
In other words, if instead of the filtering problem we were interested in the {\em interpolation} problem of recovering $x_t$ given $2n+1$ observations $y_{t-n},...,y_{t+n}$ on the left {\em and} on the right of $t$, Assumption A would imply a kind of simplicity of $x$.
On the other hand, it is clear that Assumption A is not sufficient to imply the simplicity of $x$  ``with respect to the filtering'', in the sense of the definition we use in this paper,
when we are allowed to use only observations on the left of $t$ to compute the estimation of $x_t$. Indeed, one can see, for instance, that already signals from the parametric family $\X_\alpha=\{x\in \C(\Z):\,x_\tau=c\alpha^\tau, \,c\in \C\}$, with a given $|\alpha|>1$, which form a  one-dimensional space of solutions of the equation $x_\tau=\alpha x_{\tau-1}$, cannot be estimated with small risk at $t$ using only observations on the left of $t$ (unless $c=0$), and thus are not simple in the sense of \eqref{eq:ss1}.

Of course, in the above example, the ``difficulty'' of the family $\X_\alpha$  is due to instability of solutions of the difference equation which explode when $\tau\to +\infty$. Note that signals $x\in \X_\alpha$ with $|\alpha|\leq 1$ (linear functions, oscillations, or damped oscillations) are simple.  More generally, suppose that $x$ satisfies a difference equation of degree $s$:
\begin{equation}\label{eq:deq}
0=p(\Delta)x_\tau\left[=\sum_{i=0}^s p_ix_{\tau-i}\right],
\end{equation}
where $p(z)=\sum_{i=0}^s p_iz^i$ is the corresponding characteristic polynomial and $\Delta$ is the right shift operator. When $p(z)$ is unstable -- has roots {\em inside} the unit circle -- (depending on ``initial conditions'') the set of solutions to the equation \eqref{eq:deq} contains  difficult to filter signals. Observe that stability of solutions is related to the direction of the time axis; when the characteristic polynomial $p(z)$ has roots {\em outside} the unit circle, the corresponding solutions may be ``left unstable'' -- increase exponentially when $\tau \to -\infty$. In this case ``right filtering'' -- estimating $x_\tau$ using observations on the right of $\tau$ -- will be difficult. A special situation where both interpolation and filtering are always simple arises when the characteristic polynomial of the difference equation has all its roots on the unit circle. In this case, solutions to \eqref{eq:deq} are ``generalized harmonic oscillations'' (harmonic oscillations modulated by polynomials), and such signals are known to be simple. Theorem~\ref{th:sines} summarizes the properties of the solutions of (\ref{eq:deq}) in this particular case; see Appendix~\ref{sec:prdisc} for the proof.
\begin{theorem}\label{th:sines} Let $s$ be a positive integer, and let $p=[p_0;...;p_s]\in \C^{s+1}$ be such that
the polynomial $p(z)=\sum_{i=0}^s p_i z^i$ has all its roots on the unit circle. Then for every integer $m$ satisfying
\[
m\geq m(s):=C s^2\ln(s+1),
\]
one can point out $q\in \C^{m+1}$
such that any solution to \eqref{eq:deq} satisfies
\[
x_\tau=[q*x]_\tau,\;\;\forall \tau\in \Z,
\]
and
\begin{equation}\label{propofq}
\|q\|_2\leq {\rho(s,m)}/{\sqrt{m}} \quad \mbox{where}\quad
\rho(s,m)= C'\min\left\{s^{3/2}\sqrt{\ln s},\,s\sqrt{\ln[ms]}\right\}.
\end{equation}
\end{theorem}
\section{Numerical experiments}
\label{sec:exps}
We present preliminary results on simulated data of the proposed adaptive signal recovery methods in several application scenarios.
We compare the performance of the {\em penalized $\ell_2$-recovery} of Sec.~\ref{sec:open} to that of the {Lasso} recovery of~\cite{recht1} in signal and image denoising problems.
Implementation details for the penalized $\ell_2$-recovery are given in Sec.~\ref{sec:expdetails}. Discussion of the discretization approach underlying the competing Lasso method can be found in \cite[Sec. 3.6]{recht1}. 

We follow the same methodology in both signal and image denoising experiments. For each level of the signal-to-noise ratio $\text{SNR} \in \{ 1,\, 2,\, 4,\, 8,\, 16\}$, 
we perform $N$ Monte-Carlo trials.
In each trial, we generate a random signal $x$ on a regular grid with $n$ points, corrupted by the i.i.d. Gaussian noise of variance $\sigma^2$. The signal 
is normalized: $\|x\|_2=1$ so $\text{SNR}^{-1} = \sigma\sqrt{n}$. 
We set the regularization penalty in each method as follows. For penalized $\ell_2$-recovery~\eqref{optim2pen}, we use $\lambda = 2\sigma^2\log[63n/\alpha]$ with $\alpha = 0.1$. For Lasso~\cite{recht1}, we use the common setting $\lambda = \sigma\sqrt{2\log n}$.
We report experimental results by plotting the $\ell_2$-error $\| \widehat x - x \|_2$, averaged over $N$ Monte-Carlo trials, versus the inverse of the signal-to-noise ratio $\text{SNR}^{-1}$. 

\begin{figure}[t!]
\center
\begin{minipage}{0.24\textwidth}
\centering
\includegraphics[width=1\textwidth, height=0.135\textheight, clip=true, angle=0]{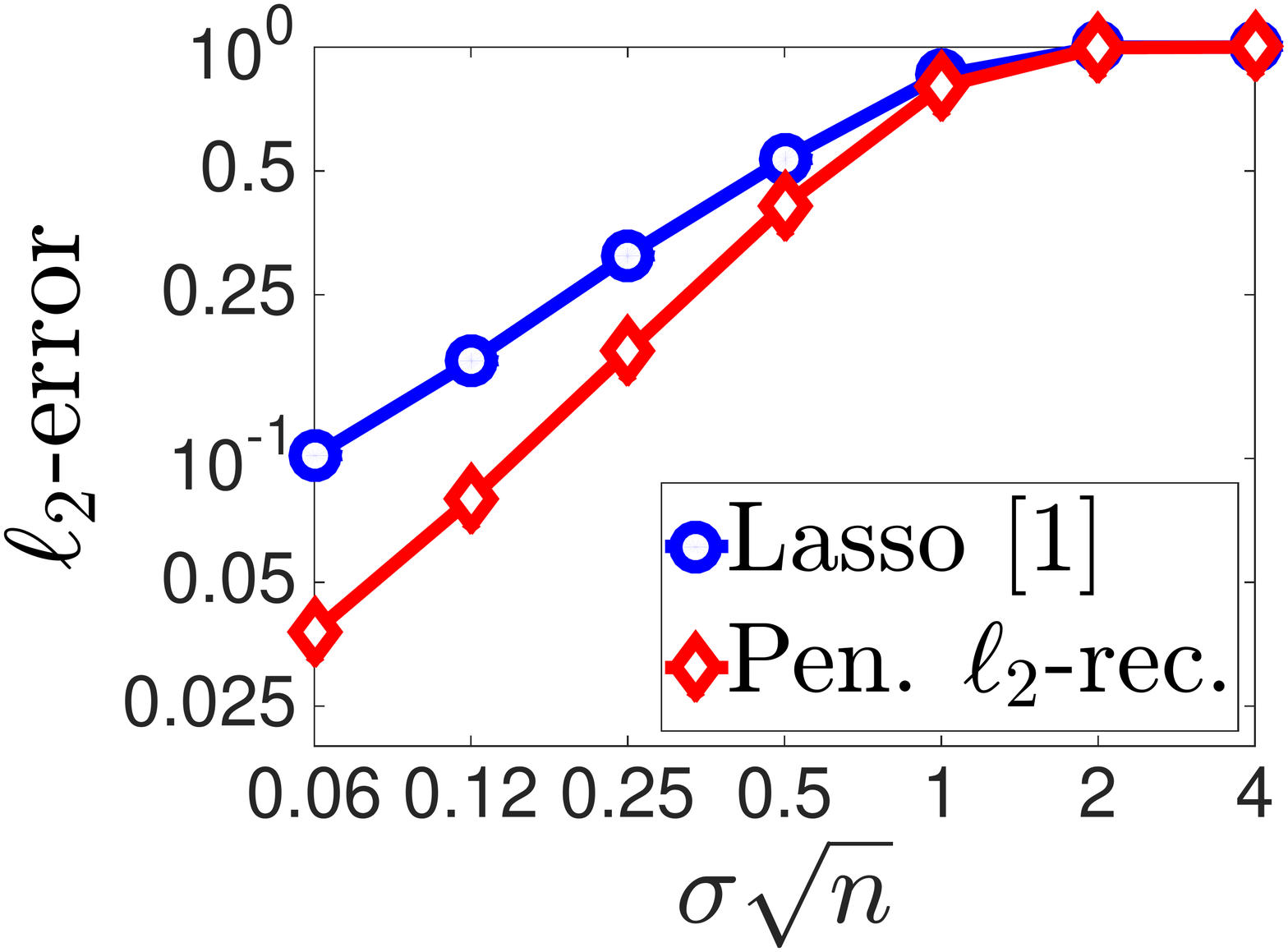}
\end{minipage}
\begin{minipage}{0.24\textwidth}
\centering
\includegraphics[width=1\textwidth, height=0.135\textheight, clip=true, angle=0]{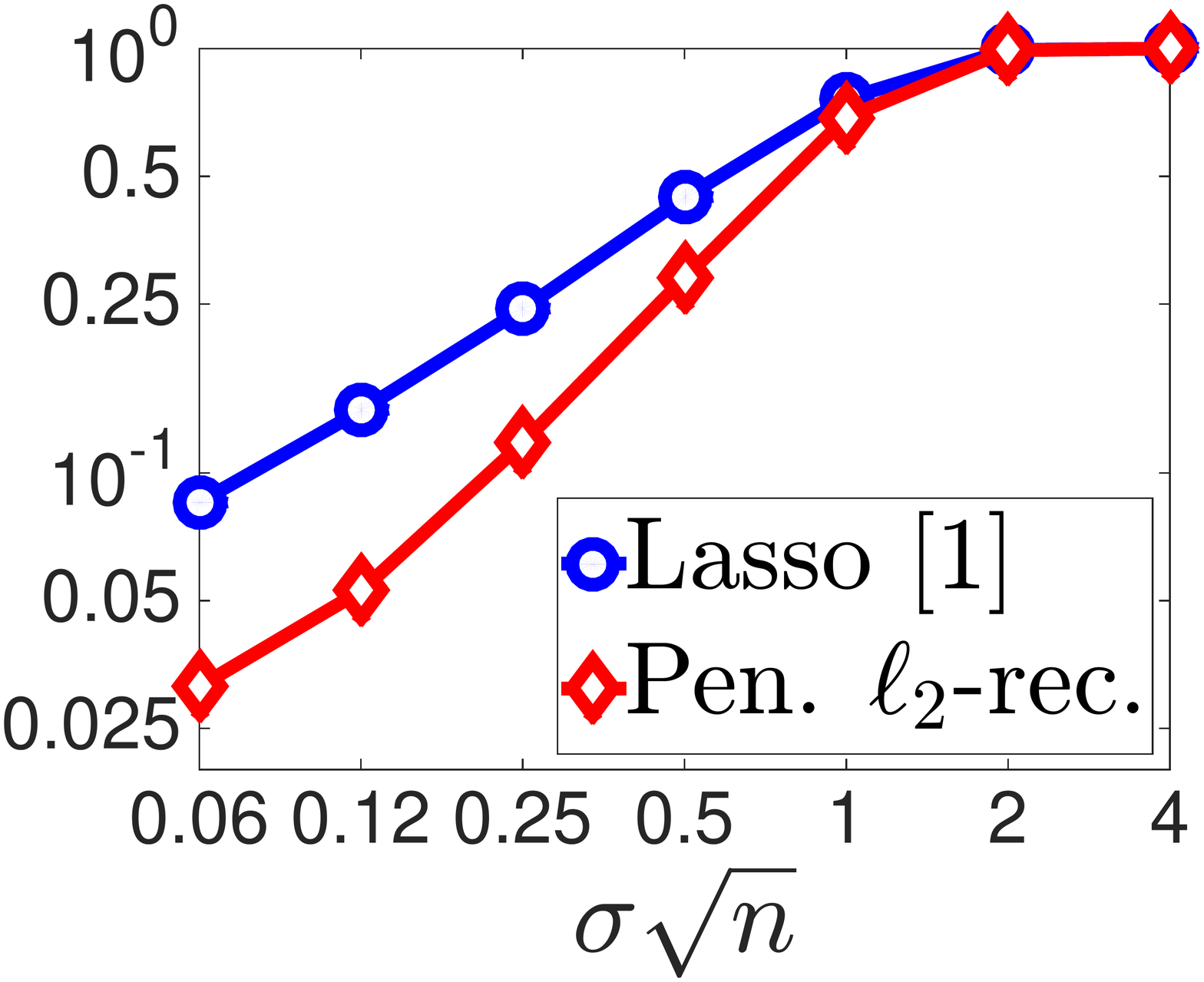}
\end{minipage}
\begin{minipage}{0.24\textwidth}
\centering
\includegraphics[width=1\textwidth, height=0.135\textheight, clip=true, angle=0]{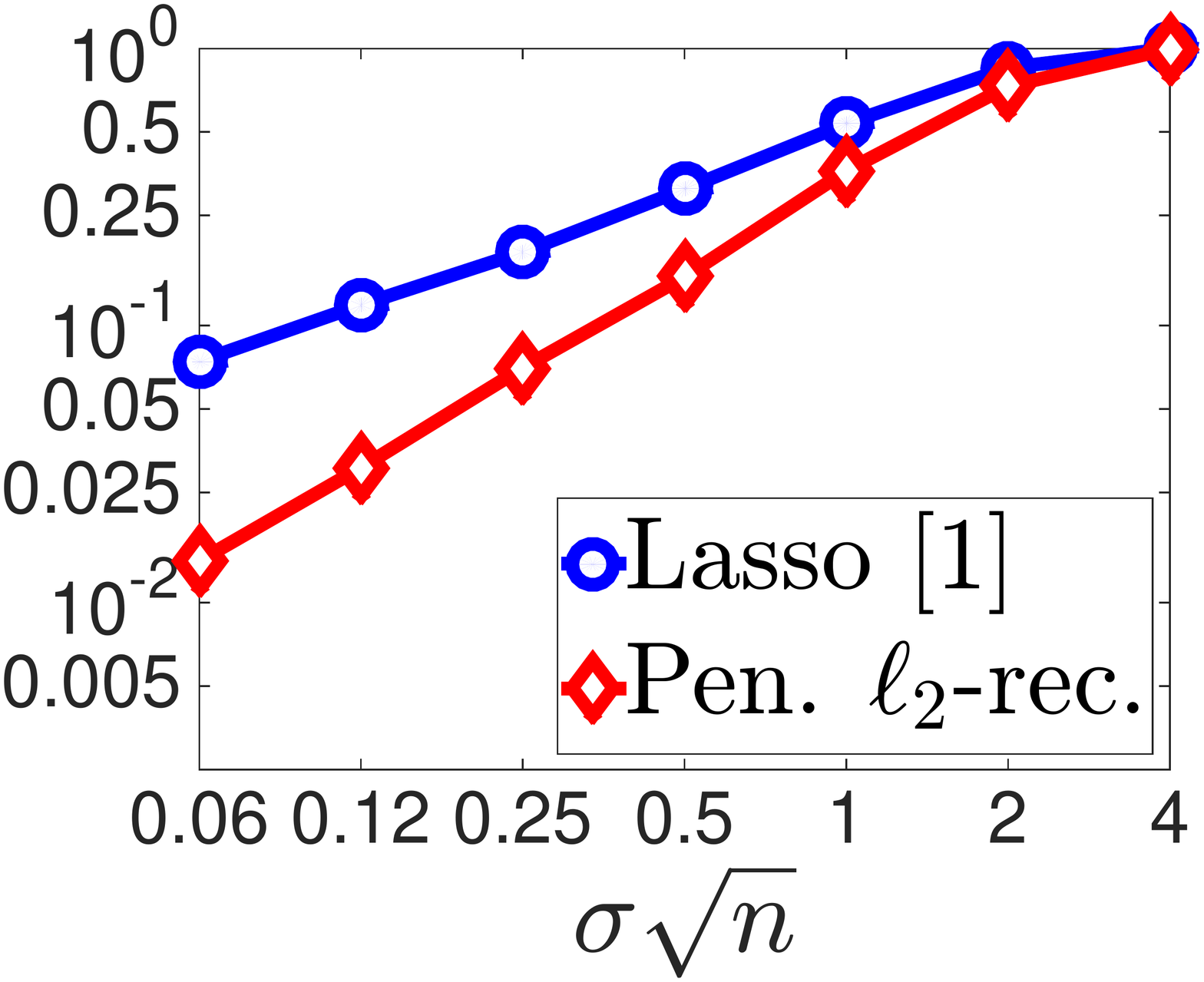}
\end{minipage}
\begin{minipage}{0.24\textwidth}
\includegraphics[width=1\textwidth, height=0.135\textheight, clip=true, angle=0]{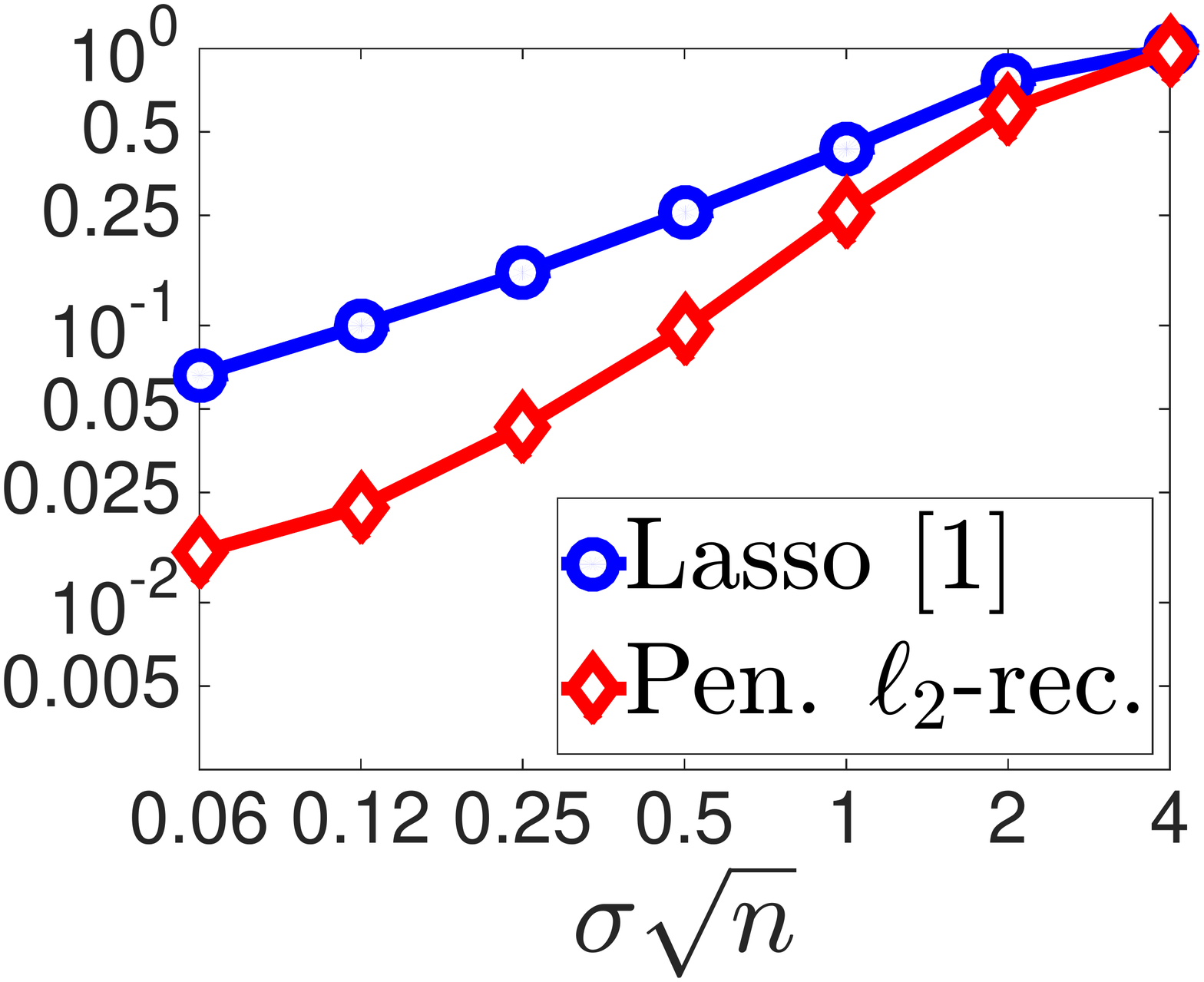}
%Scenario II
\end{minipage}
\caption{Signal and image denoising in different scenarios, left to right: \textit{RandomSpikes}, \textit{CoherentSpikes}, \textit{RandomSpikes-2D}, and \textit{CoherentSpikes-2D}. The steep parts of the curves on high noise levels correspond to observations being thresholded to zero.}
\label{fig:mc-sin}
\end{figure}

\paragraph{Signal denoising}
We consider denoising of a one-dimensional signal in two different scenarios, fixing $N=100$ and $n=100$. In the \textit{RandomSpikes} scenario, the signal is a sum of 4 harmonic oscillations, each characterized by a spike of a random amplitude at a random position in the continuous frequency domain $[0,2\pi]$. In the \textit{CoherentSpikes} scenario, the same number of spikes is sampled by pairs. Spikes in each pair have the same amplitude and are separated by only $0.1$ of the DFT bin $2\pi/n$ which could make recovery harder due to high signal coherency. However, in practice we found \textit{RandomSpikes} to be slightly harder than \textit{CoherentSpikes} for both methods, see Fig.~\ref{fig:mc-sin}. As Fig.~\ref{fig:mc-sin} shows, the proposed penalized $\ell_2$-recovery 
outperforms the Lasso method for all noise levels. The performance gain is particularly significant for high signal-to-noise ratios. 

\paragraph{Image Denoising}
We now consider recovery of an unknown regression function $f$ on the regular grid on $[0,1]^2$ given the noisy observations:
\begin{equation}\label{eq:grid_observations}
y_\tau = x_\tau + \sigma \zeta_\tau, \quad \tau \in \left\{{0}, {1}, ..., m-1\right\}^2,
\end{equation}
where $x_\tau = f(\tau/m)$. We fix $N=40$, and the grid dimension $m=40$; the number of samples is then $n=m^2$. For the penalized $\ell_2$-recovery, we implement the blockwise denoising strategy (see Appendix for the implementation details) with just one block for the entire image. We present additional numerical illustrations in the supplementary material.

We study three different scenarios for generating the ground-truth signal in this experiment. The first two scenarios, \textit{RandomSpikes-2D} and \textit{CoherentSpikes-2D}, are two-dimensional counterparts of those studied in the signal denoising experiment: the ground-truth signal is a sum of $4$ harmonic oscillations in $\R^2$ with random frequencies and amplitudes. The separation in the \textit{CoherentSpikes-2D} scenario is $0.2\pi/m$ in each dimension of the torus $[0,2\pi]^2$. The results for these scenarios are shown in Fig.~\ref{fig:mc-sin}. 
Again, the proposed penalized $\ell_2$-recovery outperforms the Lasso method for all noise levels, especially for high signal-to-noise ratios. 

In scenario \textit{DimensionReduction-2D} we investigate the problem of estimating a function with a hidden low-dimensional structure. We consider the single-index model of the regression function:
\begin{equation}\label{eq:single-index}
f(t) = g(\theta^T t),  \quad g(\cdot) \in \mathcal{S}_\beta^1(1). 
\end{equation}
Here, $\mathcal{S}_\beta^1(1) = \{g: \R \to \R, \|g^{(\beta)}(\cdot)\|_2 \le 1\}$ is the Sobolev ball of smooth periodic functions on $[0,1]$, and the unknown structure is formalized as the direction $\theta$. In our experiments we sample the direction $\theta$ uniformly at random and consider different values of the smoothness index $\beta$.
If it is known a priori that the regression function possesses the structure \eqref{eq:single-index}, and only the index is unknown, one can use estimators attaining "one-dimensional" rates of recovery; see e.g.~\cite{lepski2014adaptive} and references therein. In contrast, our recovery algorithms are not aware of the underlying structure but might still adapt to it. 
 
\begin{figure}[ht!]
\center
\begin{minipage}{0.32\textwidth}
\centering
$\beta = 2$
\includegraphics[width=1\textwidth, height=0.18\textheight, clip=true, angle=0]{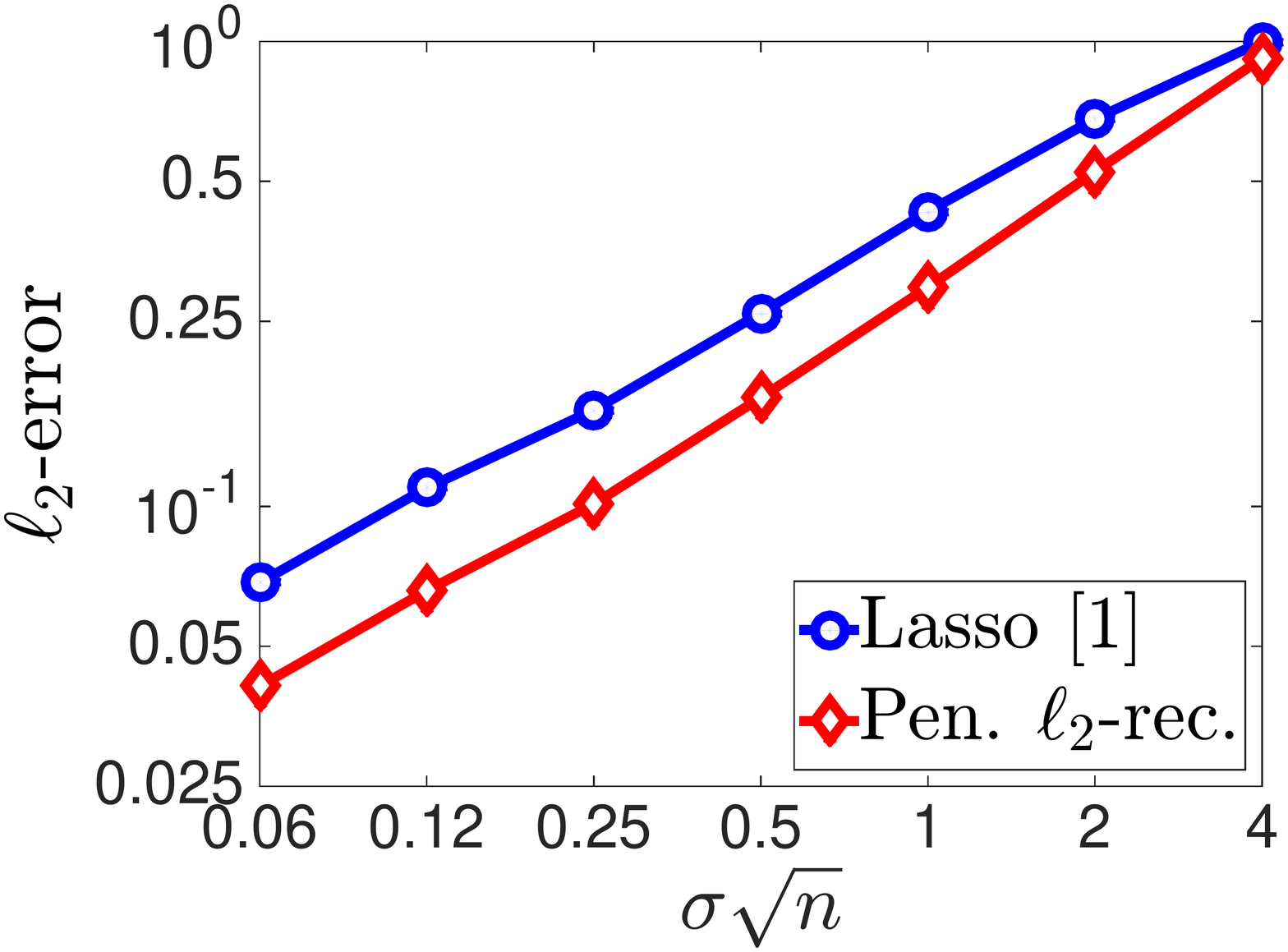}
\end{minipage}
\begin{minipage}{0.32\textwidth}
\centering
$\beta = 1$
\includegraphics[width=1\textwidth, height=0.18\textheight, clip=true, angle=0]{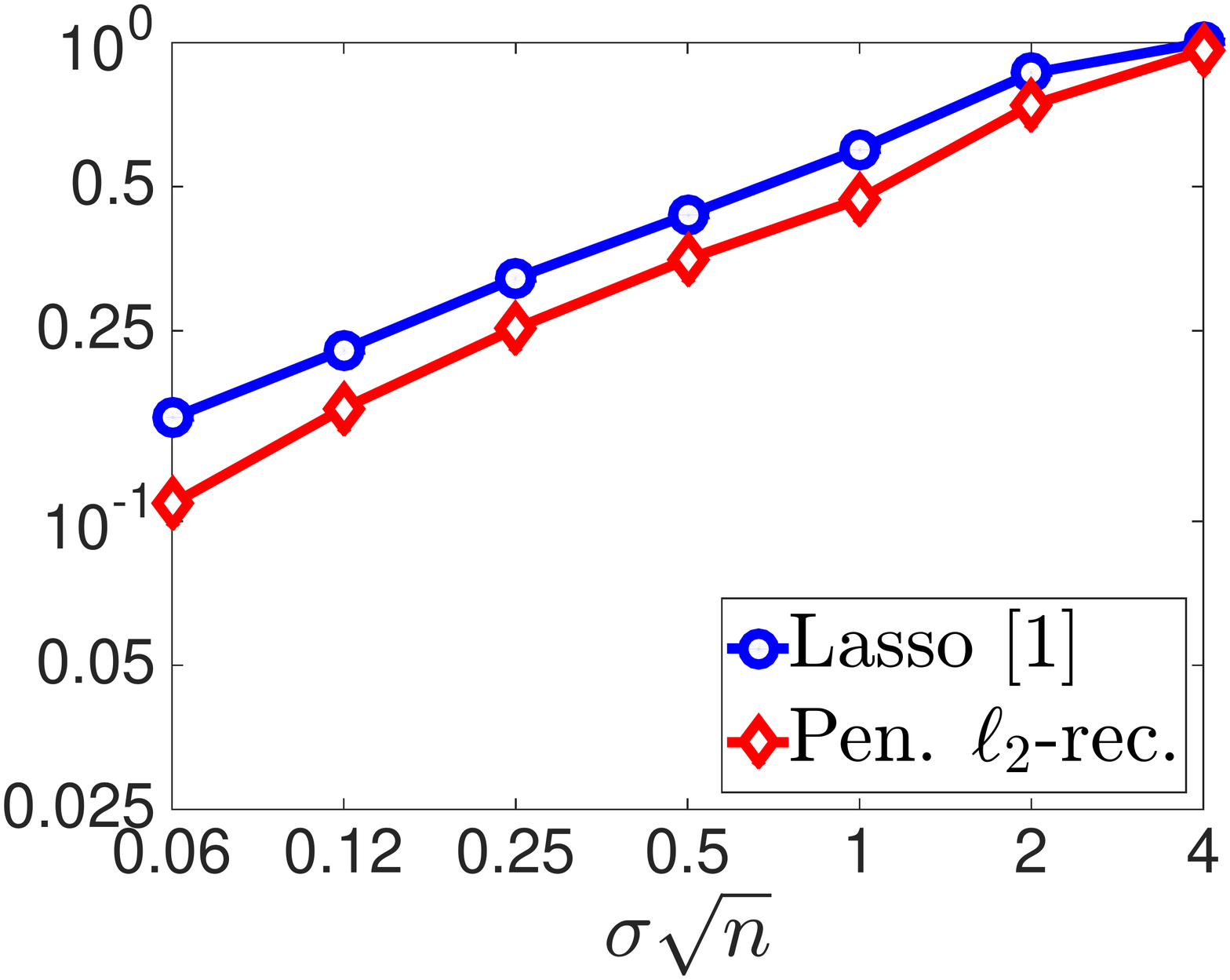}
\end{minipage}
\begin{minipage}{0.32\textwidth}
\centering
$\beta = 0.5$
\includegraphics[width=1\textwidth, height=0.18\textheight, clip=true, angle=0]{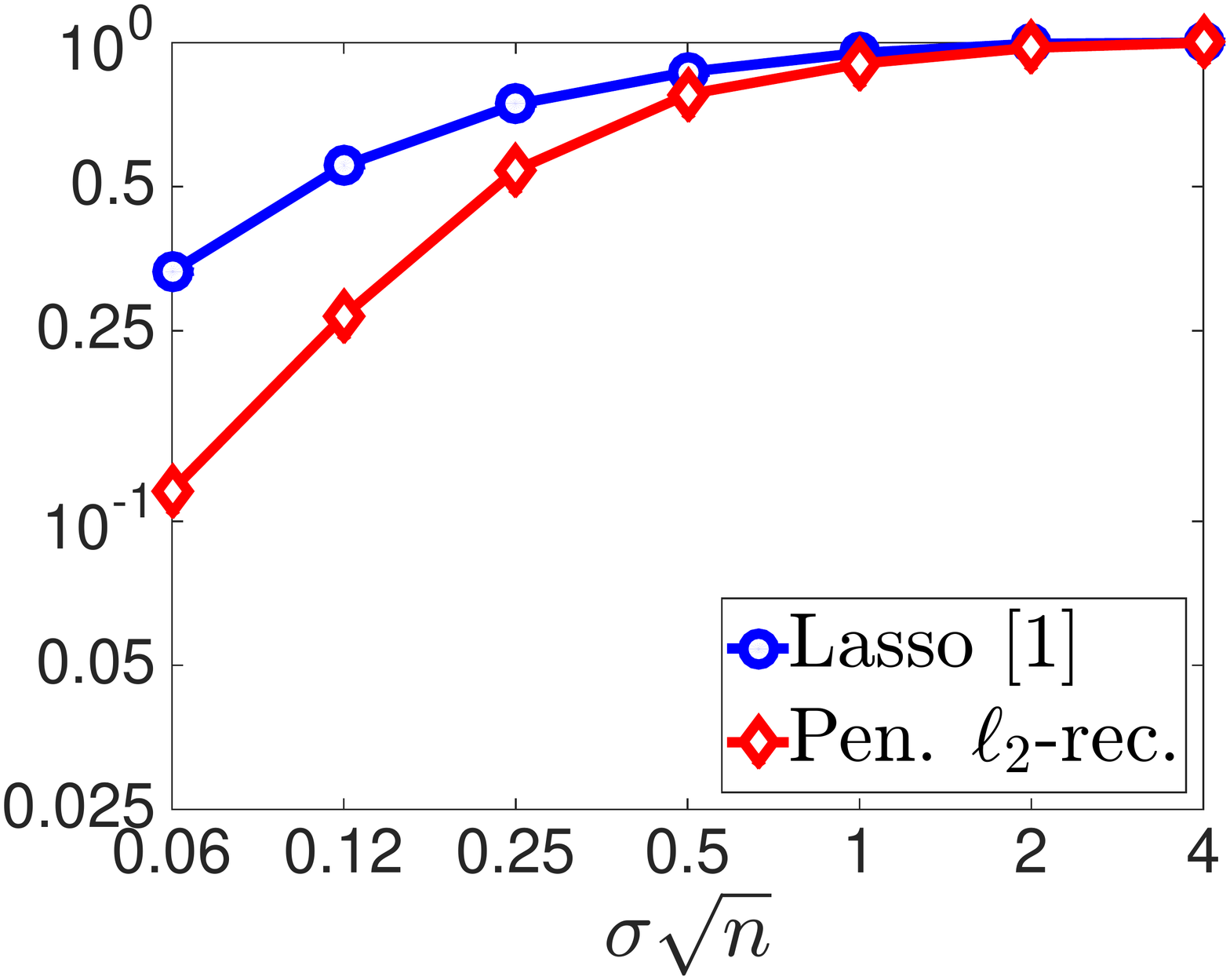}
\end{minipage}
\caption{Image denoising in \textit{DimensionReduction-2D} scenario; smoothness decreases left to right.}
\label{fig:mc-SI}
\end{figure}

As shown in Fig.~\ref{fig:mc-SI}, the $\ell_2$-recovery performs well in this scenario despite the fact that the available theoretical bounds are pessimistic.
For example, the signal \eqref{eq:single-index} with a smooth $g$ can be approximated by a small number of harmonic oscillations in $\R^2$. As follows from the proof of \cite[Proposition 10]{JN2009} combined with Theorem \ref{th:sines}, for a sum of $k$ harmonic oscillations in $\R^d$ one can point out a reproducing linear filter with $\varrho(k) = O(k^{2d})$ (neglecting the logarithmic factors), i.e. the theoretical guarantee is quite conservative for small values of $\beta$.

\section{Details of algorithm implementation}\label{sec:expdetails}
Here we give a brief account of some techniques and implementation tricks exploited in our codes.
\paragraph{Solving the optimization problems}
Note that the optimization problems~\eqref{eq:l2-1} and \eqref{optim2pen} underlying the proposed recovery algorithms are well structured  Second-Order Conic Programs (SOCP) and can be solved using Interior-point methods (IPM). However, the computational complexity of IPM applied to SOCP with {\em dense} matrices grows rapidly with problem dimension, so that large problems of this type arising in signal and image processing are well beyond the reach of these techniques. On the other hand, these problems possess nice geometry associated with complex $\ell_1$-norm. Moreover, their {\em first-order information} -- the value of objective and its gradient at a given $\varphi$ -- can be computed using Fast Fourier Transform in time which is almost linear in problem size. Therefore, we used first-order optimization algorithms, such as Mirror-Prox and Nesterov's accelerated gradient algorithm (see \cite{nesterov2013first} and references therein) in our recovery implementation. A complete description of the application of these optimization algorithms to our problem is beyond the scope of the paper; we shall present it elsewhere.
\paragraph{Interpolating recovery}
In Sec.~\ref{sec:oreg}--\ref{sec:open} we considered only recoveries which estimated the value $x_t$ of the signal via the observations at $n+1$ points $t-n,...,t$ ``on the left'' (filtering problem). To recover the whole signal, one may consider a more flexible alternative -- \textit{interpolating} recovery -- which estimates $x_t$ using observations on the left {\em and} on the right of $t$. In particular, if the objective is to recover a signal on the interval $\{-n,...,n\}$, one can apply interpolating recoveries which use the same observations $y_{-n},...,y_n$ to estimate $x_\tau$ at any $\tau\in\{-n,...,n\}$, by altering the relative position of the filter and the current point.
\paragraph{Blockwise recovery}
Ideally, when using pointwise recovery, a specific filter is constructed for each time instant $t$. This may pose a tremendous amount of computation, for instance, when recovering a high-resolution image.  Alternatively, one may split the signal into blocks, and process the points of each block using the same filter (cf. e.g. Theorem \ref{th:l2con}). For instance, a one-dimensional signal can be divided into blocks of length, say, $2m+1$, and to recover $x\in \C(\Z_{-m}^m)$ in each block one may fit one filter of length $m+1$ recovering the right ``half-block'' $x_{0}^m$ and another filter recovering the left ``half-block'' $x_{-m}^{-1}$.
\section{Conclusion}
We introduced a new family of estimators for structure-blind signal recovery that can be computed using convex optimization. The proposed estimators enjoy oracle inequalities for the $\ell_2$-risk and for the pointwise risk. Extensive theoretical discussions and numerical experiments will be presented in the follow-up journal paper. 

\section*{Acknowledgments}
We would like to thank Arnak Dalalyan and Gabriel Peyr\'{e} for fruitful discussions. DO, AJ, ZH were supported by the LabEx PERSYVAL-Lab (ANR-11-LABX-0025) and the project Titan (CNRS-Mastodons). ZH was also supported by the project Macaron (ANR-14-CE23-0003-01), the MSR-Inria joint centre, and the program ``Learning in Machines and Brains'' (CIFAR). Research of AN was supported by NSF grants CMMI-1262063, CCF-1523768.

%%%Biblio%%%
\bibliographystyle{abbrv}
\bibliography{references}
%%%
\newpage
%%%Appendix%%%
\appendix

%Commented Appendix Begins
\section{Preliminaries}
We begin by introducing several objects used in the sequel. We denote $\langle\cdot,\cdot\rangle$ the Hermitian scalar product: for $a,b\in \C^n$,
$\langle a,b\rangle=a^Hb$. For $a, b\in \C(\Z)$ we reserve the shorthand notation
\[
\langle a,b\rangle_n:=\langle a_0^n,b_0^n\rangle=[a_0^n]^Hb_0^n.
\]

\paragraph{Convolution matrices} We will extensively use the matrix-vector representation of the discrete convolution.
\begin{itemize}
\item
Given $y \in \C(\Z)$, we associate to it an $(n+1)\times(m+1)$ Toeplitz matrix
\begin{equation}\label{s(psi)}
T(y)=
\left[\begin{array}{llll}
y_{0}&y_{-1}&...&y_{-m}\\
y_{1}&y_{0}&...&y_{1-m}\\
...&...&...&...\\
y_{n}&y_{n-1}&...&y_{n-m}
\end{array}\right].
\end{equation}
such that $[\varphi*y]_{0}^n = T(y) \varphi_{0}^m$ for $\varphi \in \C(\Z_0^m)$. Its squared Frobenius norm satisfies
\begin{equation}\label{tr21}
\left\|T(y)\right\|_F^2 = \sum\limits_{\tau=0}^m\|\Delta^\tau y\|^2_{n,2}.
\end{equation}
\item
Given $\varphi \in \C(\Z_0^m)$, consider an $(n+1)\times(m+n+1)$ matrix
\begin{equation}\label{M(phi)}
M(\varphi)=
\left[\begin{array}{llllllll}
\varphi_{m}&\varphi_{m-1}&...&\varphi_{0}&0&0&...&0\\
0&\varphi_{m}&...&...&\varphi_{0}&0&...&0\\
...&...&...&...&...&...&...&...\\
0 &...&...&...&0&\varphi_{m}&...&\varphi_{0}
\end{array}\right].
\end{equation}
For $y \in \C(\Z)$ we have $[\varphi*y]_{0}^n=M(\varphi) y_{-m}^n$ and
\begin{equation}\label{tr22}
\left\|M(\varphi)\right\|_F^2= (n+1)  \|\varphi\|^2_{m,2}.
\end{equation}
\item
Given $\varphi \in \C(\Z_0^m)$, consider the following circulant matrix of size $m+n+1$:
\begin{equation}\label{C(phi)}
C(\varphi) =
\left[\begin{array}{lllllllll}
\varphi_{m}&...&...&\varphi_{0}&0&0&...&...&0\\
0&\varphi_{m}&...&...&\varphi_{0}&0&...&...&0\\
...&...&...&...&...&...&...&...&...\\
0 &...&...&...&0&\varphi_{m}&...&...&\varphi_{0}\\
\varphi_{0} &0&...&...&...&0&\varphi_{m}&...&\varphi_{1}\\
...&...&...&...&...&...&...&...&...\\
\varphi_{m-1}&...&\varphi_{0}&0&...&...&...&0&\varphi_{m}\\
\end{array}\right].
\end{equation}
One has
\begin{equation*}
\left\|C(\varphi)\right\|_F^2 = (m+n+1) \|\varphi\|^2_{m,2}.
\end{equation*}
This matrix is useful since $C(\varphi) y_{-m}^n$ encodes the circular convolution of $y_{-m}^n$ and the zero-padded filter $\varphi_0^{m+n}$ (recall that $\varphi \in \C(\Z_0^m)$) which is diagonalized by the DFT. Specifically,
\begin{equation}\label{eq:circ_diag}
C(\varphi) = F^H_{m+n} D(\varphi) F_{m+n}, \quad \text{where}\,\,\, D(\varphi) = \sqrt{m+n+1}\, \text{diag}(F_{m+n}\, \varphi_0^{m+n}).
\end{equation}
\end{itemize}

\paragraph{Deviation bounds} We use the following simple facts about Gaussian random vectors.
\begin{itemize}
\item
Let $\zeta \sim \C\N(0,I_n)$ be a standard complex Gaussian vector meaning that $\zeta = \xi_1 + \imath\xi_2$ where $\xi_{1,2}$ are two independent draws from $\N(0,I_n)$. We will use a simple bound
\begin{equation}\label{eq:complex_gaussian_max}
\prob\left\{\|\zeta\|_{\infty} \le \sqrt{2\ln n + 2u}\right\} \ge 1-e^{-u}
\end{equation}
which may be checked by explicitly evaluating the distribution since $|\zeta_1|^2_2 \sim \text{Exp}(1/2)$.\newline
\item
The following deviation bounds for $\|\zeta\|^2_2 \sim \chi^2_{2n}$ are due to \cite[Lemma 1]{lama2000}:
\begin{equation}\label{eq:chi_bound_low}
\prob\left\{\frac{\|\zeta\|^2_2}{2} \le n + \sqrt{2nu} + u \right\} \ge 1-e^{-u}, \quad \prob\left\{\frac{\|\zeta\|^2_2}{2} \ge n - \sqrt{2nu} \right\} \ge 1-e^{-u}.
\end{equation}
By simple algebra we obtain an upper bound for the norm:
\begin{equation}\label{eq:chi_bound}
\prob\left\{\|\zeta\|_2 \le \sqrt{2n} + \sqrt{2u} \right\} \ge 1-e^{-u}.
\end{equation}
\item
Further, let $K$ be an $n \times n$ Hermitian matrix with the vector of eigenvalues $\lambda = [\lambda_1;\, ...;\, \lambda_n]$. Then the real-valued quadratic form $\zeta^H K \zeta$ has the same distribution as $\xi^T B\xi$, where $\xi = [\xi_1; \xi_2] \sim \N(0,I_{2n})$, and $B$ is a real $2n \times 2n$ symmetric matrix with the vector of eigenvalues $[\lambda; \lambda]$. Hence, we have $\text{Tr}(B) = 2\text{Tr}(K)$, $\|B\|^2_F = 2\|K\|_F^2$ and $\|B\| = \|K\| \le \|K\|_F$, where $\|\cdot\|$ and $\|\cdot\|_F$ denote the spectral and the Frobenius norm of a matrix. Invoking \cite[Lemma 1]{lama2000} again (a close inspection of the proof shows that the assumption of positive semidefiniteness can be relaxed), we have
\begin{equation}\label{trq2}
\prob\left\{\frac{\zeta^H K \zeta}{2} \leq \Tr(K) + (u + \sqrt{2u})\|K\|_F\right\} \geq 1-e^{-u}.
\end{equation}
Further, when $K$ is positive semidefinite, we have $\|K\|_F \le \text{Tr}(K)$, whence
\begin{equation}\label{trq1}
\prob\left\{\frac{\zeta^H K \zeta}{2} \leq \Tr(K)(1+\sqrt{u})^2\right\} \geq 1-e^{-u}.
\end{equation}
\end{itemize}
\section{Proofs for Sections~\ref{sec:oreg} and \ref{sec:open}}\label{sec:proracle}
\subsection{Proof idea}
Despite the striking similarity with the Lasso \cite{tibshirani}, \cite{candes2007dantzig}, \cite{bickel2009simultaneous}, the recoveries of Sections~\ref{sec:oreg} and \ref{sec:open} are of quite different nature. First of all, the $\ell_1$-minimization in these methods is aimed to recover a filter but not the signal itself, and this filter is not sparse.\footnote{Unless we consider recovery of signal composed of harmonic oscillations with frequencies on the DFT grid.} The equivalent of ``regression matrices'' involved in these methods cannot be assumed to satisfy {\em Restricted Eigenvalue} or {\em Restricted Isometry} conditions, usually imposed to prove statistical properties of ``classical'' $\ell_1$-recoveries (see e.g. \cite{candes2006stable}, \cite{vandegeer2009}, and references therein). Moreover, being constructed from the noisy signal itself, these matrices depend on the noise, what introduces an extra degree of complexity in the analysis of the properties of these estimators.
Yet, proofs of  Theorem \ref{th:l2con} and \ref{th:l2pen}
rely on some simple ideas and it may be useful to expose these ideas stripped from the technicalities of the complete proof. Given $y\in\C(\Z_{-n}^n)$ let  $T(y)$ be
the $(n+1)\times(n+1)$ ``convolution matrix,'' as defined in \eqref{s(psi)}
such that for $\varphi \in \C(\Z_0^n)$ $[\varphi*y]_{0}^n = T(y) \varphi_{0}^n$. When denoting
$f=F_n\varphi$, the optimization problem in \eqref{eq:l2-1} can be recast as
a ``standard'' $\ell_1$-constrained least-squares problem with respect to $f$:
\be
\min_{f\in \C^{n+1}}\left\{\|y-A_nf\|_2^2:\;\|f\|_1\leq \overline{\varrho}/\sqrt{n+1}\right\}
\ee{eq:f2p1}
where $A_n=T(y)F^{-1}_n$. Observe that $f^o=F_n\varphi^o$ is feasible in \eqref{eq:f2p1} so that
\[
\|y-A_n\wh{f}\|_{n,2}^2\leq \|y-A_nf^o\|_{n,2}^2,
\]
where $\wh{f}=F_n\wh{\varphi}$, so that
\[\begin{array}{l}
\|x-A_n\wh{f}\|_{n,2}^2-\|x-A_nf^o\|_{n,2}^2\leq 2\sigma \left(
\Re\langle\zeta,x-A_nf^o\rangle_n
-\Re\langle \zeta,x-A_n\wh{f}\rangle_n\right)\\\leq
2\sigma \big|\langle\zeta,A_n(f^o-\wh{f})\rangle_n\big|\leq
2\sigma \|A_n^H\zeta_{0}^n\|_\infty\|f^o-\wh{f}\|_1\leq 4\sigma \|A_n^H\zeta_{0}^n\|_\infty{\frac{\overline{\varrho}}{\sqrt{n+1}}}
\end{array}\]
In the ``classical'' situation, where $\zeta_{0}^n$ is independent of $A_n$ (see, e.g., \cite{judnem2000}), the norm $\|A_n^H\zeta_{0}^n\|_\infty$ is bounded by $c_\alpha \sqrt{\ln n}\max_{j}\|[A_n]_j\|_2\leq c_\alpha \sqrt{n\ln n}\|A_n\|_\infty$ where $\|A\|_\infty=\max_{i,j}|A_{ij}|$ and $c_\alpha$ is a logarithmic in $\alpha^{-1}$ factor. This would rapidly lead to the bound \eqref{eq:th1} of the theorem.
In the case we are interested in, where $A_n$ incorporates observations $y_{-n}^n$ and thus depends on $\zeta_0^n$, curbing the cross term is more involved and requires extra assumptions, e.g. Assumption A.
\subsection{Extended formulation}
We will prove a simple generalization of Theorems~\ref{th:l2con} and \ref{th:l2pen} in the case where the length $n$ of ``validation sample'' may be different from the length $m$ of the adjusted filter.
We consider the ``skewed'' sample
\begin{align}\label{prob_model_mn}
y_\tau = x_\tau + \sigma \zeta_\tau \quad -m \le \tau \le n,
\end{align}
with $\kappa_{m,n} := \sqrt{\frac{n+1}{m+1}}$.
Accordingly, we assume that the regular recovery $\widehat\varphi * y$ uses the filter
\be
\wh\varphi\in \Argmin_{\varphi\in \C(\Z^{m}_0)}\big\{\|y-\varphi*y\|_{n,2}:\;
\|\varphi\|_{m,1}^*\leq \overline{\varrho}/\sqrt{m+1}\big\}.
\ee{eq:l2-1-mn}
The corresponding modification of Assumption A is as follows:
\paragraph{Assumption A$'$}
{\em Let $\mathcal{S}$ be a (unknown) shift-invariant linear subspace of $\C(\Z)$, $\Delta \mathcal{S} = \mathcal{S}$, of dimension $s$, $1 \le s \le n+1$.
We suppose that $x$ admits the decomposition:
\[
x=x^{\mathcal{S}}+\varepsilon,
\]
where $x^{\mathcal{S} }\in \mathcal{S}$, and $\varepsilon$ is ``small'', namely,
\[
{\left\|\Delta^{\tau} \varepsilon\right\|_{n,2}
}\leq \sigma\varkappa, \quad 0 \le \tau \le m.
\]
}
In what follows we use the following convenient reformulation of Assumption A$'$ (reformulation of Assumption A when $m=n$):
\begin{quotation}\noindent
{\em There exists an $s$-dimensional (complex) subspace $\mathcal{S}_{n} \subset \C^{n+1}$
and  an idempotent Hermitian $(n+1)\times (n+1)$ matrix $\Pi_{\mathcal{S}_{n}}$ of rank $s$ -- the projector on $\S_n$ -- such that
\be
\big\|\left(I_{n+1}-\Pi_{\mathcal{S}_{n}}\right) \left[\Delta^\tau x\right]_0^{n}\big\|_2 \, \Big[={\left\|\Delta^{\tau} \varepsilon\right\|_{n,2}
}\Big] \leq \sigma\varkappa,\;\tau=0, ...,m
\ee{subspace1}
where $I_{n+1}$ is the $(n+1)\times (n+1)$ identity matrix.
}
\end{quotation}
\begin{theorem}\label{th:l2conl}
Let $m, n\in \Z_+$, $\overline{\varrho}\geq 1$, and let $\varphi^o\in\C(\Z^{m}_0)$ be such that
\[
\|\varphi^o\|^*_{1}\leq {\overline{\varrho}/\sqrt{m+1}}.
\]
Suppose that Assumption A$'$ holds for some $s\in \Z_+$ and $\varkappa<\infty$. Then for any $\alpha, 0<\alpha\leq 1$, there is a set $\Xi\subset \C^{m+n+1}$,
$\Prob\{\zeta_{-m}^n\in \Xi \}\geq 1-\alpha$, of ``good realisations'' of $\zeta$ such that whenever $\zeta_{-m}^n\in \Xi$,
\begin{equation}\label{eq:l2con-stat}
\|x - \widehat{\varphi}*y \|_{n,2} \le \|x-\varphi^{o}*y\|_{n,2}
+ 2 \sigma \left[ \sqrt{\overline\varrho V^2_\alpha + (\overline{\varrho}+1)c_\alpha \varkappa} + \sqrt{2s}+ c_\alpha \right],
\end{equation}
where $c_\alpha := \sqrt{2\ln [16/\alpha]}$, and $$V^2_\alpha = 2\left(1 + 4\kappa_{m,n}\right)^2\ln\left[{55(m+n+1)}/{\alpha}\right].$$
\end{theorem}

\begin{theorem}\label{th:l22penl}
Let $m, n \in \Z_+$, and let $\varphi^o\in\C(\Z^{m}_0)$ be such that $$\|\varphi^o\|^*_{m,1} \le \varrho/\sqrt{m+1}$$ for some $\varrho \ge 1$. Denote $\widehat\varrho = \sqrt{m+1}\,\|\widehat\varphi\|^*_{m,1}$.\\
$1^o$. Suppose that  Assumption A$'$ holds for some $s\in \Z_+$ and $\varkappa<\infty$, and the regularization parameter of penalized recovery with $q=2$ %$\wh{x}^{\LS2}_\lasso$
satisfies $\lambda \ge \sigma^2 \Lambda_{\alpha}$, where
\begin{equation*}
\Lambda_{\alpha} := 4\sqrt{2}\left(1+2.25\kappa_{m,n}\right)^2 \ln\left[{21(m+n+1)}/{\alpha}\right].
\end{equation*}
Then for any $\alpha, 0<\alpha\leq 1$, there is a set $\Xi\subset \C^{m+n+1}$,
$\Prob\{\zeta_{-m}^n\in \Xi \}\geq 1-\alpha$, of ``good realisations'' of $\zeta$ such that whenever $\zeta_{-m}^n\in \Xi$, for the same $c_\alpha$ as in Theorem~\ref{th:l2conl}, it holds
\begin{equation}\label{eq:l22pen-stat}
\|x - \widehat{\varphi}*y \|_{n,2} \le \|x-\varphi^{o}*y\|_{n,2}
+ 2\sqrt{\varrho {\lambda}} + 2 \sigma \left[\sqrt{(\widehat\varrho + 1) c_\alpha \varkappa} + \sqrt{2s} +  c_\alpha \right].
\end{equation}
$2^o$. Moreover, if $\varkappa\leq \frac{\Lambda_{\alpha}}{4c_\alpha},$ and $\lambda \ge 2\sigma^2\Lambda_{\alpha},$ one has
\begin{equation}\label{eq:l22pen-stat-simple}
\|x - \widehat{\varphi}*y \|_{n,2}  \le \|x-\varphi^{o}*y\|_{n,2}
+ \sqrt{3\varrho\lambda} +2\sigma[\sqrt{2s} + c_\alpha].
\end{equation}
\end{theorem}

\subsection{Proof of Theorem \ref{th:l2conl}}
\paragraph{1$^o$.} The oracle filter $\varphi^{o}$ is feasible in \eqref{eq:l2-1-mn}, hence,
\begin{eqnarray}\label{first0}
\|x-\wh\varphi*y\|_{n,2}^2&\ci{=}{\leq}& \|(1-\varphi^{o})*y\|_{n,2}^2-{\sigma^2\|\zeta}\|_{n,2}^2-2\,\Re{\sigma\langle \zeta}, x-\wh\varphi*y\rangle_n \nn
&=&\|x-\varphi^{o}*y\|_{n,2}^2  - 2\underbrace{\Re{\sigma\langle \zeta}, x-\wh\varphi*y\rangle_n }_{{\delta^{(1)}}} + 2\underbrace{\Re{\sigma\langle \zeta}, x-\varphi^{o}*y\rangle_n }_{{\delta^{(2)}}}.
\end{eqnarray}
Let us bound ${\delta^{(1)}}$.
Denote for brevity $I := I_{n+1}$, and recall that $\Pi_{\cS_n}$ is the projector on $\mathcal{S}_n$ from \eqref{subspace1}. We have the following decomposition:
\begin{eqnarray}\label{I_pi}
{\delta^{(1)}}
&=&\underbrace{\Re{\sigma\langle \zeta_0^n}, \Pi_{\cS_n}[x-\wh\varphi*y]_0^n\rangle}_{{\delta^{(1)}_1}} +\underbrace{\Re{\sigma\langle \zeta_0^n}, (I-\Pi_{\cS_n})[x-\wh\varphi*x]_{0}^n\rangle}_{{\delta^{(1)}_2}}\nn &&-\underbrace{\Re{\sigma^2\langle \zeta_0^n}, (I-\Pi_{\cS_n})[\wh\varphi * {\zeta}]_{0}^n \rangle}_{{\delta^{(1)}_3}}
\end{eqnarray}
One can easily bound $\delta^{(1)}_1$ under the premise of the theorem:
\[
\left|{\delta^{(1)}_1} \right|
\leq {\sigma}\big\|\Pi_{\cS_n} {\zeta}_0^n\big\|_2  \big\|\Pi_{\cS_n}[x-\wh\varphi*y]_0^n\big\|_{2} \leq {\sigma}\big\|\Pi_{\cS_n} {\zeta}_0^n\big\|_2  \big\|x-\wh\varphi*y\big\|_{n,2}.
\]
Note that $\Pi_{\cS_n}\zeta_0^n \sim \C\N(0,I_s)$, and by \eqref{eq:chi_bound} we have
\[
\prob\left\{\big\|\Pi_{\cS_n}\zeta_0^n\big\|_2 \geq \sqrt{2s}+\sqrt{2u}\right\}\leq e^{-u},
\]
obtaining the bound
\begin{eqnarray}\label{i11}
\text{Prob}\left\{\big|{\delta^{(1)}_1}\big|
\leq \sigma\big\|x-\wh\varphi*y\big\|_{n,2} \left(\sqrt{2s}+\sqrt{2\ln\left[1/\alpha_1\right]}\right) \right\} \ge 1-\alpha_1.
\end{eqnarray}
\paragraph{2$^o$.} We are to bound the second term of \eqref{I_pi}. To this end, note first that
\begin{equation*}
{\delta^{(1)}_2} =\Re {\sigma\langle \zeta_0^n}, (I-\Pi_{\cS_n})x_0^n \rangle - \Re{\sigma\langle \zeta_0^n}, (I-\Pi_{\cS_n})[\wh\varphi*x]_0^n\rangle.
\end{equation*}
By \eqref{subspace1}, $\left\|(I-\Pi_{\cS_n})x_0^n\right\|_{2}\leq \sigma \varkappa$, thus with probability $1-\alpha$,
\begin{equation}\label{i22-1}
\left|{\langle \zeta_0^n}, (I-\Pi_{\cS_n})x_0^n\rangle \right| \le \sigma\varkappa\sqrt{2\ln[1/\alpha]}.
\end{equation}
On the other hand, using the notation defined in \eqref{s(psi)}, we have
 $[\wh\varphi*x]_0^n=T(x)\wh\varphi_{0}^m$, so that
\[
{\langle \zeta_0^n}, (I-\Pi_{\cS_n})[\wh\varphi*x]_0^n\rangle = {\langle \zeta_0^n}, (I-\Pi_{\cS_n})T(x)\wh\varphi_{0}^m\rangle.
\]
Note that $[T(x)]_{\tau} = x^{-\tau+n}_{-\tau}$ for the columns of $T(x)$, $0\le \tau \le m$. By \eqref{subspace1},
$(I-\Pi_{\cS_n})T(x)=T(\varepsilon)$, and by \eqref{tr21},
\[
\left\|(I-\Pi_{\cS_n})T(x)\right\|_F^2 = \left\| T(\varepsilon)\right\|_F^2
= \sum_{\tau=0}^{m}\big\|\varepsilon_{-\tau}^{-\tau+n}\big\|^2_2 \leq (m+1)\sigma^2\varkappa^2.
\]
Due to \eqref{trq1} we conclude that
\[
\left\|T(x)^H(I-\Pi_{\cS_n})\zeta_0^n\right\|_2^2 \le 2(m+1)
\sigma^2\varkappa^2\left(1+\sqrt{\ln[1/\alpha]}\right)^2
\]
with probability at least $1-\alpha$. Since
\[\left|\left\langle\zeta_0^n, (I-\Pi_{\cS_n})T(x)\wh\varphi_{0}^m\right\rangle\right|\leq \frac{\overline\varrho}{\sqrt{m+1}} \left\|T(x)^H(I-\Pi_{\cS_n})\zeta_0^n\right\|_{2},
\]
we arrive at the bound with probability $1-\alpha$:
\[
\left|\left\langle\zeta_0^n, (I-\Pi_{\cS_n})T(x)\wh\varphi_{0}^m\right\rangle\right|\le\sqrt{2}\sigma\varkappa\overline{\varrho} \left(1+\sqrt{\ln[1/\alpha]}\right).
\]
Along with \eqref{i22-1} this results in the following bound:
\begin{eqnarray}\label{i22}
\text{Prob}\left\{\big|{\delta^{(1)}_2}\big|
\leq \sqrt{2}\sigma^2\varkappa(\overline{\varrho}+1) \left(1+\sqrt{\ln\left[2/\alpha_{2}\right]}\right)\right\} \ge 1-\alpha_{2}.
\end{eqnarray}
\paragraph{3$^o$.} Let us rewrite ${\delta^{(1)}_3}$ as follows:
\begin{equation*}
{\delta^{(1)}_3} =\Re{\sigma^2\langle \zeta_0^n}, (I-\Pi_{\cS_n})M(\wh\varphi){\zeta}_{-m}^{n} \rangle =
\Re {\sigma^2\langle \zeta_{-m}^n}, Q M(\wh\varphi){\zeta}_{-m}^{n}\rangle,
\end{equation*}
where $M(\wh\varphi)\in \C^{(n+1)\times (m+n+1)}$ is defined by \eqref{M(phi)}, and $Q\in \C^{(m+n+1)\times (n+1)}$ is given by
\begin{equation*}
Q=\left[
\begin{array}{c c}
0_{m,n+1}; & I - \Pi_{\cS_n}
\end{array}
\right];
\end{equation*}
hereinafter we denote $0_{m,n}$ the $m\times n$ zero matrix.
Now, by the definition of $\widehat\varphi$ and since the mapping $\varphi \mapsto M(\varphi)$ is linear,
\begin{eqnarray}\label{maxR}
{\delta^{(1)}_3}
&=&{\half}({\zeta}_{-m}^{n})^H(\underbrace{Q^{\vphantom{*}}M(\wh\varphi)+M(\wh\varphi)^H Q^H}_{K_{1}(\wh\varphi)}) z_{-m}^{n}
\leq\frac{\sigma^2\overline\varrho}{2\sqrt{m+1}}\max\limits_{\scriptsize\begin{array}{c} u \in \C(\Z_0^m),\nn \|u\|^*_{m,1}\leq 1\end{array}} (\zeta_{-m}^{n})^H K_{1}(u) \zeta_{-m}^{n}\\
&=&\frac{\sigma^2\overline\varrho}{\sqrt{m+1}}\,\,\max\limits_{\tiny{1\le j \le m+1}}\,\,\max\limits_{\theta \in [0,2\pi]} \half (\zeta_{-m}^{n})^H K_{1}(e^{\imath\theta} u^j) \zeta_{-m}^{n},
\end{eqnarray}
where $u^j \in \C(\Z_0^m)$, and $[u^j]_0^m = F_m^{-1} \text{e}_j$, $\text{e}_j$ being the $j$-th canonical orth of $\R^{m+1}$. Indeed, $M(\varphi)$ attains its maximum over the convex set
\begin{equation}\label{eq:hyperoct}
\mathcal{B}_{m,1}^* = \{u \in \C(\Z_0^m), \quad \|u\|^*_{m,1} \le 1\}.
\end{equation}
at an extremal point $e^{\imath\theta} u_j$, $\theta \in [0,2\pi]$. It is easy to verify that $$K_1(e^{\imath\theta}u) = K_1(u)\cos\theta  + K_2(u)\sin\theta $$ for the Hermitian matrix
\begin{equation*}
K_2(u) = \imath\left(Q^{\vphantom{*}}M(u) - M(u)^H Q^H \right).
\end{equation*}
Denoting $q^j_i(\zeta) = \half(\zeta_{-m}^n)^H K_i(u^j) \zeta_{-m}^n$ for $i=1,2$, we have
\begin{equation}\label{eq:max_complex_to_real}
\begin{array}{ll}
&\max\limits_{\theta \in [0,2\pi]} \half (\zeta_{-m}^{n})^H K_{1}(e^{\imath\theta} u^j) \zeta_{-m}^{n} =\max\limits_{\theta \in [0,2\pi]}\left[q^j_1(\zeta)\cos\theta   + q^j_2(\zeta)\sin\theta\right]\\
&=\sqrt{\big|q_1^j(\zeta)\big|^2 + \big|q_2^j(\zeta)\big|^2} \le \sqrt{2}\max \big(\big|q_1^j(\zeta)\big|, \big|q_2^j(\zeta)\big|\big).
\end{array}
\end{equation}
By simple algebra and using \eqref{tr22}, we get
\begin{equation*}\label{trA1}
\Tr\left[K_{i}(u^j)^2\right] \leq
4\,\Tr[M(u^j)M(u^j)^H]= 4(n+1)\|u^j\|^2_{m,2}\leq 4(n+1), \quad i=1,2.
\end{equation*}
Now let us bound  $\Tr[K_{i}(u)]$, $i=1,2$, on the set \eqref{eq:hyperoct}.
One may check that for the circulant matrix $C(u)$, cf. \eqref{C(phi)}, it holds:
\begin{equation*}
QM(u) = \underbrace{Q^{\vphantom{H}}Q^H}_{R} C(u),
\end{equation*}
where $R = Q^{\vphantom{H}}Q^H$ is an $(m+n+1)\times (m+n+1)$ projection matrix of rank $s$ defined by
\begin{equation*}
R = \left[
\begin{array}{l|l}
0_{m,m} & 0_{m,n+1}\\
\hline
0_{n+1,m} & I-\Pi_{\mathcal{S}_n}
\end{array}
\right].
\end{equation*}
Hence, denoting $\|\cdot\|_*$ the nuclear norm, we can bound $\Tr[K_{i}(u)]$, $i=1,2$, as follows:
\begin{equation*}
\big|\Tr[K_{i}(u)]\big| \le 2\big|\Tr[R C(u)]\big| \le 2\|R\|  \left\|C(u)\right\|_*
\le 2\|C(u)\|_*= 2\sqrt{m+n+1}\|u\|^*_{m+n, 1},
\end{equation*}
where in the last transition we used \eqref{eq:circ_diag}.
The following technical lemma gives an upper bound on the norm of a zero padded filter (see Appendix \ref{proof:zero_padding} for the proof):
\begin{lemma}\label{lemma:zero_padding}
For any $u \in \mathcal{B}^*_{m,1}$, see \eqref{eq:hyperoct}, and $n \ge 1$, we have
\begin{equation*}
\|u\|^*_{m+n,1}\le \sqrt{1 + \kappa^2_{m,n}}  (\ln[m+n+1]+3).
\end{equation*}
\end{lemma}
\noindent Thus we arrive at
$$\left|\Tr[K_{i}(u^j)]\right| \le 2 \sqrt{m+1}(\kappa_{m,n}^2 + 1)  (\ln[m+n+1]+3), \quad i=1,2.$$
By \eqref{trq2} we conclude that for any fixed pair $(i,j) \in \{1,2\} \times \{1,\, ...,\, m+1\}$, with probability $1-\alpha$,
\[ \big|q^j_i(\zeta)\big| \leq \left|\Tr[K_{i}(u^j)]\right| +\left\|K_{i}(u^j)\right\|_F \left(1 + \sqrt{\ln[2/\alpha]}\right)^2.
\]
With $\alpha_0=2(m+1)\alpha$,  by the union bound together with \eqref{maxR} and \eqref{eq:max_complex_to_real} we get
\begin{equation}\label{i23-2}
\Prob\left\{{\delta^{(1)}_3}
\leq
2\sqrt{2}\sigma^2\overline{\varrho}\left[(\kappa_{m,n}^2+1)(\ln[m+n+1]+3)+\kappa_{m,n}\left(1+\sqrt{\ln\left[{4(m+1)}/{\alpha_0}\right]}\right)^2\right]\right\} \ge 1-\alpha_0.
\end{equation}
\paragraph{4$^o$.} Bounding ${\delta^{(2)}}$ is a relatively simple task since $\varphi^{o}$ does not depend on the noise. We decompose
\[
{\delta^{(2)}} = \sigma \Re\langle \zeta, x-\varphi^{o}*x\rangle_n - \sigma^2\Re\langle \zeta,\varphi^{o}*\zeta\rangle_n.
\]
Note that $\Re\langle \zeta,x-\varphi^{o}*x\rangle_n \sim \C\N(0,\|x-\varphi^{o}*x\|_{n,2}^2)$, therefore, with probability $1-\alpha$,
\begin{equation}\label{simple0}
\Re\langle\zeta,x-\varphi^{o}*x\rangle_n\leq\sqrt{2\ln[1/\alpha]}  \|x-\varphi^{o}*x\|_{n,2}.
\end{equation}
On the other hand,
\begin{align}\label{simple1}
\|x-\varphi^{o}*x\|_{n,2}
&\leq \|x-\varphi^{o}*y\|_{n,2}+\sigma \|\varphi^{o}*\zeta\|_{n,2}\notag\\
&\leq \|x-\varphi^{o}*y\|_{n,2}+\sqrt{2} \sigma \varrho \kappa_{m,n} \left(1+\sqrt{\ln [1/\alpha]}\right)
\end{align}
with probability $1-\alpha$. Indeed, one has
\[
\|\varphi^{o}*\zeta\|^2_{n,2}=\left\|M(\varphi^{o})\zeta_{-m}^n \right\|_2^2,
\]
where for $M(\varphi^{o})$ by $\eqref{tr22}$ we have
\begin{equation}\label{eq:trmm}
\left\|M(\varphi^{o})\right\|_F^2  = (n+1)\|\varphi^{o}\|_{m,2}^2 \leq\kappa_{m,n}^2\varrho^2.
\end{equation}
Using \eqref{trq1} we conclude that, with probability at least $1-\alpha$,
\begin{equation*}
\|\varphi^{o}*\zeta\|^2_{n,2} \leq 2\kappa_{m,n}^2\varrho^2\left(1+\sqrt{\ln[1/\alpha]}\right)^2,
\end{equation*}
which implies \eqref{simple1}. Using \eqref{simple0} and \eqref{simple1}, we get that with probability at least $ 1-\alpha_{3}$,
\begin{eqnarray}\label{simple2}
\lefteqn{\Re\langle\zeta,x-\varphi^{o}*x\rangle_n}\nn
&\leq& \sqrt{2\ln \left[2/\alpha_{3}\right]} \left[\|x-\varphi^{o}*y\|_{n,2}+\sqrt{2}\sigma\varrho\kappa_{m,n}
\left(1+\sqrt{\ln\left[2/\alpha_{3}\right]}\right)\right]\nn
&\leq& \|x-\varphi^{o}*y\|_{n,2}\sqrt{2\ln\left[2/\alpha_{3}\right]}+2\sigma\varrho \kappa_{m,n}\left(1+\sqrt{\ln \left[2/\alpha_{3}\right]}\right)^2.
\end{eqnarray}
The indefinite quadratic form
\begin{equation*}
\Re\langle\zeta,\varphi^{o}*\zeta\rangle_n=\frac{(\zeta_{-m}^n)^H{K_0(\varphi^{o})}\zeta_{-m}^n}{2},
\end{equation*}
where $K_0(\varphi^{o}) =
\left[0_{m,m+n+1};\, M(\varphi^{o})\right]+
\left[0_{m,m+n+1}; \, M(\varphi^{o})\right]^H$,
can be bounded similarly to 3$^o$. We get
\[
|\Tr[{K_0(\varphi^{o})}]| \le 2(n+1)\left|\varphi^{o}_m\right|\leq 2\kappa_{m,n}^2\varrho.
\]
Indeed, for $e_{m+1} = [0;\, ...;\,0;\,1] \in \R^{m+1}$ one has
\begin{equation*}
\left|\varphi^{o}_m\right| = |\langle [\varphi^{o}]_0^m, {e}_{m+1} \rangle| \le \|\varphi^{o}\|^*_{m,1} \|F_m {e}_{m+1}\|_{\infty} \le \frac{\varrho}{m+1}
\end{equation*}
since $\|F_m {e}_{m+1}\|_{\infty} = 1/\sqrt{m+1}$.
By \eqref{eq:trmm}, $\left\|K_0(\varphi^{o})\right\|_F^2\leq 4\left\|M(\varphi^{o})\right\|_F^2 \le 4\kappa_{m,n}^2\varrho^2.$ Hence by \eqref{trq2},
\begin{equation}\label{i112}
\text{Prob}\left\{-\Re \langle\zeta,\varphi^{o}*\zeta\rangle_n \leq 2\kappa_{m,n}^2\varrho+2\kappa_{m,n}\varrho \left(1+\sqrt{2\ln\left[1/\alpha_4\right]}\right)^2 \right\} \ge 1-\alpha_4.
\end{equation}
\paragraph{5$^o$.} It remains to combine the bounds obtained in $1^o$-$4^o$. For any $\alpha \in (0,1]$, putting $\alpha_0 = \alpha_1 = \alpha_4 = \alpha/4$, $\alpha_{2} = \alpha_{3} = \alpha/8$, and using the union bound, we get from \eqref{first0} with probability $\geq 1-\alpha$:
\begin{eqnarray}\label{eq:l2con-long}
\|x-\wh\varphi*y\|_{n,2}^2&\leq& \|x-\varphi^{o}*y\|_{n,2}^2+2{\delta^{(2)}}-2{\delta^{(1)}}\nn
\mbox{[by \eqref{simple2}]}&\leq& \|x-\varphi^{o}*y\|_{n,2}^2+2\sigma\|x-\varphi^{o}*y\|_{n,2}\sqrt{2\ln[16/\alpha]}\nn
\mbox{[by \eqref{simple2}, \eqref{i112}]}&&+4\sigma^2\varrho\Big[\kappa_{m,n}^2+2\kappa_{m,n}\big(1+\sqrt{2\ln[16/\alpha]}\big)^2\Big]
\nn
\mbox{[by \eqref{i11}]}&&+2\sigma\|x-\wh\varphi*y\|_{n,2}\big(\sqrt{2s}+\sqrt{2\ln [16/\alpha]}\big)
\nn
\mbox{[by \eqref{i22}]}&&+2\sqrt{2}\sigma^2(\overline{\varrho}+1)\big(1+\sqrt{\ln[16/\alpha]}\big)\varkappa
\nn
\mbox{[by \eqref{i23-2}]}&&+
{4\sqrt{2}\sigma^2\overline{\varrho}\left[(\kappa_{m,n}^2+1)(\ln[m+n+1]+3)
+\kappa_{m,n}\left(1+\sqrt{\ln\left[{16(m+1)}/{\alpha}\right]}\right)^2\right]}.
\nn
\end{eqnarray}
Hence, denoting
\begin{eqnarray}
u(\varrho,\varkappa)&=&2\sqrt{2}\sigma^2({\varrho}+1)\big(1+\sqrt{\ln[16/\alpha]}\big)\varkappa,\label{eq:part-u}\\
v_1(\varrho)&=&{4\sqrt{2}\sigma^2{\varrho}\left[(\kappa_{m,n}^2+1)(\ln[m+n+1]+3)
+\kappa_{m,n}\left(1+\sqrt{\ln\left[{16(m+1)}/{\alpha}\right]}\right)^2\right]},\label{eq:part-v1}\\
v_2(\varrho) &=& 4\sigma^2\varrho\Big[\kappa_{m,n}^2+2\kappa_{m,n}\big(1+\sqrt{2\ln[16/\alpha]}\big)^2\Big],\label{eq:part-v2}
\end{eqnarray} we obtain
\begin{eqnarray*}
\|x-\wh\varphi*y\|_{n,2}^2 &\le& \|x-\varphi^{o}*y\|_{n,2}^2
+ 2\sigma\big(\sqrt{2s}+\sqrt{2\ln [16/\alpha]}\big)  \big(\|x-\wh\varphi*y\|_{n,2} + \|x-\varphi^{o}*y\|_{n,2}\big)\\&&+ u(\overline\varrho,\varkappa) + v_1(\overline\varrho) + v_2(\overline\varrho).
\end{eqnarray*}
The latter implies that
\begin{align*}
\|x-\wh\varphi*y\|_{n,2} \le \|x-\varphi^{o}*y\|_{n,2}+2\sqrt{2}\sigma\big(\sqrt{s}+\sqrt{\ln [16/\alpha]}\big) + \sqrt{u(\overline\varrho,\varkappa) + v_1(\overline\varrho) +  v_2(\overline\varrho)}
\end{align*}
Finally, we arrive at \eqref{eq:l2con-stat} using the bound
\begin{equation}\label{eq:part-u-bound}
u(\overline\varrho,\varkappa) \le 4 \sqrt{2} \sigma^2 (\overline\varrho+1) \sqrt{\ln[16/\alpha]} \varkappa
\end{equation}
and
\begin{align*}%\label{eq:final_stoch_term}
v_1(\overline\varrho) + v_2(\overline\varrho)
&\le\sigma^2\overline{\varrho} \left(4\sqrt{2} (\kappa_{m,n}^2+1)(\ln[m+n+1]+4) + 4.5(4\sqrt{2}+8)\kappa_{m,n}\ln\left[16(m+1)/\alpha\right]\right)\nn
&\le 8\sigma^2\overline{\varrho} \left(1 + 4\kappa_{m,n}\right)^2 \ln\left[55(m+n+1)/\alpha\right].\tag*{\qed}
\end{align*}
%(we used $\ln[8(m+1)] \ge \ln 16 \ge 2$ as long as $m \ge 1$, and that $(1 + \sqrt{2u})^2 \le 4.5u$ when $u \ge 2$; in the second line we used $e^4 < 55$ and $4\sqrt{2} < 6$. To obtain the claim it suffices to note that

\subsection{Proof of Lemma~\ref{lemma:zero_padding}}
\label{proof:zero_padding}
The function $\|u\|^*_{m+n,1}$ is convex on \eqref{eq:hyperoct}, so its maximum over this set is attained at one of the extreme points $F_m [u^j]_0^m = e^{\imath \theta} {e}_j$ where ${e}_j$ is the $j$-th canonical orth of $\R^{m+1}$ and $\theta \in [0,2\pi]$. Since $u^j_\tau = \frac{1}{\sqrt{m+1}} \exp\left[\imath\theta-\frac{2\pi\imath \tau j}{m+1}\right],$ we obtain
\begin{align*}
\big\|u^j\big\|^*_{m+n,1}
= \frac{1}{\gamma} \sum_{k=0}^{m+n} \Bigg| \sum_{\tau=0}^m \exp\Bigg[\imath\underbrace{2\pi\left(\frac{k}{m+n+1} - \frac{j}{m+1}\right)}_{\omega_{jk}}\tau\Bigg] \Bigg|
&= \frac{1}{\gamma} \sum_{k=0}^{m+n} \left| D_m\left(\frac{\omega_{jk}}{2}\right)\right| ,
\end{align*}
where $\gamma =\sqrt{(m+n+1)(m+1)}$, and the Dirichlet kernel $D_m(\cdot)$ is defined as
\begin{align*}
D_m(x) := \left\{
\begin{array}{ll}
\frac{\sin ((m+1)x)}{\sin(x)}, & \quad x \ne \pi l,\\
m+1, & \quad x = \pi l.
\end{array}\right.
\end{align*}
Hence, $\gamma\|u^j\|^*_{m+n,1} \le \max_{\epsilon \in [0,\pi]} S_{m+n}(\epsilon)$, where
\begin{align}\label{eq:dirichlet_sum}
S_{m,n}(\epsilon) = \sum_{k=0}^{m+n} \left| D_m\left(\frac{\pi k}{m+n+1} -\epsilon\right)\right|, \quad \epsilon \in [0, \pi].
\end{align}
Note that $|D_m(x)|$ is upper bounded  by the following (positive) function on the circle $\R/\pi \Z$:
\begin{align*}
B_m(x) = \left\{
\begin{array}{ll}
\frac{\pi}{2 \min (x, \pi-x)}, & \quad x \in (0,\pi),\\
m+1, & \quad x = 0.
\end{array}
\right.
\end{align*}
For any $\epsilon \in [0,\pi]$, the summation in \eqref{eq:dirichlet_sum} is over a regular  $(m+n+1)$-grid on $\R/\pi \Z$. The contribution to the sum of each of two closest to zero points of the grid is at most $\max_{x \in [0,\pi]} D_m(x) = m+1$. For the remaining points, we can upper bound $|D_m(x)| \le B_m(x)$ noting that $B_m(x)$ decreases over $[\frac{\pi}{m+n+1}, \frac{\pi}{2}]$ as long as $n \ge 1$. These considerations result in
\begin{align*}
S_{m,n}(\epsilon)
&\le 2(m+1) + \sum_{k=1}^{\left\lceil \frac{m+n-1}{2}\right\rceil} \frac{m+n+1}{k}
\le 2(m+1) + (m+n+1) \left(\ln\left(\frac{m+n+1}{2}\right)+1\right)
\end{align*}
where in the last transition we used the simple bound $H_n \le \ln n + 1$ for harmonic numbers.
\qed

\subsection{Proof of Theorem \ref{th:l22penl}}
We will use the same notation as in the proof of Theorem \ref{th:l2conl}. Due to feasibility of $\varphi^{o}$, we have the following counterpart of  \rf{first0}:
\begin{equation*}
\|x-\wh\varphi*y\|_{n,2}^2 + \lambda\widehat\varrho
\leq\|x-\varphi^{o}*y\|_{n,2}^2  - 2\delta^{(1)} + 2\delta^{(2)} + \lambda \varrho.
\end{equation*}
Thus, repeating steps $1^o-4^o$ of the proof of Theorem~\ref{th:l2conl}, we obtain, cf. \eqref{eq:l2con-long}, \eqref{eq:part-u}, \eqref{eq:part-v1}, \eqref{eq:part-v2},
\begin{eqnarray*}%\label{eq:l22pen-long}
\|x-\wh\varphi*y\|_{n,2}^2 \leq\, &\|x-\varphi^{o}*y\|_{n,2}^2+2\sigma\|x-\varphi^{o}*y\|_{n,2}\sqrt{2\ln[16/\alpha]} - \lambda\widehat\varrho + \lambda \varrho\nn
&+2\sigma\|x-\wh\varphi*y\|_{n,2}\big(\sqrt{2s}+\sqrt{2\ln [16/\alpha]}\big)\nn
&+u(\widehat\varrho,\varkappa)+v_1(\widehat\varrho) + v_2(\varrho).
\end{eqnarray*}
 Using that $(1+\sqrt{2x})^2 \le 4.5x$ when $x \ge 2$, one may check that as long as $m, n \ge 1$,
\begin{equation*}%\label{eq:l2pen-tokill}
v_1(\widehat{\varrho})\le 4\sqrt{2}\sigma^2\widehat\varrho \left(\kappa_{m,n}^2 + 4.5\kappa_{m,n} + 1\right) \ln\left[{21(m+n+1)}/{\alpha}\right] \le \sigma^2\widehat\varrho\Lambda_{\alpha}.
\end{equation*}
Hence, choosing $\lambda \geq \sigma^2\Lambda_{\alpha}$, one guarantees $v_1(\widehat\varrho) - \lambda\widehat\varrho \le 0$. Using also
\begin{equation*}%\label{eq:final_long-1-bound}
v_2(\varrho) \le 4\sigma^2\varrho(\kappa_{m,n}^2 + 9\kappa_{m,n}\ln[16/\alpha]) \le 3 \sigma^2\varrho\Lambda_{\alpha},
\end{equation*}
one arrives at
\begin{eqnarray*}
\|x-\wh\varphi*y\|_{n,2}^2 &\le& \|x-\varphi^{o}*y\|_{n,2}^2 + 2\sigma\left(\|x-\wh\varphi*y\|_{n,2} + \|x-\varphi^o*y\|_{n,2}\right)\big(\sqrt{2s}+\sqrt{2\ln [16/\alpha]}\big)\\
&&+ u(\widehat\varrho,\varkappa) + \lambda \varrho + 3\sigma^2\varrho\Lambda_{\alpha},
\end{eqnarray*}
whence \eqref{eq:l22pen-stat} follows by $\eqref{eq:part-u-bound}$ and $\lambda \ge \sigma^2\Lambda_{\alpha}$. To prove \eqref{eq:l22pen-stat-simple}, note that if $\varkappa$ is bounded from above as in the premise of the theorem, one has, by $\eqref{eq:part-u-bound}$,
\begin{equation*}%\label{eq:varkappa-bias-bound}
u(\widehat\varrho,\varkappa) \leq 4 \sqrt{2} \sigma^2 (\widehat\varrho+1) \sqrt{\ln[16/\alpha]} \varkappa \le \sigma^2\Lambda_{\alpha} (\widehat\varrho+1) \le \sigma^2\Lambda_{\alpha}(\widehat\varrho + \varrho),
\end{equation*}
thus arriving at
\begin{eqnarray*}
\|x-\wh\varphi*y\|_{n,2}^2 &\le& \|x-\varphi^{o}*y\|_{n,2}^2 + 2\sigma\left(\|x-\wh\varphi*y\|_{n,2} + \|x-\varphi^o*y\|_{n,2}\right)\big(\sqrt{2s}+\sqrt{2\ln [16/\alpha]}\big)\\
&&+ (\lambda  + 4\sigma^2\Lambda_{\alpha})\varrho + (2\sigma^2\Lambda_{\alpha} - \lambda)\widehat\varrho.
\end{eqnarray*}
Whence \eqref{eq:l22pen-stat-simple} follows by simple algebra using $\lambda \ge 2\sigma^2\Lambda_{\alpha}$.
\qed

\subsection{Proof of Theorem \ref{th:point}}
We decompose
\begin{eqnarray}\label{eq:l2point-long}
\lefteqn{|[x - \widehat\varphi * y]_n|=|[(\phi^o+(1-\phi^o))*(x - \widehat\varphi * y)]_n|}\nn&\leq&
|[\phi^o*(x - \widehat\varphi * y)]_n| + |[(1-\widehat \varphi) * (1 - \phi^o) * x]_n|
+ \sigma|[\widehat \varphi * \zeta]_n| + \sigma|[\widehat \varphi * \phi^o * \zeta]_n|\nn
&:=&\delta^{(1)}+\delta^{(2)}+\delta^{(3)}+\delta^{(4)}.
\end{eqnarray}
We have
\[
\delta^{(1)}\leq \|\phi^o\|_2\|x - \widehat\varphi * y)\|_{m,2}\leq{\frac{\rho}{\sqrt{m+1}}} \|x - \widehat\varphi * y\|_{m,2}.
\]
When using the bound
of Corollary  \ref{cor:l2con} with $\overline{\varrho}=2\rho^2$, we conclude that, with probability $\ge 1-\alpha/3$,
\[
\delta^{(1)}\leq  c{\sigma\rho\over \sqrt{n}}\left[\rho^2\sqrt{\ln[1/\alpha]}+ \rho\sqrt{ \big(\varkappa\ln\big[{1/ \alpha}\big]+\ln\big[{n/ \alpha}\big]\big)}+
\sqrt{s}\right].
\]
Next we get
\[
\delta^{(2)}\le \left(1 + \left\|\widehat\varphi\right\|_1\right) \left\|(1 - \phi^o) * x\right\|_{n,\infty} \le (1+2\rho^2){\sigma\rho\over \sqrt{m+1}}.
\]
By the Parseval identity,
\begin{align*}
\delta^{(3)} =\sigma|\langle F_n[\widehat\varphi^*]_{0}^n,F_n[\zeta]_{0,n}\rangle|\leq \sigma\|\widehat{\varphi}\|^*_{n,1} \|\zeta\|^*_{n,\infty} \le \frac{2\sigma\rho^2}{\sqrt{m+1}} \sigma\sqrt{2 \ln \left[3(n+1)/\alpha\right]},
\end{align*}
where the last inequality, holding with probability $\ge 1-\alpha/3$, is due to \eqref{eq:complex_gaussian_max}.
\par
Finally, observe that, with probability $\geq 1-\alpha/3$ (cf. \eqref{eq:varphizeta}),
\[
\|\phi^o * \zeta\|_{n,2}\leq \sqrt{2}\rho\left(1+\sqrt{\ln[3/\alpha]}\right).
\]
Therefore, we have for $\delta^{(4)}$:
\[
\delta^{(4)} \le \sigma\|\widehat \varphi\|_{n,2}\|\phi^o * \zeta\|_{n,2} \le
\sigma\frac{2\rho^2}{\sqrt{m+1}} \sqrt{2}\rho\left(1+\sqrt{\ln[3/\alpha]}\right)=
\sigma\frac{2\sqrt{2}\rho^3}{\sqrt{m+1}}\left(1+\sqrt{\ln[3/\alpha]}\right)
\]
with probability $1-\alpha/3$. When substituting the bound for $\delta^{(k)},\,k=1,...,4$, into \eqref{eq:l2point-long} we arrive at the result of the theorem.\qed

\section{Miscellaneous proofs}\label{sec:prdisc}
\paragraph{Proof of relations \eqref{eq:l1f} and \eqref{eq:var2f}}
Let $n=2m$, $\phi\in \C(\Z_0^m)$, and let $\varphi\in \C(\Z_{0}^n)$ satisfy $\varphi=\phi*\phi$. Then
 \begin{eqnarray*}
 \|\varphi\|^*_{n,1} &=& (2m+1)^{-1/2}\sum_{k=0}^{2m} |[F_{2m}\varphi_{0}^{2m}]_k| =  \sqrt{2m+1}\sum_{k=0}^{2m}  \left({|[F_{2m}\phi_{0}^{2m}]_k|\over{\sqrt{2m+1}}}\right)^2\nn
&=&  \sqrt{2m+1} \|\phi\|^{*2}_{2m,2} =   \sqrt{2m+1}\|\phi\|^2_{2m,2}=   \sqrt{2m+1} \|\phi\|_{m,2}^2
\leq {\sqrt{2m+1} \, \rho^2\over m+1},
 \end{eqnarray*}
implying \eqref{eq:l1f}. Moreover, since $1-\phi*\phi=(1+\phi)*(1-\phi)$, for all $x\in \C(\Z)$ one has for all ${\tau}\in \Z$:
\begin{eqnarray*}
|x_{{\tau}} - [\varphi^o* x]_{{\tau}}| &=& |[(1 +\phi^o)*(1 - \phi^o)*x]_{\tau}| = \left|\sum_{j=0}^m
 [1 + \phi^o]_{j} \,[x - \phi^o*x]_{\tau - j}\right|\\
&\le& \|1 + \phi^o\|_1 \max_{0\leq j\leq m} |[x - \phi^o* x]_{{\tau - j}}|\leq {\sigma(1+\rho)\rho\over \sqrt{m+1}}
\end{eqnarray*}
(we have used \eqref{eq:ss1} to obtain the last inequality), and
\[\|x - [\varphi^o* x]\|_{n,2}\leq \sigma(1+\rho)\rho.
\] Next note that
\[
\|\varphi^o*\zeta\|_{n,2}^2=\langle \zeta,M(\varphi^o)\zeta\rangle_n,
\]
where $M(\varphi)$ is defined as in \eqref{M(phi)}. When taking into account that
\[
\|\varphi\|_2\leq \|\varphi\|^*_{n,1}\leq{2\rho^2\sqrt{n+1}\over n+2},
\]
we get (cf. \eqref{tr22})
$\|M(\varphi^o)\|_F^2=(n+1)\|\varphi\|_2^2\leq 4\rho^4$,
so that the concentration inequality \eqref{trq1} now implies that, given $0<\alpha\leq 1$, with probability at least $1-\alpha$,
\begin{equation}\label{eq:varphizeta}
\|\varphi^o*\zeta\|_{n,2}\leq 2\sqrt{2}\rho^2\left(1+\sqrt{\ln[1/\alpha]}\right),
\end{equation}
and we arrive at \eqref{eq:var2f}.

\paragraph{Proof of Lemma \ref{lemma:shift_invariant}}
As a precursory remark, note that if a finite-dimensional subspace $\S$ is shift-invariant, i.e. $\Delta \S \subseteq \S$, then necessarily $\Delta$ is a bijection on $\S$, and $\Delta\S = \S$. Indeed, when restricted on $\S$, $\Delta$ obviously is a linear transformation with a trivial kernel, and hence a bijection.

$1^\circ.$ To prove the direct statement, note that the solution set of \eqref{eq:diff_eq} with $\deg(p(\cdot)) = s$ is a shift-invariant subspace of $\C(\Z)$ -- let us call it $\S'$. Indeed, if $x \in \C(\Z)$ satisfies $\eqref{eq:diff_eq}$, so does $\Delta x$, so $\S'$ is shift-invariant. To see that $\dim(\S') = s$, note that $x \mapsto x_1^s$ is a bijection $\S' \to \C^s$: under this map arbitrary $x_1^s \in \C^{s}$ has a unique preimage. Indeed, as soon as one fixes $x_1^s$, \eqref{eq:diff_eq} uniquely defines the next samples $x_{s+1}, x_{s+2}, ...$ (note that $p(0) \ne 0$); dividing \eqref{eq:diff_eq} by $\Delta^s$, one can retrieve the remaining samples of $x$ since $\deg(p(\cdot))=s$ (here we used that $\Delta$ is bijective on $\S$).

$2^\circ.$ For the converse statement, first note that any polynomial $p(\cdot)$ with $\deg(p(\cdot)) = s$, $p(0)=1$ is uniquely expressed as $p(z) = \prod_{k=1}^s (1-z/z_k)$ where $z_1, ..., z_s$ are its roots. Since $\S$ is shift-invariant, we have $\Delta \S = \S$ as discussed above, i.e. $\Delta$ is a bijective linear transformation on $\S$. Let us fix some basis $\boldsymbol{e}^T = \begin{bmatrix} e^1,\, ...,\, e^s\end{bmatrix}$ on $\S$, and denote $A$ the $s \times s$ matrix of $\Delta$ in it, i.e. $\Delta(e^j) = \sum_{i=1}^s a_{ij} e^{i}$. The basis $\boldsymbol{e}$ might be chosen such that $A$ is upper-triangular (say, by passing to its Jordan normal form). Then, any vector $x \in \S$ satisfies
$
q(\Delta) x \equiv 0,$
where $q(\Delta) = \prod_{i=1}^s (\Delta - a_{ii}) = \det(\Delta I - A)$ is the characteristic polynomial of $A$.
Note that $\prod_{i=1}^s a_{ii} = \det A \ne 0$ since $\Delta$ is a bijection on $\S$. Hence, introducing $p(\Delta) = q(\Delta) / \det A$, we obtain
$\prod_{i=1}^s (1-\Delta c_i) x \equiv 0$ for some complex numbers $c_i \ne 0$. This means that $\S$ is contained in the solution set $\S'$ of \eqref{eq:diff_eq} with $\deg(p(\cdot))=s$, $p(0)=1$, which is itself a shift-invariant subspace of dimension $s$ (by the direct statement). Since $\dim(\S) = \dim(\S') = s$, the two subspaces coincide.
The uniqueness of $p(\cdot)$ is implied by the fact that $q(\cdot)$ is a characteristic polynomial of $A$.
\qed

\paragraph{Proof of \eqref{eq:interps}} Assume that Assumption A holds true for some $n\geq s$ and $\varepsilon\equiv 0$. Let $\Pi_{\S_n}$ be the Euclidean projector
on the space $\S_n$ of elements of $\mathcal{S}$ restricted on $\C(\Z_0^n)$. Since $\dim(\S_n)\leq s$,
$\|\Pi_{\S_n}\|_2^2=\Tr(\Pi_{\S_n})\leq s$, there is $\iota\in \{0,...,n\}$ such that the $\iota$-th column $[\Pi_{\S_n}]_\iota$ of  $\Pi_{\S_n}$ satisfies $\|[\Pi_{\S_n}]_\iota\|_2\leq \sqrt{s/(n+1)}$. Note that one has $x_\iota-\langle [\Pi_{\S_n}]_\iota, x_0^n\rangle=0$, and $\S_n$ is time-invariant, implying that
\[
x_\tau-\langle [\Pi_{\S_n}]_\iota, x_{\tau-\iota}^{\tau-\iota+n}\rangle=0, \;\;\forall \tau \in\Z.
\]
We conclude that there is $\phi^o\in \C(\Z_{-n}^n)$,
$\phi^o=[0;...;0;[\Pi_{\S_n}]_\iota; 0;...;0]$ (i.e., vector $[\Pi_{\S_n}]_\iota$ completed with zeros in such a way that $\iota$-th element of $[\Pi_{\S_n}]_\iota$ becomes the central ($n+1-$th) entry of $\phi^o$) such that
\[
\|\phi^o\|_2\leq \sqrt{s/(n+1)},\;\;\mbox{and}\;\;x_\tau-[\phi^o*x]_\tau=0, \;\;\forall \tau \in\Z.\eqno{\square}
\]
\subsection{Proof of Theorem \ref{th:sines}}\label{sec:prsines}
Note that to prove the theorem we have to exhibit a vector $q\in \C^{m+1}$ of small $\ell_2$-norm and such that the polynomial
$1-q(z)=1-\left[\sum_{i=0}^m q_iz^i\right]$ is divisible by $p(z)$, i.e., that there is a polynomial $r(z)$ of degree $m-s$ such that
\[
1-q(z)=r(z)p(z).
\]
Indeed, this would imply that
\[
x_t-[q*x]_t=[1-q(\Delta)]x_t=r(\Delta)p(\Delta)x_t=0
\]
due to $p(\Delta)x_t=0$,
\par
The bound
$
\|q\|_2\leq C's^{3/2}\sqrt{\ln s\over m}
$
of \eqref{propofq} is proved in \cite[Lemma 6.1]{jn-2014}. Our objective is to prove the ``remaining'' inequality
\[
\|q\|_2\leq C's\sqrt{\ln[ms]\over m}.
\]
So, let $\theta_1,...,\theta_s$ be complex numbers of modulus 1 -- the roots of the polynomial $p(z)$. Given $\delta=1-\epsilon\in(0,1)$, let us set $\bar{\delta}={2\delta/(1+\delta)}$, so that
\begin{equation}\label{eq1}
{\bar{\delta}\over\delta}-1=1-\bar{\delta}>0.
\end{equation}
Consider the function
\[%begin{equation}\label{eq2}
\bar{q}(z)=\prod_{i=1}^s{z-\theta_i\over\delta z-\theta_i}.
\]%end{equation}
\def\B{{\cal B}}
Note that $\bar{q}(\cdot)$ has no singularities in the circle
\[
\B=\{z:|z|\leq {1/\bar{\delta}}\};
\]
besides this, we have
$%\begin{equation}\label{eq4}
\bar{q}(0)=1.
$%e\nd{equation}
Let $|z|=1/\bar{\delta}$, so that $z=\bar{\delta}^{-1}w$ with $|w|=1$. We have
$$
{|z-\theta_i|\over|\delta z-\theta_i|}={1\over\delta}{|w-\bar{\delta}\theta_i|\over |w-{\bar{\delta}\over \delta}\theta_i|}.
$$
 We claim that when $|w|=1$, $|w-\bar{\delta}\theta_i|\leq |w-{\bar{\delta}\over\delta}\theta_i|$.
\begin{quote}
Indeed, assuming w.l.o.g. that $w$ is not proportional to $\theta_i$, consider triangle $\Delta$ with the vertices $A=w$, $B=\bar{\delta}\theta_i$ and $C={\bar{\delta}\over\delta}\theta_i$. Let also $D=\theta_i$. By (\ref{eq1}), the segment $\overline{AD}$ is a median in $\Delta$, and $\angle CDA $ is $\geq{\pi\over 2}$ (since $D$ is the closest to $C$ point in the unit circle, and the latter contains $A$), so that  $|w-\bar{\delta}\theta_i|\leq |w-{\bar{\delta}\over\delta}\theta_i|$.
\end{quote}
As a consequence, we get
\begin{equation}\label{eq10}
z\in \B\; \Rightarrow\; |\bar{q}(z)|\leq\delta^{-s},
\end{equation}
whence also
\begin{equation}\label{eq5}
|z|=1\;\Rightarrow\; |\bar{q}(z)|\leq\delta^{-s}.
\end{equation}
Now, the polynomial $p(z)=\prod_{i=1}^s(z-\theta_i)$ on the boundary of $\B$ clearly satisfies
$$
|p(z)|\geq \left[{1\over \bar{\delta}}-1\right]^s=\left[{1-\delta\over 2\delta}\right]^s,
$$
which combines with (\ref{eq10}) to imply that the modulus of the holomorphic in $\B$ function
$$
\bar{r}(z)=\left[\prod_{i=1}^s(\delta z-\theta_i)\right]^{-1}
$$
is bounded with $\delta^{-s}\left[{1-\delta\over 2\delta}\right]^{-s}=\left[{2\over 1-\delta}\right]^s$ on the boundary of $\B$. It follows that the coefficients $r_j$ of the Taylor series of $\bar{r}$ satisfy
$$
|r_j|\leq  \left[{2\over 1-\delta}\right]^s\bar{\delta}^{j},\;\;j=0,1,2,...
$$
When setting
\begin{equation}\label{eq:qldiv}
q^\ell(z)=p(z)r^\ell(z),\;\;r^\ell(z)=\sum_{j=1}^\ell r_j z^j,
\end{equation}
for $|z|\leq 1$, utilizing the trivial upper bound $|p(z)|\leq 2^s$, we get
\begin{eqnarray}\label{eq12}
|q^\ell(z)-\bar{q}(z)|&=&|p(z)|[r^\ell(z)-\bar{r}(z)]|\leq |p(z)||r^\ell(z)-\bar{r}(z)|\nn
&\leq&2^s\left[{2\over 1-\delta}\right]^s \sum_{j=\ell+1}^\infty |r_j|\leq \left[{4\over 1-\delta}\right]^s{\bar{\delta}^{\ell+1}\over 1-\bar{\delta}}.
\end{eqnarray}
Note that $q^\ell(0)=p(0)r^\ell(0)=p(0)\bar{r}(0)=1$, that $q^\ell$ is a polynomial of degree $\ell+s$, and that $q^\ell$ is divisible by $p(z)$.  Besides this, on the unit circumference we have, by (\ref{eq12}),
\begin{equation}\label{eq13}
|q^\ell(z)|\leq |\bar{q}(z)|+\left[{4\over 1-\delta}\right]^s{\bar{\delta}^{\ell+1}\over 1-\bar{\delta}}\leq \delta^{-s}+\underbrace{\left[{4\over 1-\delta}\right]^d{\bar{\delta}^{\ell+1}\over 1-\bar{\delta}}}_{R}
\end{equation}
(we have used (\ref{eq5})).
Now,
$$
\bar{\delta}={2\delta\over1+\delta}={2-2\epsilon\over 2-\epsilon}={1-\epsilon\over 1-\epsilon/2}\leq 1-\epsilon/2\leq \e^{-\epsilon/2},
$$
and
$$
{1\over 1-\bar{\delta}}={1+\delta\over 1-\delta}={2-\epsilon\over \epsilon}\leq {2\over\epsilon}.
$$
We can upper-bound $R$:
\[
R=\left[{4\over 1-\delta}\right]^s{\bar{\delta}^{\ell+1}\over 1-\bar{\delta}}
\leq {2^{2s+1}\over\epsilon^{s+1}}\e^{-\epsilon \ell/2}
\]Now, given positive integer $\ell$ and positive $\alpha$ such that
\begin{equation}\label{alphasmall}
{\alpha\over \ell}\leq\four,
\end{equation}
let $\epsilon={\alpha\over 2\ell s}$.
Since $0<\epsilon\leq {1\over 8}$, we have $-\ln(\delta)=-\ln(1-\epsilon)\leq 2\epsilon={\alpha\over \ell s},$ implying that $\bar{\delta}\leq \e^{-\epsilon/2}=\e^{-{\alpha\over 4\ell s}}$, and
\[
R\leq \left[{8\ell s\over\alpha}\right]^{s+1}\exp\{-{\alpha\over 4s}\}.
\]
Now let us put
\[%begin{equation}\label{alpha}
\alpha=\alpha(\ell,s)=4s(s+2)\ln(8\ell s);
\]%end{equation}
observe that this choice of $\alpha$ satisfies (\ref{alphasmall}), provided that
\begin{equation}\label{nislarge}
\ell\geq O(1)s^2\ln(s+1)
\end{equation}
with properly selected absolute constant $O(1)$. With this selection of $\alpha$, we have $\alpha\geq1$, whence
\begin{eqnarray*}
R\left[{\alpha\over \ell}\right]^{-1}&\leq& \exp\left\{-{\alpha\over 4s}\right\}\left[{8\ell s\over\alpha}\right]^{s+1}{\ell\over \alpha}\leq\exp\left\{-{\alpha\over 4s}\right\}[8\ell s]^{s+2}\\
&\leq& \exp\{-(s+2)\ln(8\ell s)\}\exp\{(s+2)\ln(8\ell s)\}=1,
\end{eqnarray*}
that is,
\begin{equation}\label{RR}
R\leq {\alpha\over \ell}\leq {\four}.
\end{equation}
Furthermore,
\begin{equation}\label{firstterm}
\begin{array}{rcll}
\delta^{-s}&=&\exp\{-s\ln(1-\epsilon)\}\leq \exp\{2\epsilon s\}=\exp\{{\alpha\over \ell}\}\leq 2,\\
\delta^{-2s}&=&\exp\{-2s\ln(1-\epsilon)\}\leq \exp\{4\epsilon s\}=\exp\{{2\alpha\over \ell }\}\leq 1+\exp\{{1\over 2}\}{2\alpha\over \ell}\leq 1+{4\alpha\over \ell}.
\end{array}
\end{equation}
When invoking (\ref{eq13}) and utilizing (\ref{firstterm}) and (\ref{RR}) we get
$$
{1\over 2\pi}\oint_{|z|=1}|q^\ell(z)|^2|dz|\leq \delta^{-2s}+2\delta^{-s}R+R^2\leq
1+4{\alpha\over \ell}+4R+{1\over 4}R\leq 1+10{\alpha\over \ell}.
$$
On the other hand, denoting by $q_0$, $q_1$,...,$q_{\ell+s}$ the coefficients of the polynomial $q^\ell$ and taking into account that $\bar{q}_0=q^\ell(0)=1$, we have
\begin{equation}\label{eq:qnorm}
1+\sum_{i=1}^{\ell+s}|q_i|^2=|q_0|^2+...+|q_{\ell+s}|^2={1\over 2\pi}\oint_{|z|=1}|q^\ell(z)|^2|dz|\leq 1+10{\alpha\over \ell}.
\end{equation}
We are done: when denoting $m=\ell+s$, and $q(z)=\sum_{i=1}^m q_j z^j$, we have the vector of coefficients $q=[0;q_1;...;q_m]\in \C^{m+1}$ of $q(z)$ such that, by \eqref{eq:qnorm},
\[
\|q\|^2_2\leq {40s(s+2)\ln[8s(m-s)]\over m-s},
\]
and such that the polynomial $q^\ell(z)=1+q(z)$ is divisible by $p(z)$ due to \eqref{eq:qldiv}.\qed

\section{Additional numerical illustrations}\label{sec:demoexp}
In these demonstration experiments, we compare the penalized $\ell_2$-recovery of Sec.~\ref{sec:open} to the Lasso as in  \cite{recht1}. We use the same setting of the penalization parameters as in the Monte-Carlo experiments Sec. \ref{sec:exps} with regularization parameters set to the theoretically recommended value~\cite{recht1}.

\begin{figure}[b!]
\center
\begin{minipage}{0.24\textwidth}
\centering
\includegraphics[width=1\textwidth, height=0.135\textheight, clip=true, angle=0]{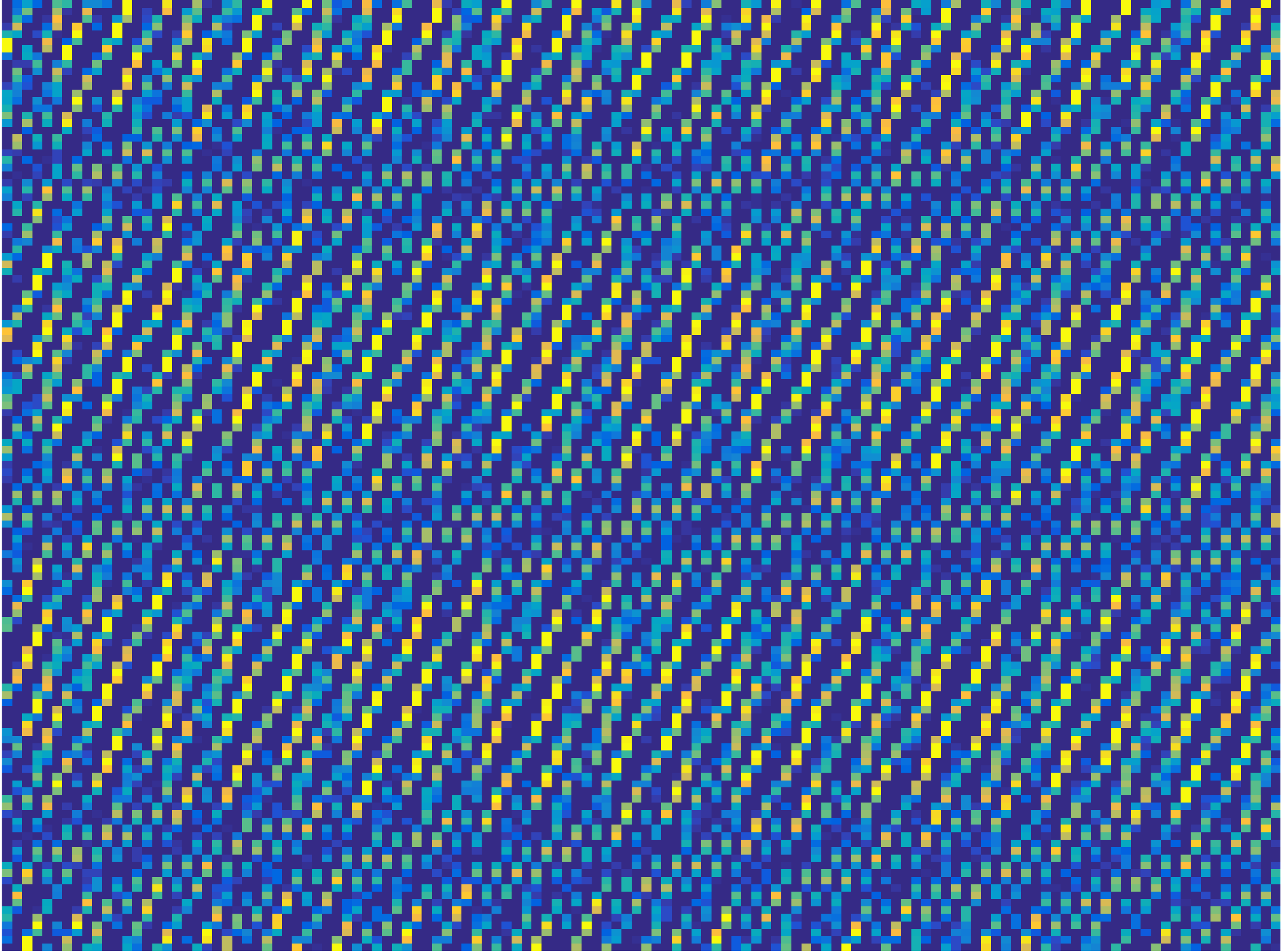}\vspace{0.1cm}
\includegraphics[width=1\textwidth, height=0.135\textheight, clip=true, angle=0]{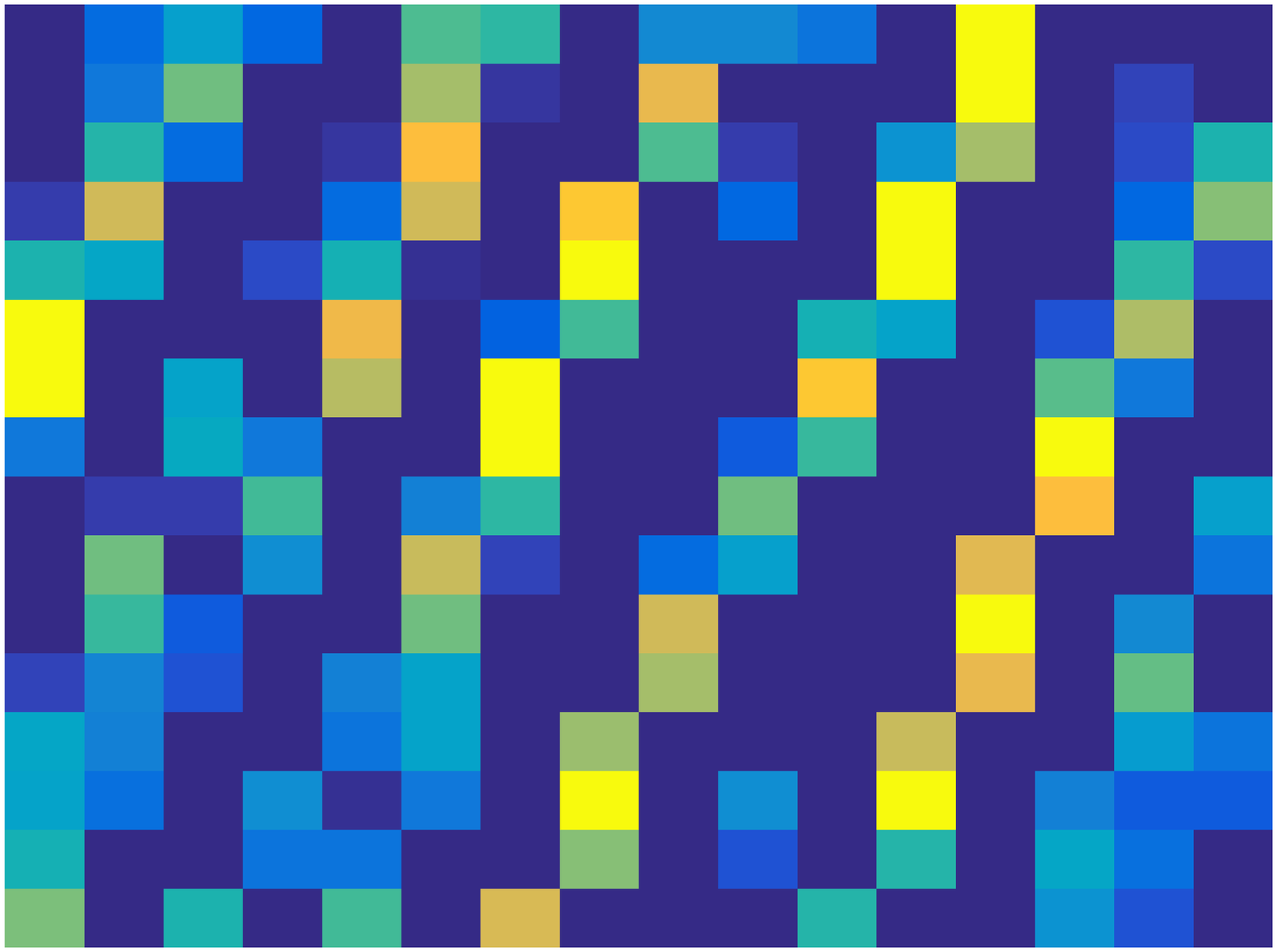}
True signal
\end{minipage}
\begin{minipage}{0.24\textwidth}
\centering
\vspace{-0.07cm}
\includegraphics[width=1\textwidth, height=0.135\textheight, clip=true, angle=0]{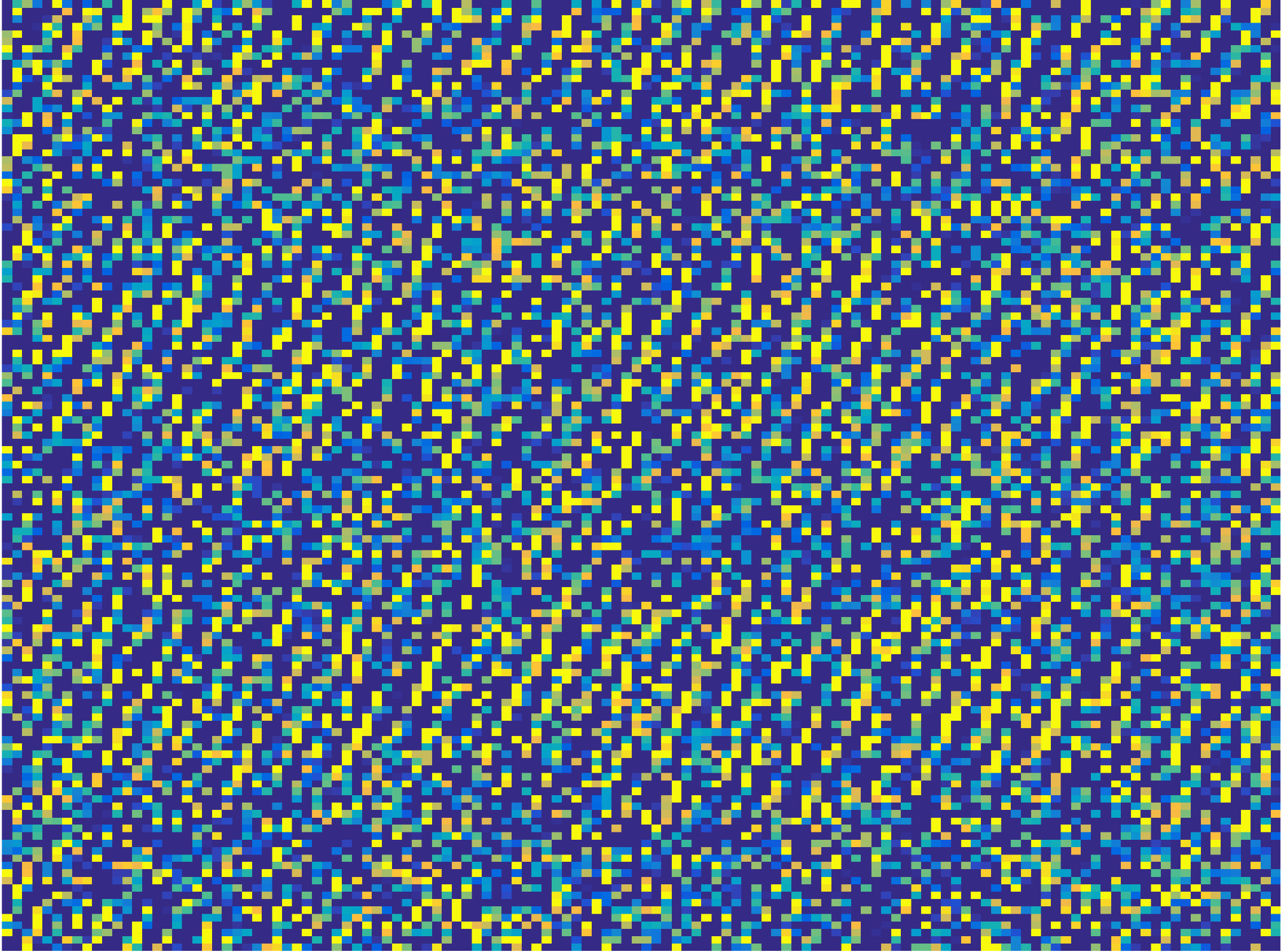}\vspace{0.1cm}
\includegraphics[width=1\textwidth, height=0.135\textheight, clip=true, angle=0]{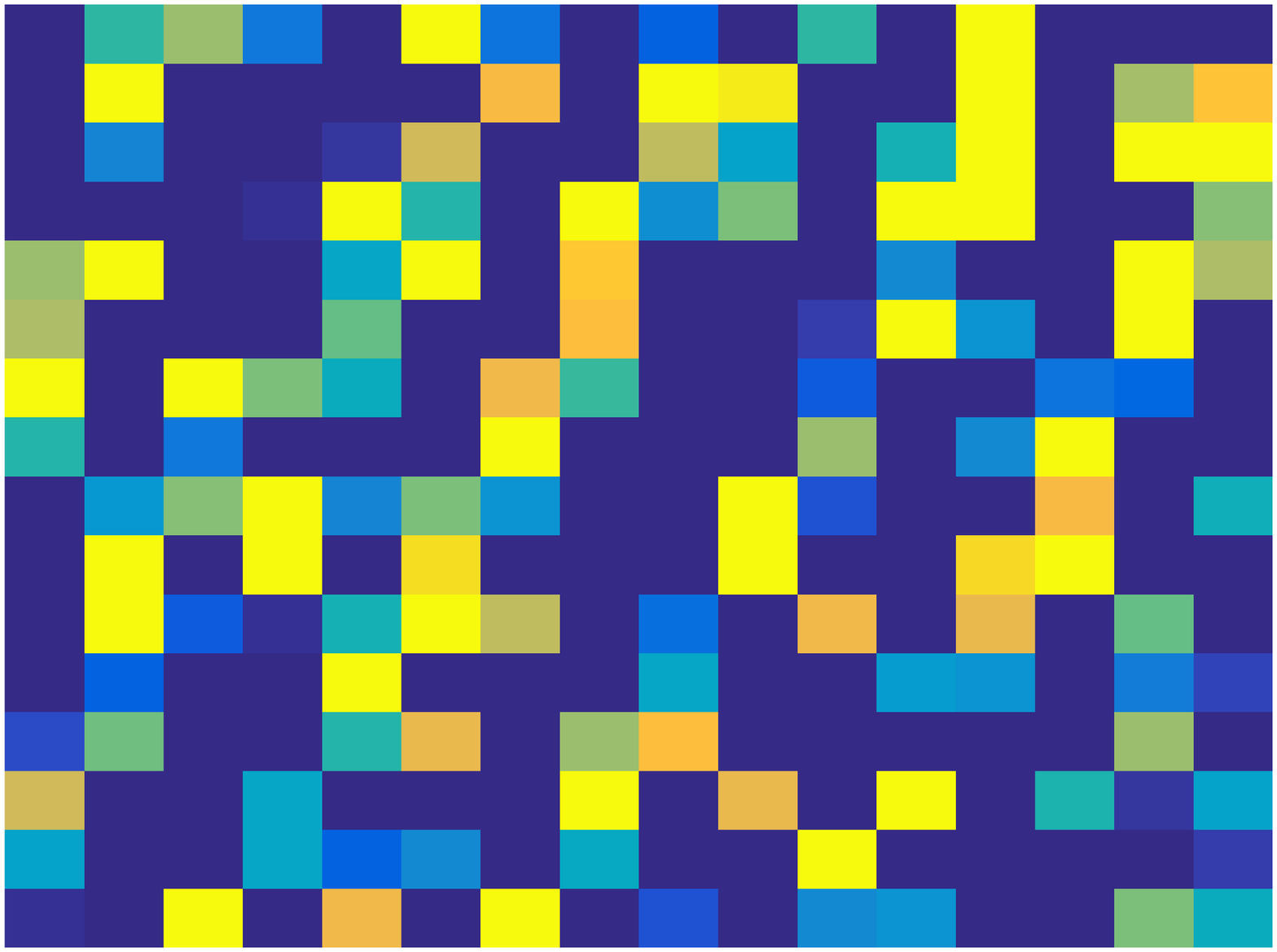}
Observations
\end{minipage}
\begin{minipage}{0.24\textwidth}
\centering
\includegraphics[width=1\textwidth, height=0.135\textheight, clip=true, angle=0]{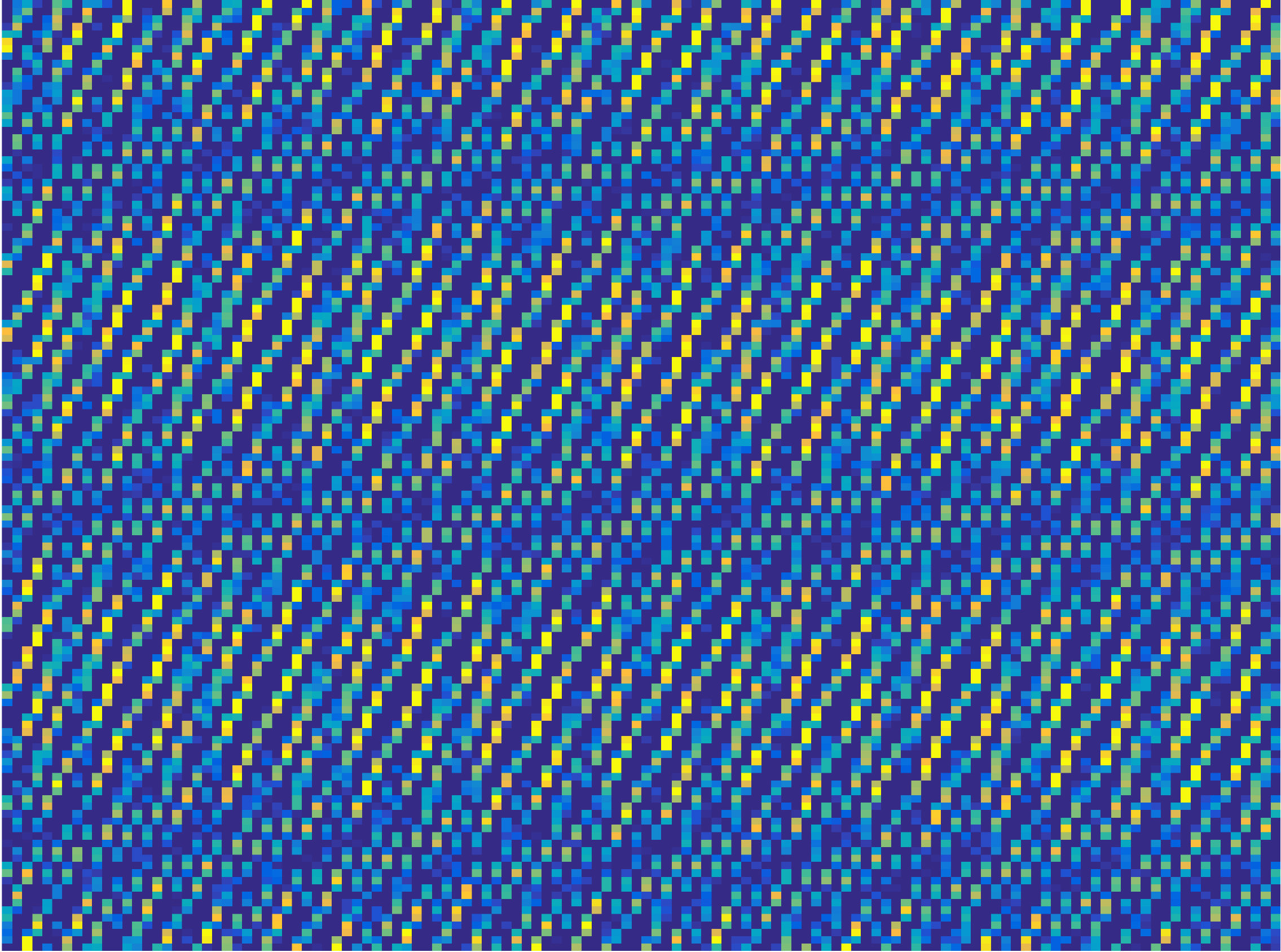}\vspace{0.1cm}
\includegraphics[width=1\textwidth, height=0.135\textheight, clip=true, angle=0]{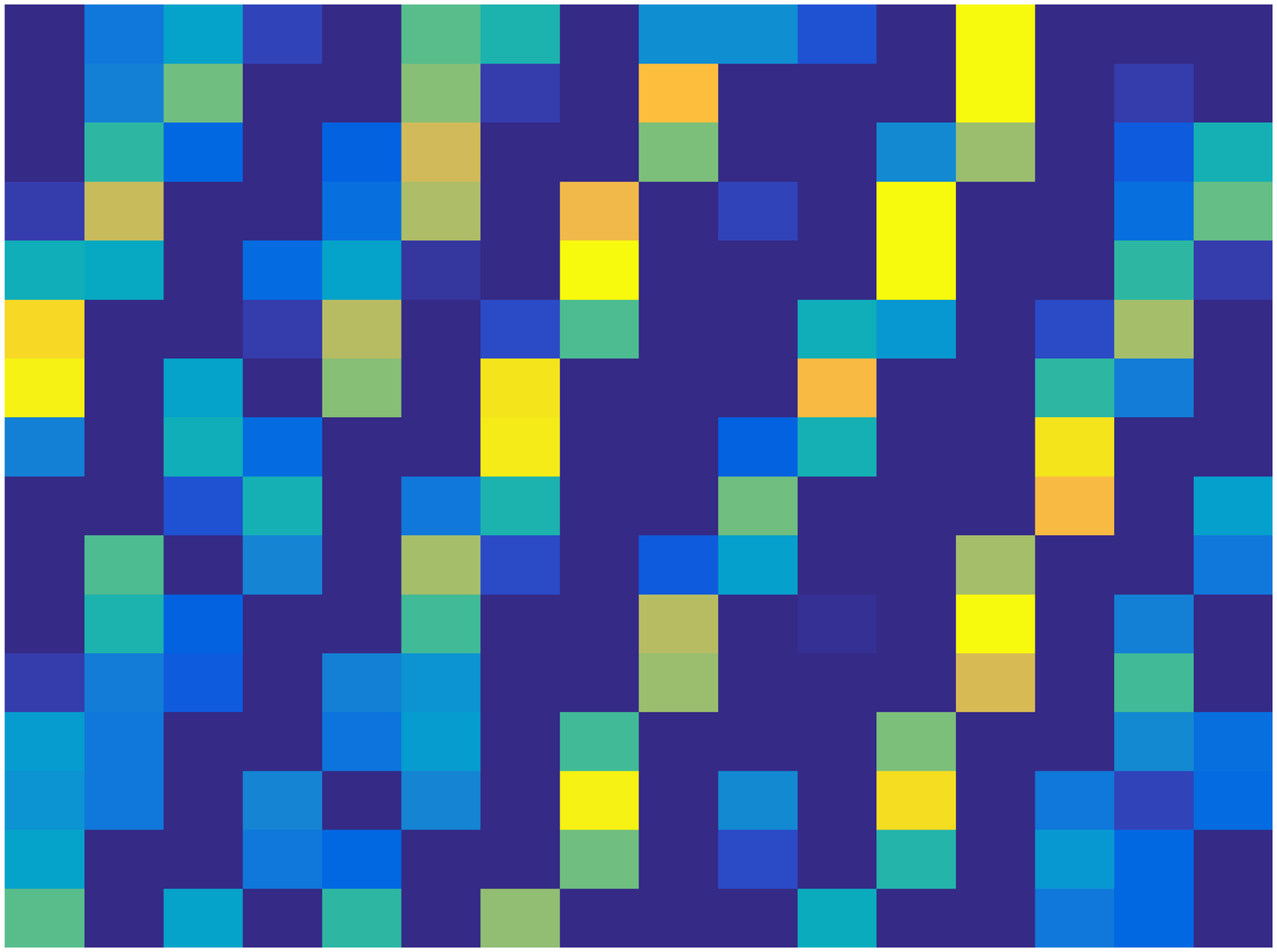}
Penalized recovery
\end{minipage}
\begin{minipage}{0.24\textwidth}
\centering
\vspace{-0.07cm}
\includegraphics[width=1\textwidth, height=0.135\textheight, clip=true, angle=0]{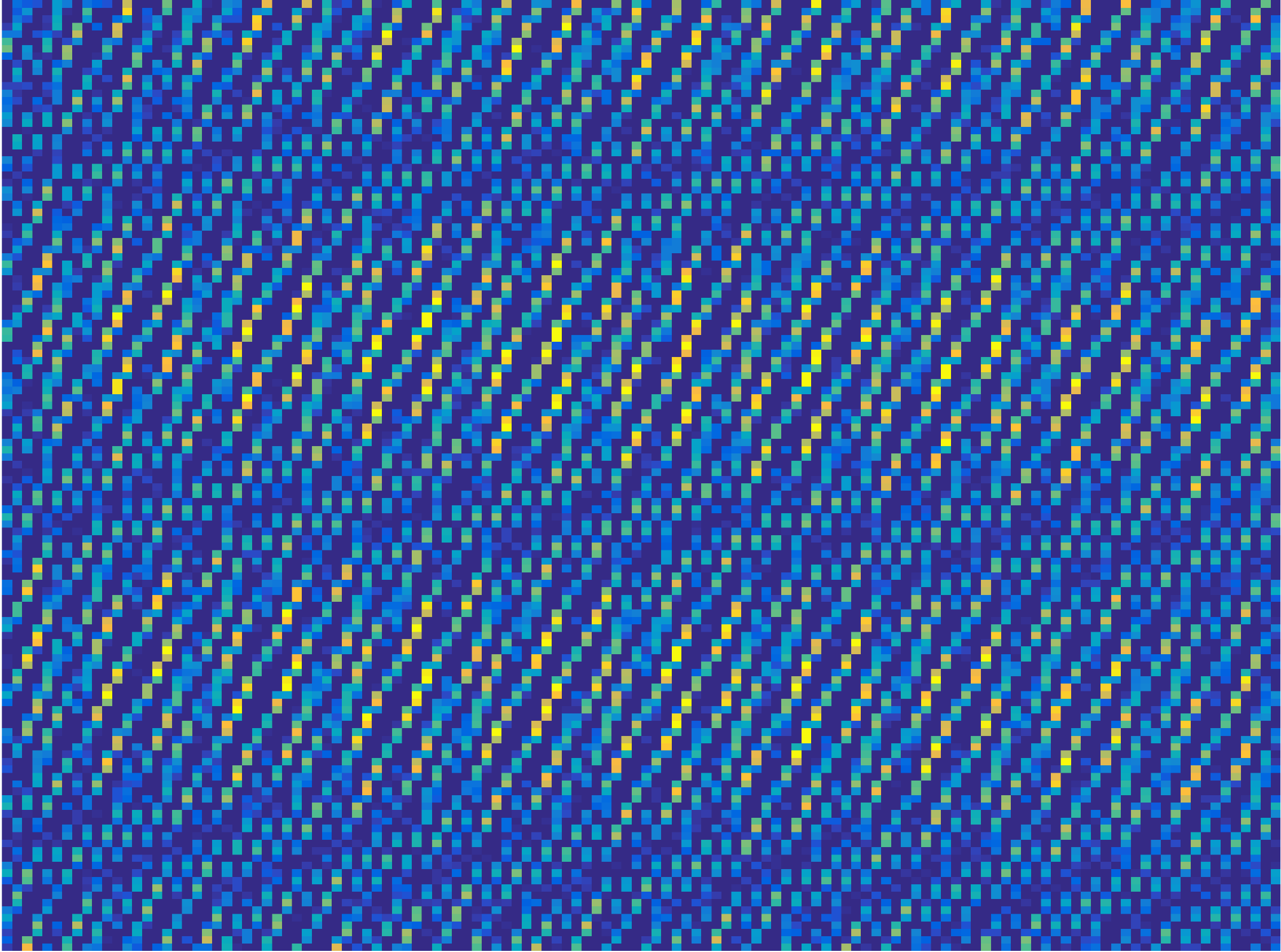}\vspace{0.1cm}
\includegraphics[width=1\textwidth, height=0.135\textheight, clip=true, angle=0]{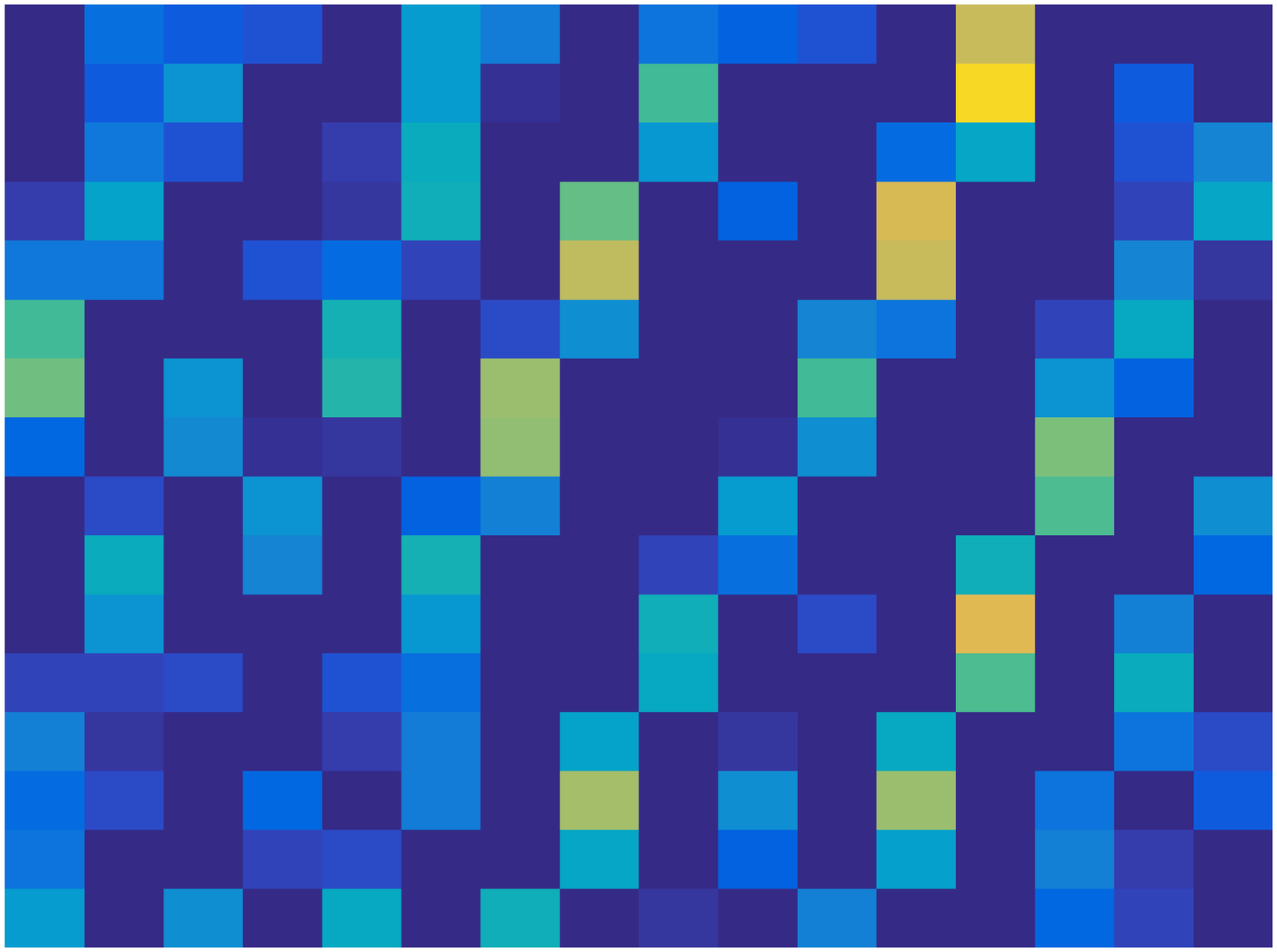}
Lasso
\end{minipage}
\caption{Recovery of a sum of $4$ harmonic oscillations observed with $\text{SNR}= {0.5}$, i.e. $\|y-x\|_2 \approx 2$. 2nd row:  magnified upper left corner of the image.
$\ell_2$-error: 0.18 for the penalized $\ell_2$-recovery, 0.45 for the Lasso \cite{recht1}.}
\label{fig:sines2}
\end{figure}

\paragraph{2-D harmonic oscillations}
In this experiment (see Fig.~\ref{fig:sines2}) we recover a sum of $4$ harmonic oscillations in $\R^2$ with random frequencies, observed with $\text{SNR} = 0.5$ (the signal is normalized in the $\ell_2$-norm).

\begin{figure}[b!]
\center
\begin{minipage}{0.24\textwidth}
\centering
\includegraphics[width=1\textwidth, height=0.135\textheight, clip=true, angle=0]{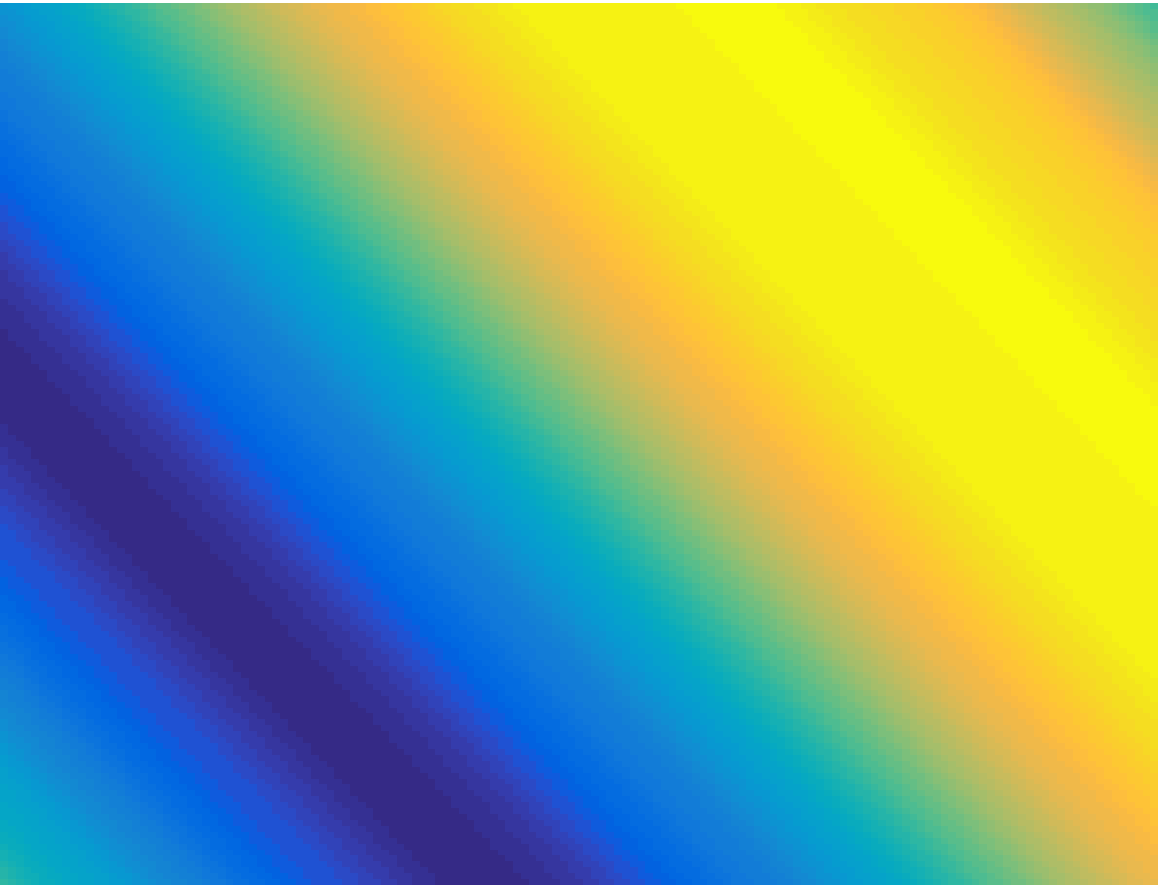}\vspace{0.1cm}
\includegraphics[width=1\textwidth, height=0.135\textheight, clip=true, angle=0]{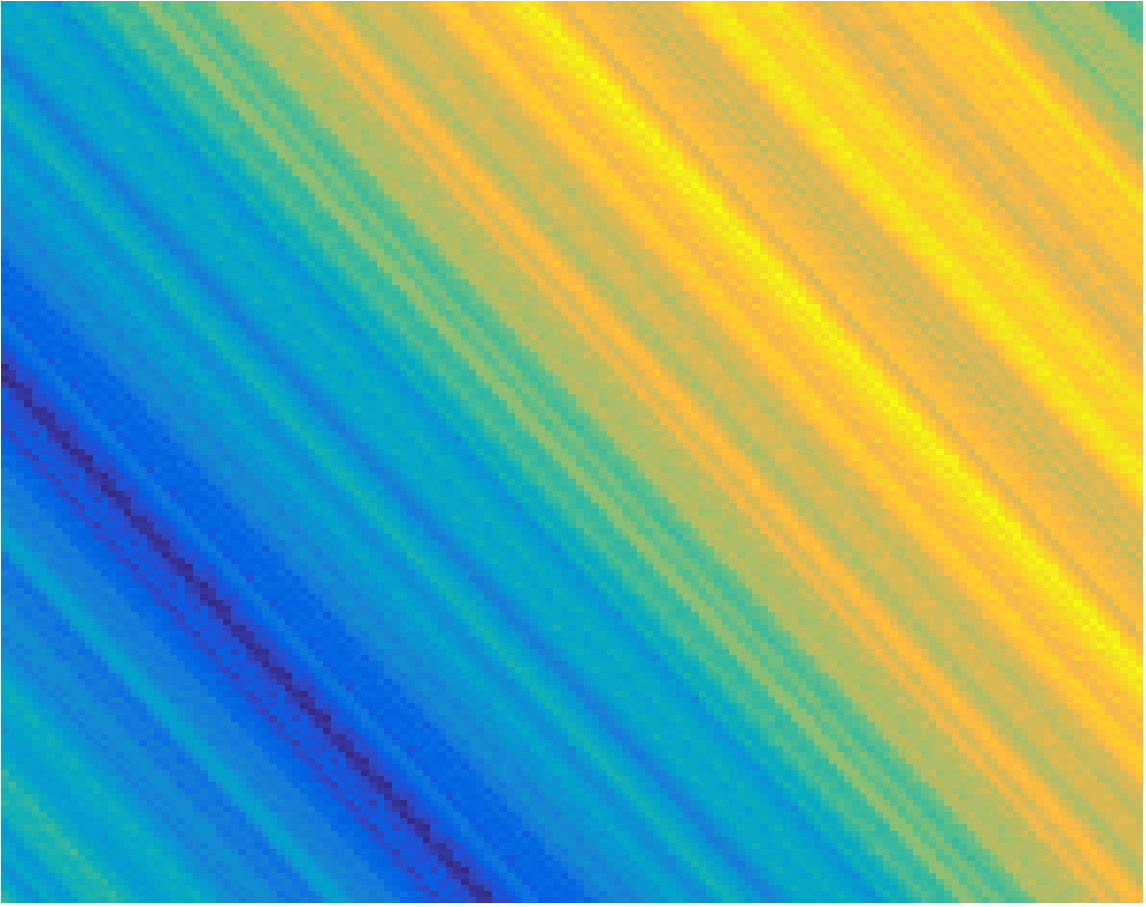}
True signal
\end{minipage}
\begin{minipage}{0.24\textwidth}
\centering
\vspace{-0.07cm}
\includegraphics[width=1\textwidth, height=0.135\textheight, clip=true, angle=0]{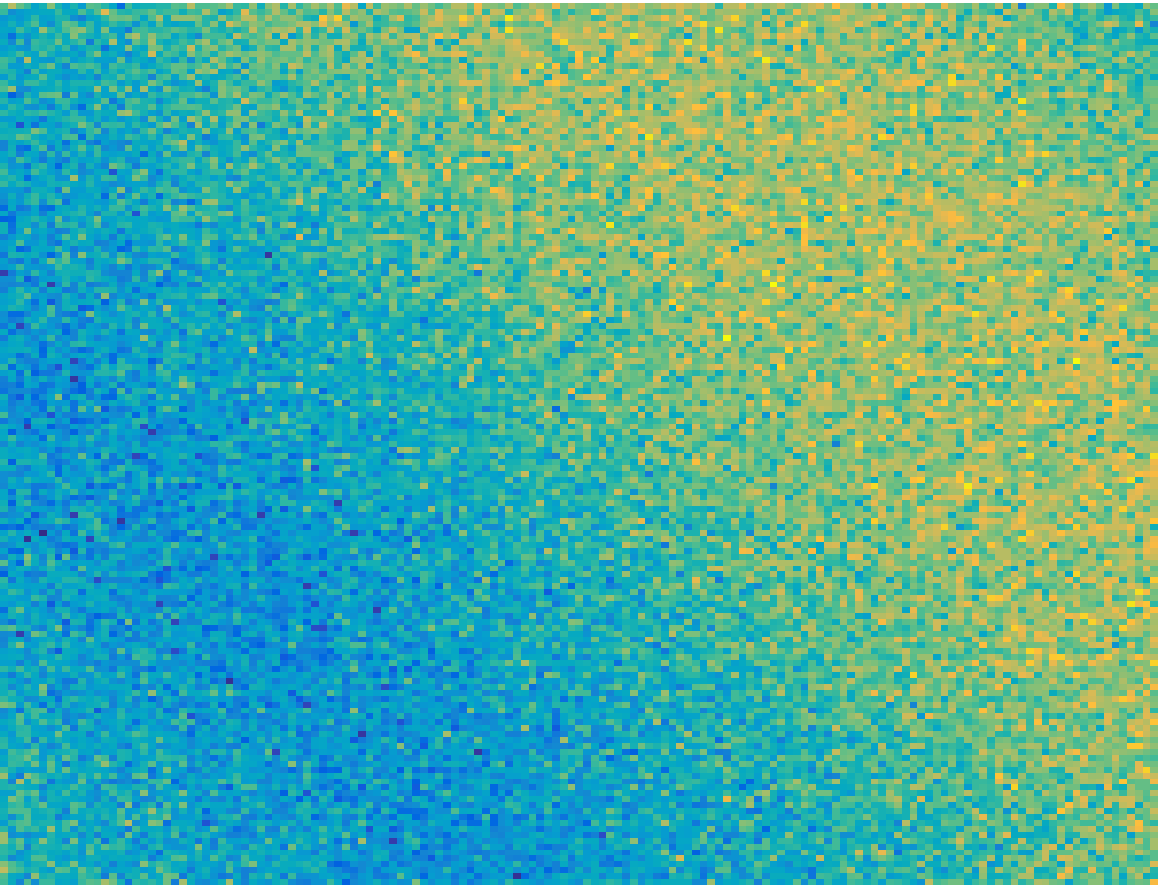}\vspace{0.1cm}
\includegraphics[width=1\textwidth, height=0.135\textheight, clip=true, angle=0]{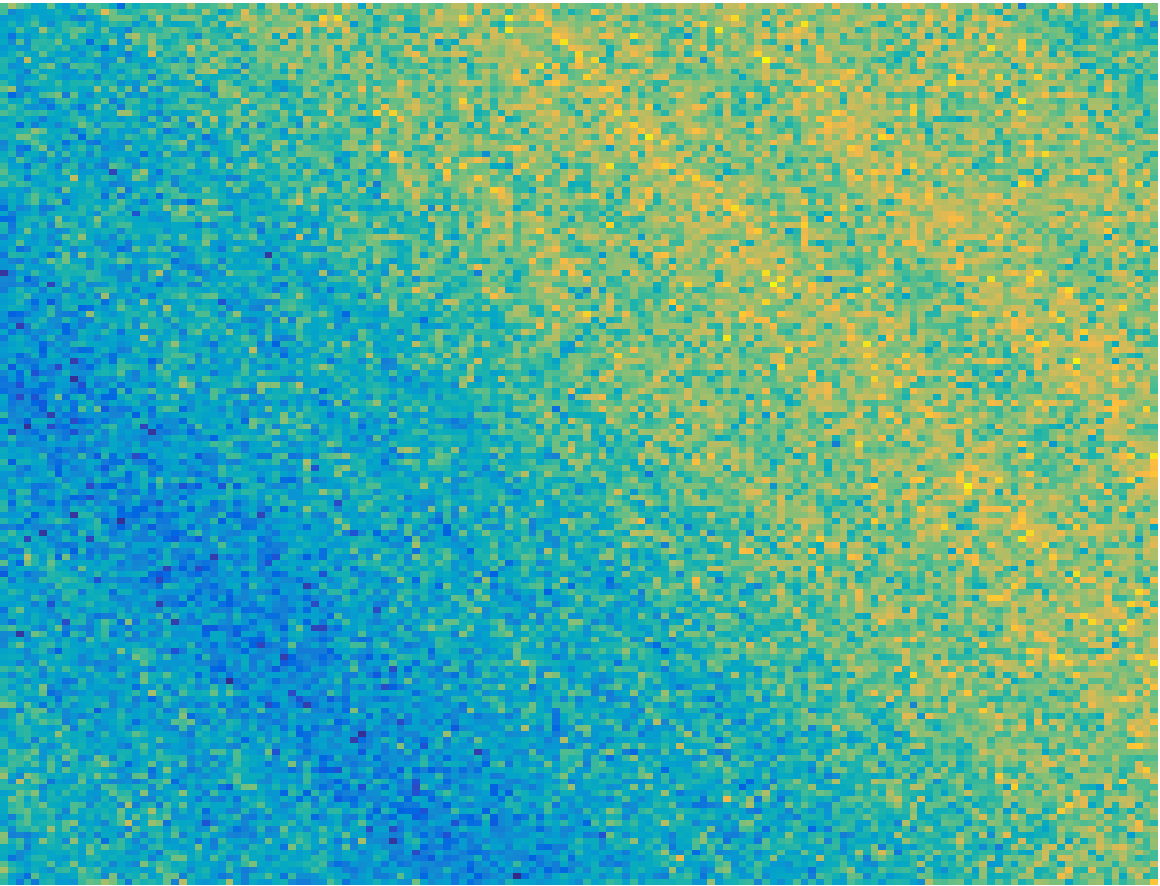}
Observations
\end{minipage}
\begin{minipage}{0.24\textwidth}
\centering
\includegraphics[width=1\textwidth, height=0.135\textheight, clip=true, angle=0]{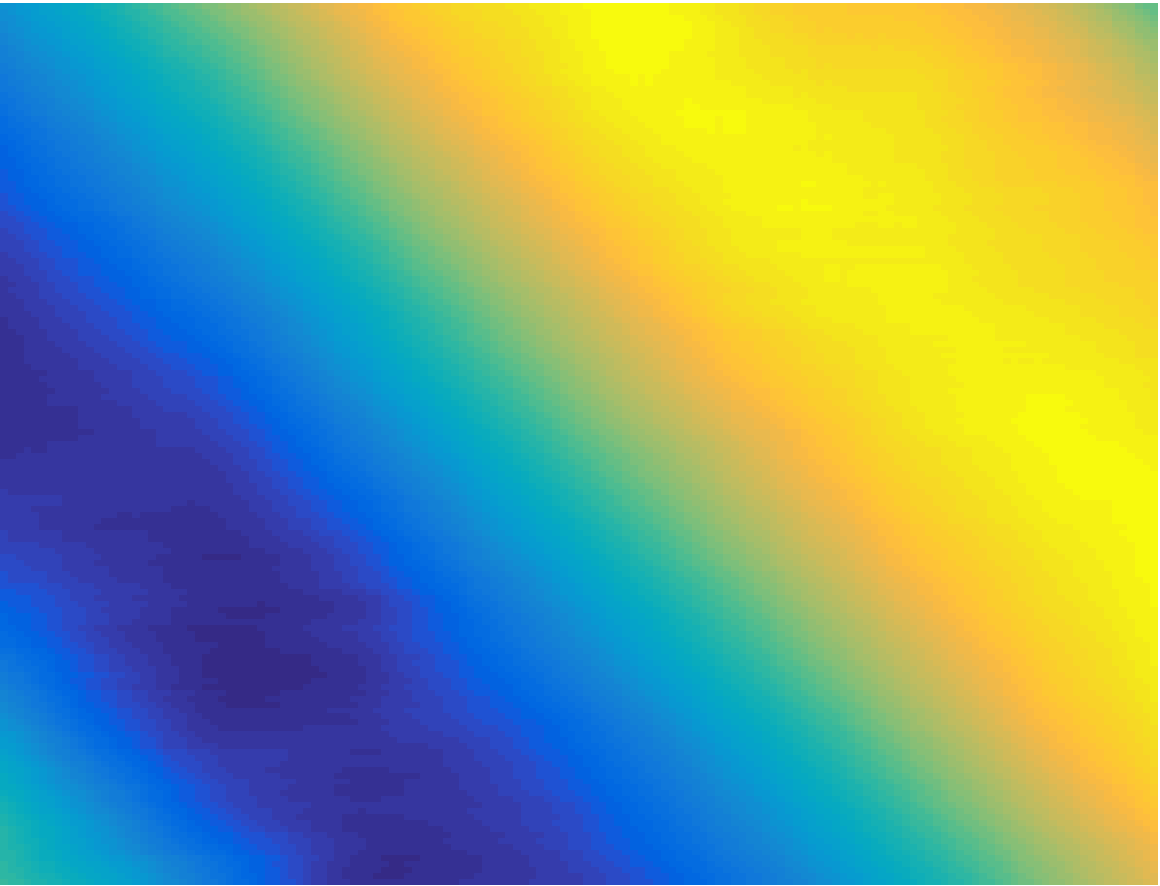}\vspace{0.1cm}
\includegraphics[width=1\textwidth, height=0.135\textheight, clip=true, angle=0]{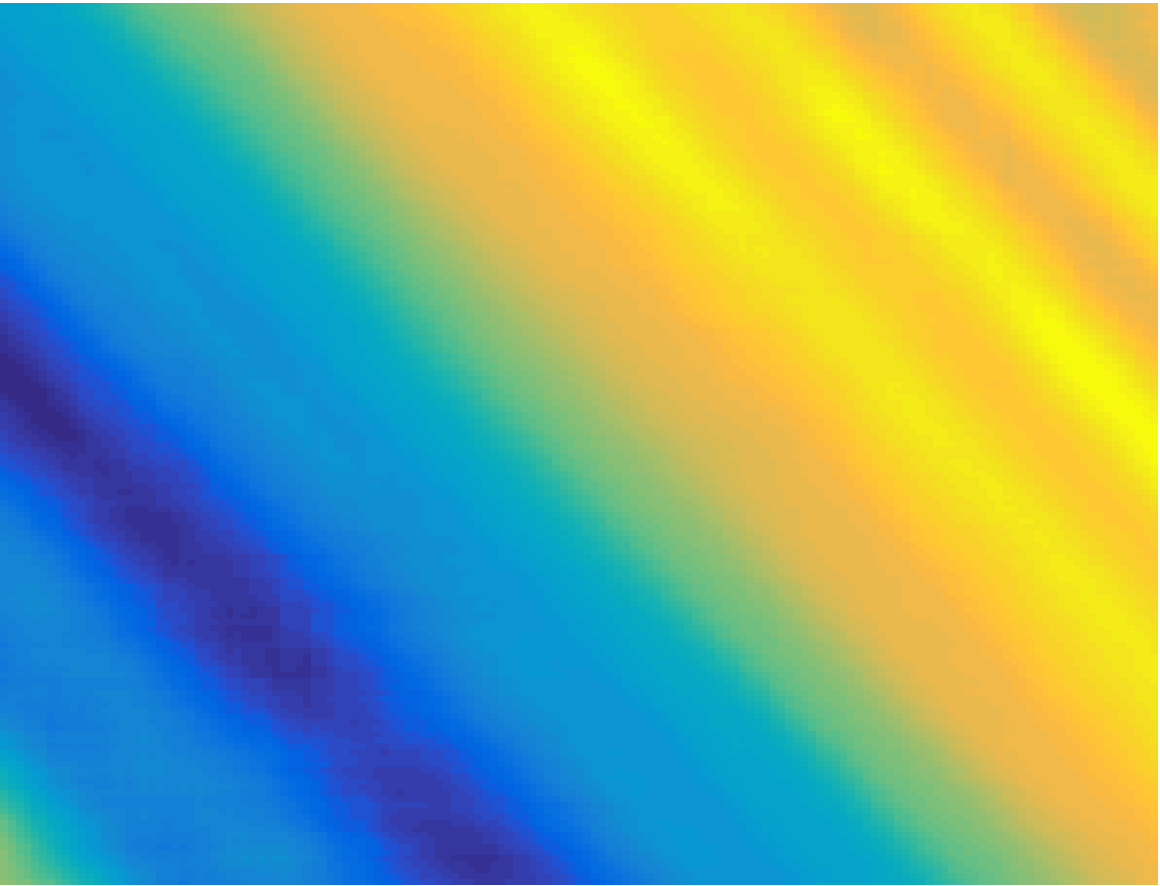}
Penalized recovery
\end{minipage}
\begin{minipage}{0.24\textwidth}
\centering
\vspace{-0.07cm}
\includegraphics[width=1\textwidth, height=0.135\textheight, clip=true, angle=0]{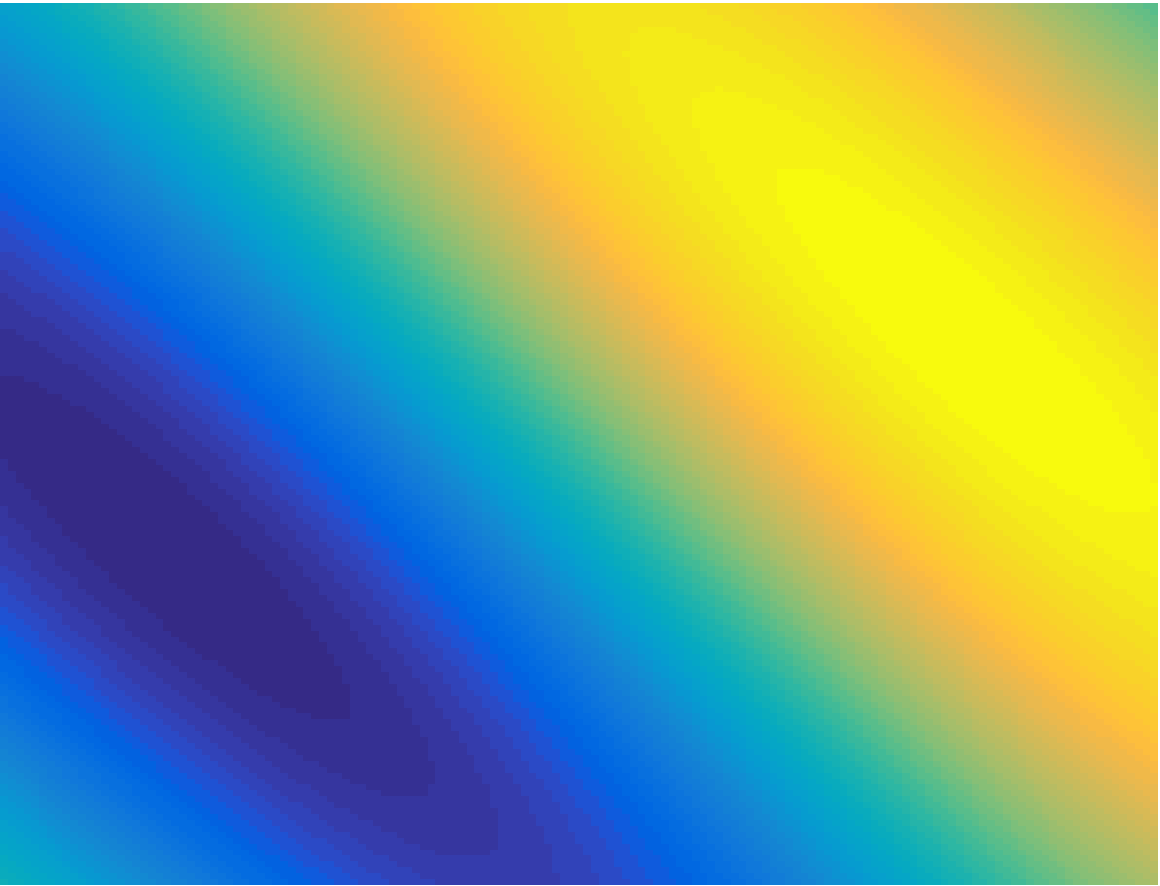}\vspace{0.1cm}
\includegraphics[width=1\textwidth, height=0.135\textheight, clip=true, angle=0]{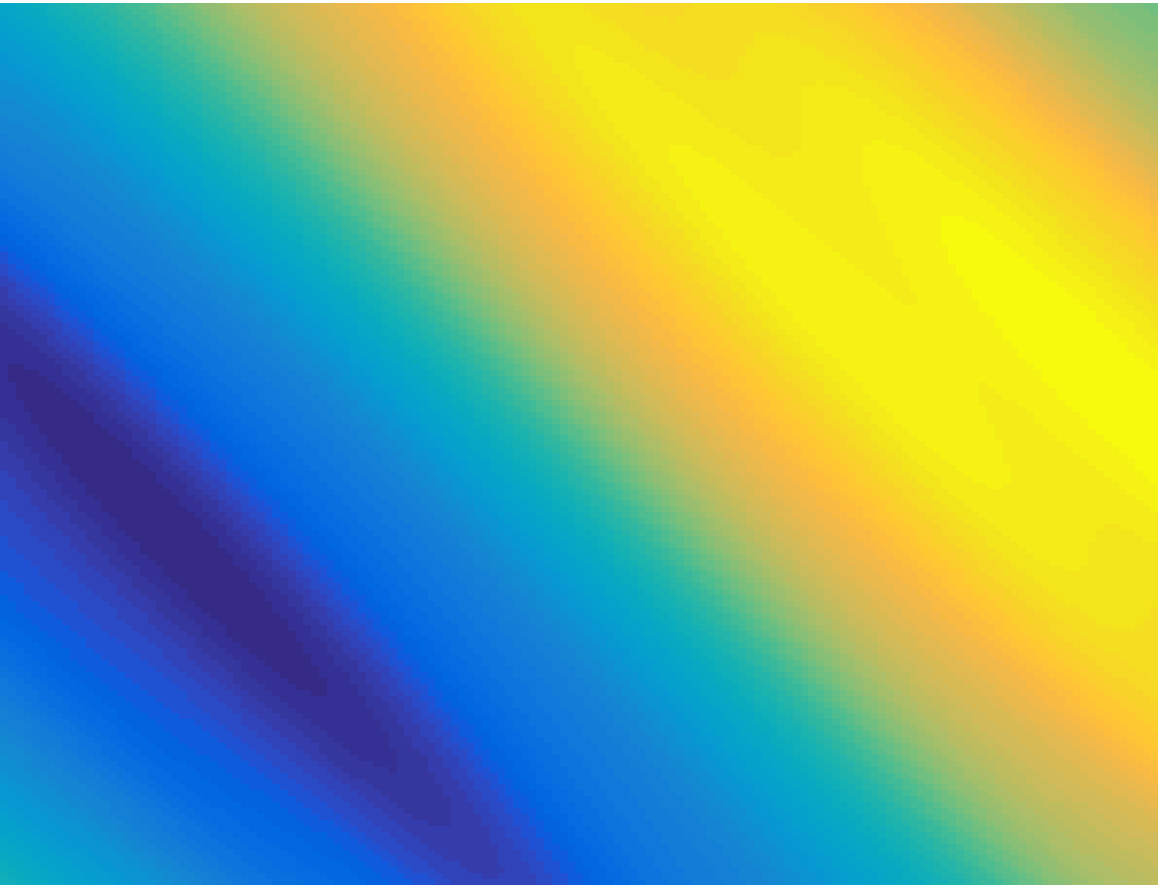}
Lasso
\end{minipage}
\caption{Recovery of a single-index signal  \eqref{eq:single-index}, observed with $\text{SNR}=1$, for $\beta=2$, first row, and $\beta=1$, second row. $\ell_2$-error ($\ell_2$-recovery vs. the Lasso): $0.07$ vs. $0.13$ for $\beta = 2$; $0.25$ vs. $0.31$ for $\beta = 1$. % Results are in Tab.~\ref{tab:single}.
}
\label{fig:single}
\end{figure}

\paragraph{Dimension reduction}
Here we illustrate denoising of a single-index signal \eqref{eq:single-index}, $\text{SNR} = 1$, with the direction $\theta$ close to the diagonal $(1,1)$, for two values of the smoothness index $\beta \in \{1,2\}$.  The results are presented in Fig.~\ref{fig:single}. One can see that Lasso tends to over-smooth the signal.

\begin{figure}[h!]
\center
\begin{minipage}{0.24\textwidth}
\centering
\includegraphics[width=1\textwidth, height=0.135\textheight, clip=true, angle=0]{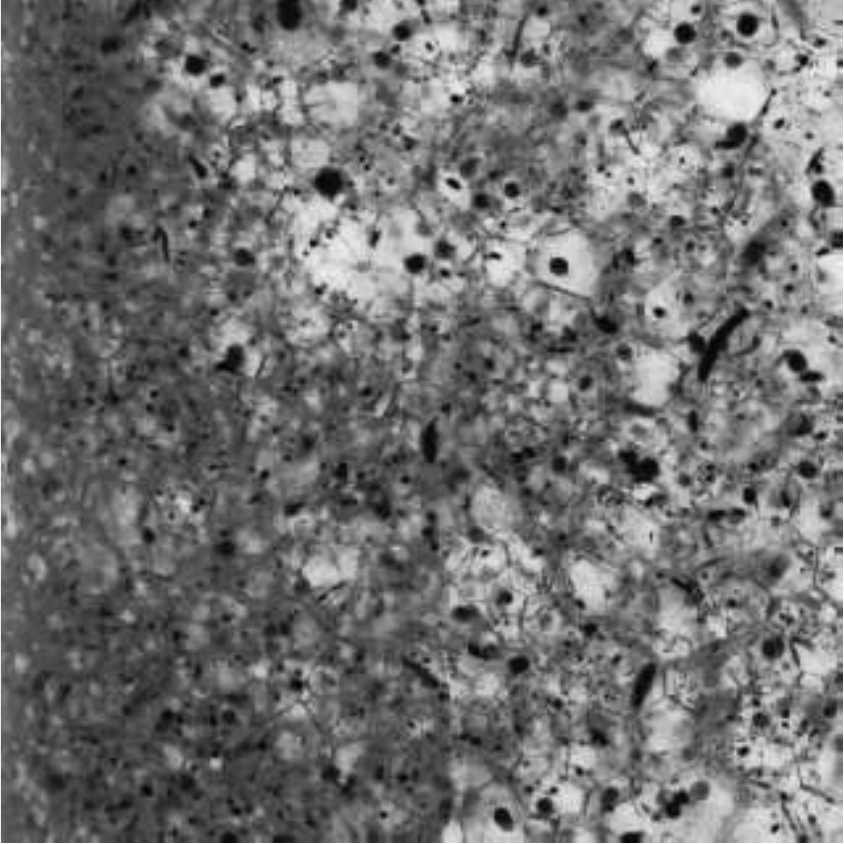}\vspace{0.1cm}
\includegraphics[width=1\textwidth, height=0.135\textheight, clip=true, angle=0]{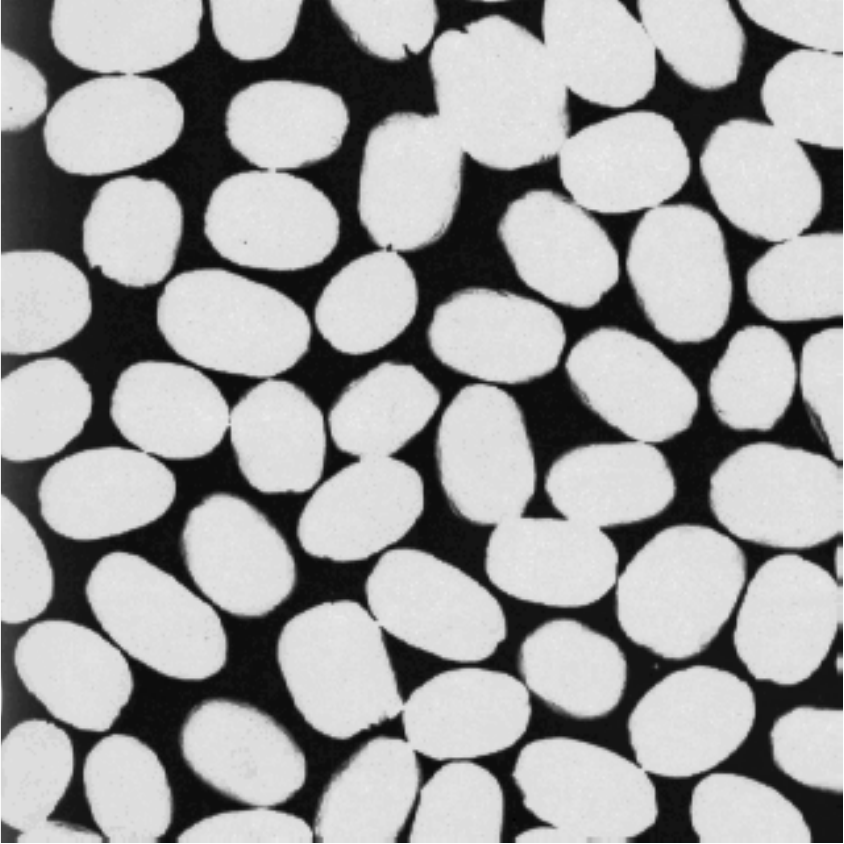}
True signal
\end{minipage}
\begin{minipage}{0.24\textwidth}
\centering
\vspace{-0.07cm}
\includegraphics[width=1\textwidth, height=0.135\textheight, clip=true, angle=0]{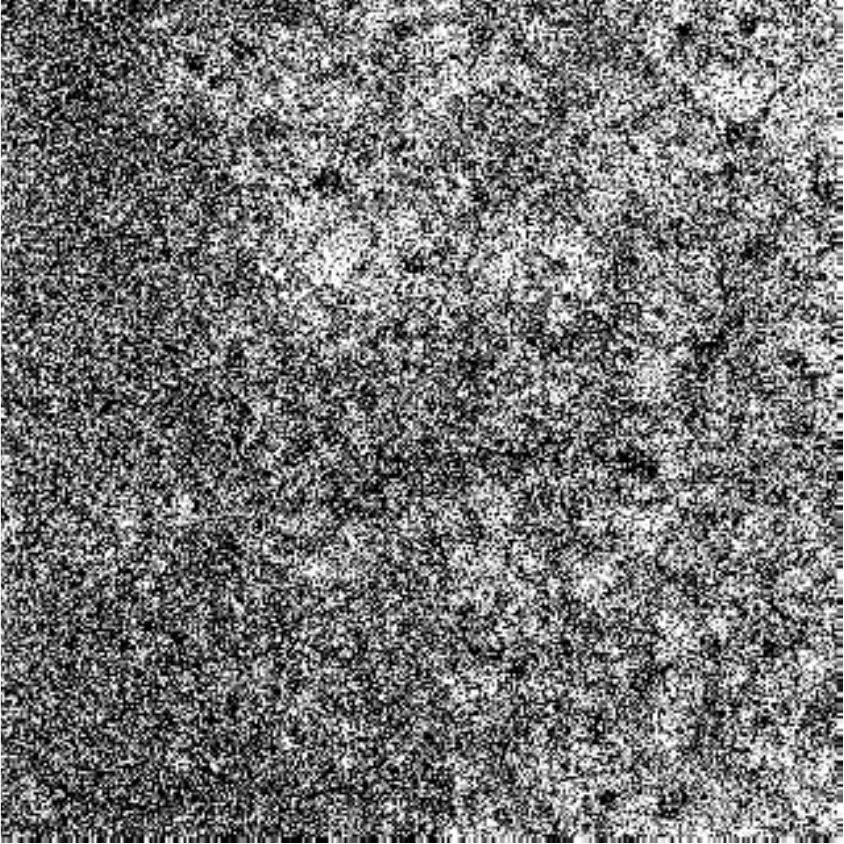}\vspace{0.1cm}
\includegraphics[width=1\textwidth, height=0.135\textheight, clip=true, angle=0]{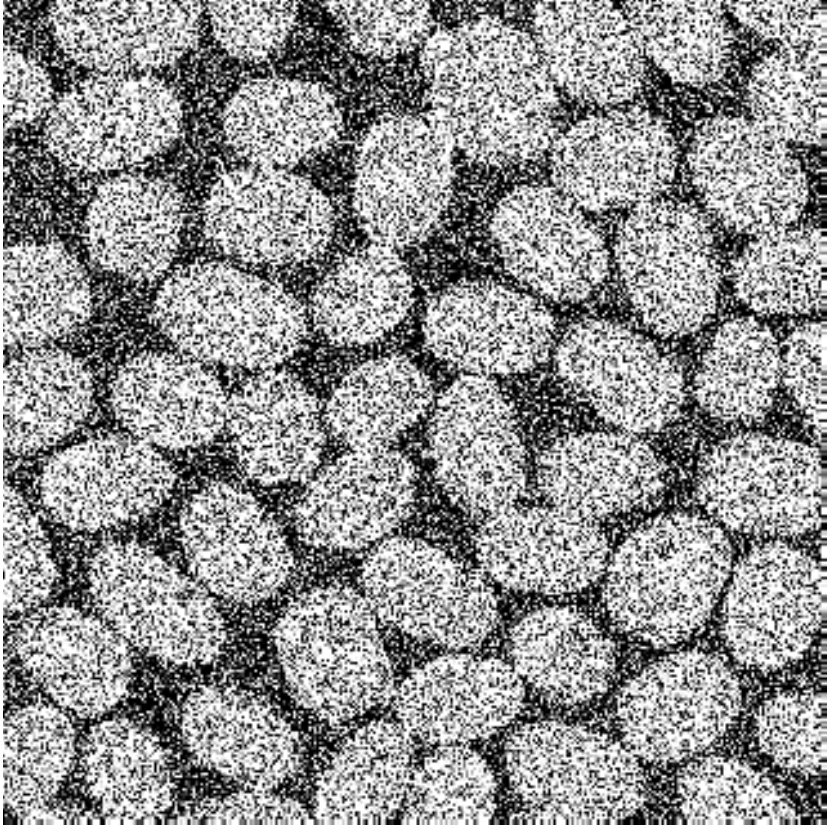}
Observations
\end{minipage}
\begin{minipage}{0.24\textwidth}
\centering
\includegraphics[width=1\textwidth, height=0.135\textheight, clip=true, angle=0]{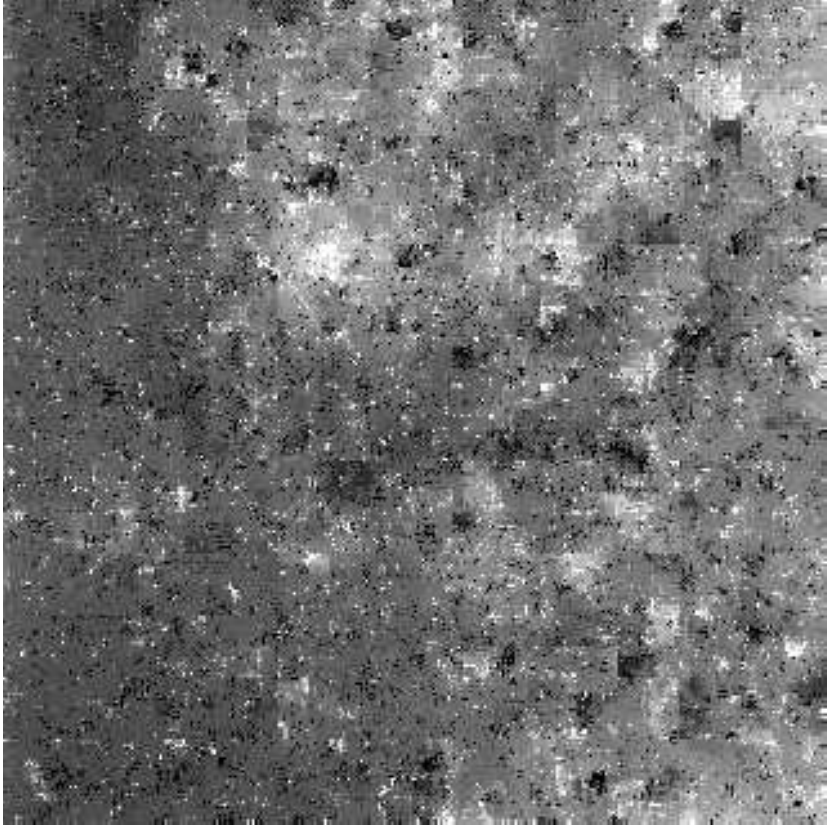}\vspace{0.1cm}
\includegraphics[width=1\textwidth, height=0.135\textheight, clip=true, angle=0]{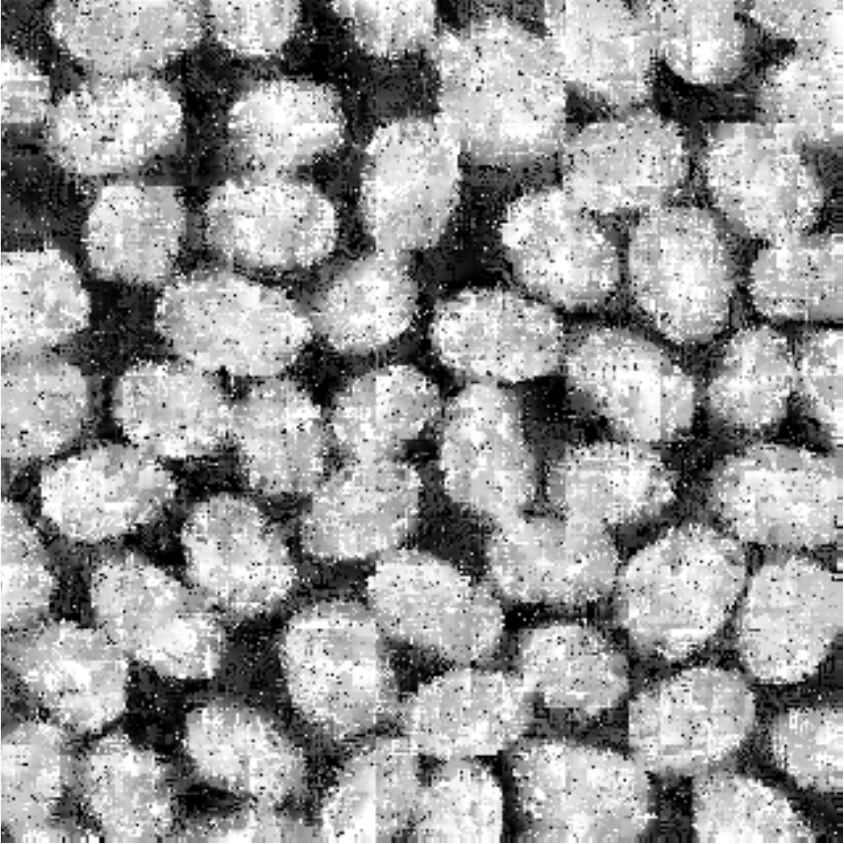}
Constrained recovery
\end{minipage}
\begin{minipage}{0.24\textwidth}
\centering
\vspace{-0.07cm}
\includegraphics[width=1\textwidth, height=0.135\textheight, clip=true, angle=0]{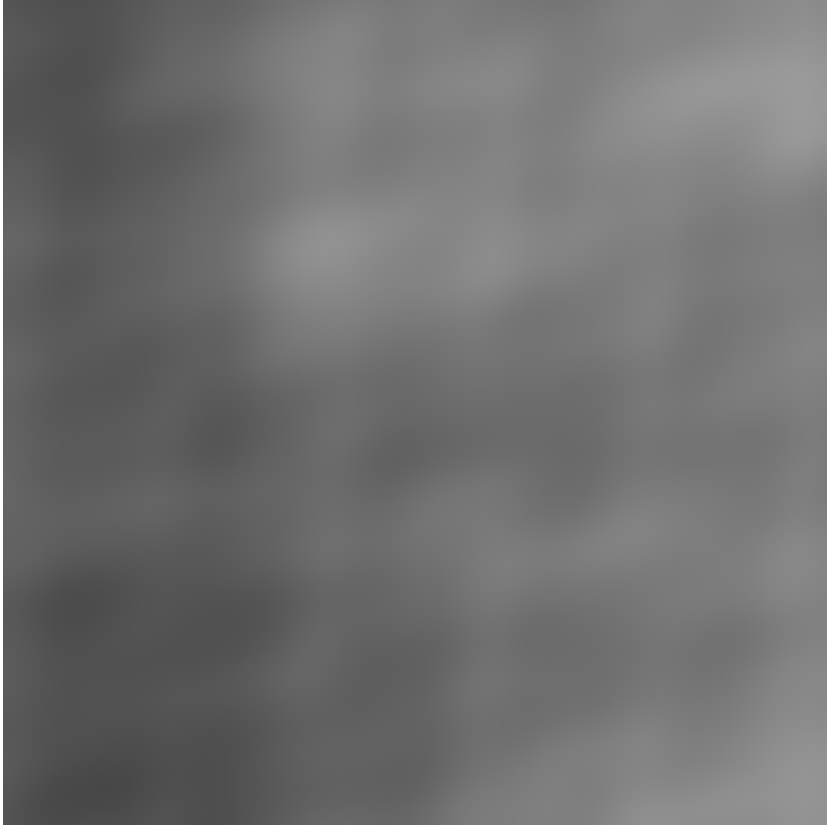}\vspace{0.1cm}
\includegraphics[width=1\textwidth, height=0.135\textheight, clip=true, angle=0]{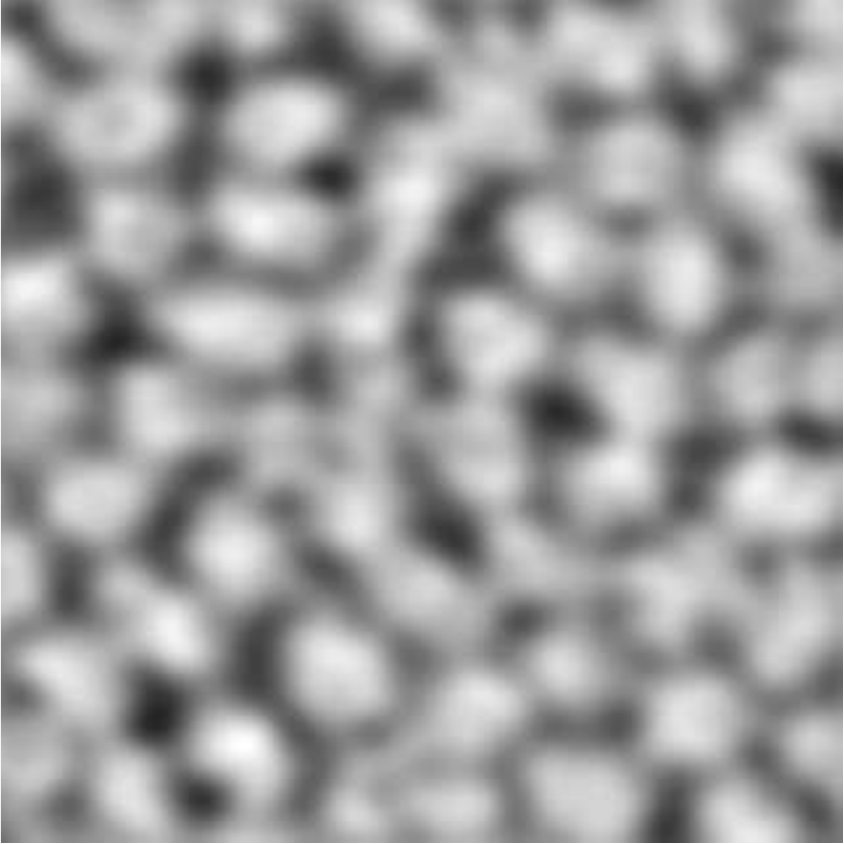}
Lasso
\end{minipage}
\caption{Recovery of two instances of the Original Brodatz database, cut by half to $320\times 320$ and observed with $\text{SNR}=1$. $\ell_2$-error:  $1.35e4$ for the constrained $\ell_2$-recovery, $1.25e4$ for the Lasso in the first row (inst. D73); $1.97e4$ for the constrained $\ell_2$-recovery,  $2.02e4$ for the Lasso method in the second row (inst. D75).
}
\label{fig:brodatz}
\end{figure}

\paragraph{Denoising textures}
In this experiment, we apply the proposed recovery methods to denoise two images from the original Brodatz texture database, observed in white Gaussian noise.
We set $\text{SNR} = 1$. We use the blockwise implementation of the constrained $\ell_2$-recovery algorithm, as described in Sec.~\ref{sec:expdetails}. We set the constraint parameter to $\overline{\varrho} = 4$, and we use the adaptation procedure
of~\cite[Sec. 3.2]{harchaoui2015adaptive} to define the estimation bandwidth. As in the above experiments, we use the Lasso~\cite{recht1} with $\lambda = \sigma\sqrt{2\log n}$, $n$ being the number of pixels. 
The resulting images are presented in Fig.~\ref{fig:brodatz}. Despite comparable quality in the mean square sense, the two methods significantly differ in their local behaviour. In particular, the constrained $\ell_2$-recovery better restores the local signal features, whereas the Lasso tends to oversmooth. 

\end{document}